\documentclass[a4paper]{article}
\usepackage{authblk}
\title{Polyhedral completeness of intermediate logics: the Nerve Criterion}
\author[1]{Sam Adam-Day}
\author[2]{Nick Bezhanishvili}
\author[3]{David Gabelaia}
\author[4]{Vincenzo Marra}
\affil[1]{\small Mathematical Institute, University of Oxford, United Kingdom}
\affil[2]{\small Institute for Logic, Language and Computation, University of Amsterdam, The Netherlands}
\affil[3]{\small A. Razmadze Mathematical Institute of I. Javakhishvili Tbilisi State University, Georgia}
\affil[4]{\small Dipartimento di Matematica ``Federigo Enriques'', Universit\`a degli Studi di Milano, Italy}
\date{\today}

\RequirePackage{silence}
\usepackage[utf8]{inputenc}
\usepackage[UKenglish]{babel}
\usepackage[margin=4cm]{geometry}
\usepackage[hidelinks]{hyperref}
\usepackage{amsmath}
\usepackage{mathtools}
\usepackage{amsthm}
\usepackage{thmtools}
\usepackage{thm-restate}
\usepackage{titling}
\usepackage{csquotes}
\usepackage[inline]{enumitem}
\usepackage{nicefrac}
\usepackage[capitalize,noabbrev]{cleveref}
\usepackage{fonttable}
\usepackage{tikz}
\usepackage{tikz-3dplot}
\usepackage[status=draft]{fixme}
\usepackage[backend=biber,style=alphabetic,maxbibnames=99]{biblatex}

\usepackage[charter]{mathdesign}
\usepackage[T1]{fontenc}
\usepackage{caption}
\usetikzlibrary{positioning,arrows,graphs,arrows.meta,decorations.pathreplacing,calligraphy}

\fxusetheme{color}
\FXRegisterAuthor{sam}{asam}{Sam}
\FXRegisterAuthor{nick}{anick}{Nick}
\FXRegisterAuthor{enzo}{aenzo}{Enzo}
\FXRegisterAuthor{david}{adavid}{David}

\newtheorem{theorem}{Theorem}[section]
\newtheorem{lemma}[theorem]{Lemma}
\newtheorem{proposition}[theorem]{Proposition}
\newtheorem{corollary}[theorem]{Corollary}
\theoremstyle{definition}
\newtheorem{definition}[theorem]{Definition}
\theoremstyle{remark}
\newtheorem{remark}[theorem]{Remark}

\tikzset{
	world/.style={draw,circle,outer sep=0pt,inner sep=0,minimum size=15},
	point/.style={draw=black,fill=black,opacity=1,circle,outer sep=0pt,inner sep=0,minimum size=2},
	facefill/.style={fill=blue!20,fill opacity=0.4}
}
\tikzgraphsset{
	poset/.style={grow up=#1, branch right=#1, empty nodes, simple, nodes=point},
	poset/.default=1cm
}

\makeatletter
\tikzdeclarecoordinatesystem{bc}
{%
  {%
    \pgf@xa=0pt
    \pgf@ya=0pt%
    \pgf@xb=0pt
    \tikz@bc@dolist#1,,%
    \pgfmathparse{1/\the\pgf@xb}%
    \global\pgf@x=\pgfmathresult\pgf@xa%
    \global\pgf@y=\pgfmathresult\pgf@ya%
  }%
}%

\def\tikz@bc@dolist#1,{%
  \def\tikz@temp{#1}%
  \ifx\tikz@temp\pgfutil@empty%
  \else
    \pgf@process{\pgfpointanchor{#1}{center}}%
    \pgfmathparse{1}%
    \advance\pgf@xa by\pgfmathresult\pgf@x%
    \advance\pgf@ya by\pgfmathresult\pgf@y%
    \advance\pgf@xb by\pgfmathresult pt%
    \expandafter\tikz@bc@dolist%
  \fi%
}
\makeatother

\newcommand{\starlike}[1]{\ensuremath{\ab{#1}}}

\newcommand{\Starlikes}{\ensuremath{\mathcal S}}
\newcommand{\jsl}[1]{\ensuremath{\chi(\starlike{#1})}}

\newcommand{\SL}{\ensuremath{\mathbf{SL}}}
\newcommand{\SFL}{\ensuremath{\mathbf{SFL}}}
\newcommand{\IPC}{\ensuremath{\mathbf{IPC}}}
\newcommand{\CPC}{\ensuremath{\mathbf{CPC}}}
\newcommand{\KC}{\ensuremath{\mathbf{KC}}}
\newcommand{\LC}{\ensuremath{\mathbf{LC}}}

\newcommand{\Sf}{\ensuremath{\mathbf{S4}}}

\newcommand{\SfGrz}{\ensuremath{\mathbf{S4.Grz}}}

\newcommand{\BD}{\ensuremath{\mathbf{BD}}}
\newcommand{\BW}{\ensuremath{\mathbf{BW}}}
\newcommand{\BTW}{\ensuremath{\mathbf{BTW}}}
\newcommand{\BC}{\ensuremath{\mathbf{BC}}}

\DeclareMathOperator{\depth}{\mathsf{depth}}

\DeclareMathOperator{\height}{\mathsf{height}}
\DeclareMathOperator{\width}{\mathsf{width}}

\DeclareMathOperator{\Trunk}{\mathsf{Trunk}}
\DeclareMathOperator{\Top}{\mathsf{Top}}

\DeclareMathOperator{\Up}{\mathrm{Up}}

\DeclareMathOperator{\Logic}{\mathrm{Logic}}
\DeclareMathOperator{\Frames}{\mathrm{Frames}}
\DeclareMathOperator{\FramesFin}{\Frames_\mathrm{fin}}
\DeclareMathOperator{\FramesRoot}{\Frames_\bot}
\DeclareMathOperator{\FramesFinRoot}{\Frames_{\bot,\mathrm{fin}}}

\DeclareMathOperator{\uset}{\ua}
\DeclareMathOperator{\dset}{\da}

\DeclareMathOperator{\Tree}{\mathrm{Tree}}
\DeclareMathOperator{\Succ}{\mathrm{Succ}}

\DeclareMathOperator{\Opens}{\mathcal{O}}
\DeclareMathOperator{\Closeds}{\mathcal{C}}
\DeclareMathOperator{\Int}{\mathrm{Int}}
\DeclareMathOperator{\Cl}{\mathrm{Cl}}
\newcommand{\comp}[1]{\ensuremath{#1^\mathsf{C}}}

\DeclareMathOperator{\Conv}{\mathrm{Conv}}
\DeclareMathOperator{\Relint}{\mathrm{Relint}}

\DeclareMathOperator{\Sd}{\mathrm{Sd}}
\newcommand{\R}{\ensuremath{\mathbb{R}}}
\newcommand{\Q}{\ensuremath{\mathbb{Q}}}

\newcommand{\Z}{\ensuremath{\mathbb{Z}}}
\DeclareMathOperator{\lcm}{\mathrm{lcm}}
\DeclareMathOperator{\Den}{\mathrm{Den}}

\newcommand{\Subo}{\ensuremath{\mathrm{Sub}_\mathrm{o}}}

\newcommand{\Sub}{\ensuremath{\mathrm{Sub}}}
\newcommand{\Po}{\ensuremath{\mathrm{P}_\mathrm{o}}}
\newcommand{\Pc}{\ensuremath{\mathrm{P}_\mathrm{c}}}

\newcommand{\Lo}{\ensuremath{\mathcal{L}}}
\newcommand{\CC}{\ensuremath{\mathbf{C}}}
\newcommand{\D}{\ensuremath{\mathbf{D}}}
\newcommand{\W}{\ensuremath{\mathcal{W}}}

\newcommand{\NN}{\mathbb N}
\newcommand{\NNp}{\mathbb N^{>0}}
\newcommand{\inv}[1]{{#1}^{-1}}
\newcommand{\invfs}[1]{\inv f\{#1\}}

\newcommand{\N}{\ensuremath{\mathcal{N}}}

\newcommand{\ostar}{\ensuremath{\mathrm{o}}}

\newcommand{\contype}{\ensuremath{\mathrm{ConType}}}
\newcommand{\concomps}{\ensuremath{\mathrm{ConComps}}}
\newcommand{\vDN}{\mathrel{{\vDash}_{\!\!\!\N}}}
\newcommand{\nvDN}{\mathrel{{\nvDash}_{\!\!\!\N}}}

\makeatletter
\newcommand{\@usstar}[1]{{\ua}\left(#1\right)}
\newcommand{\@usnostar}[1]{{\ua}(#1)}
\newcommand{\us}{\@ifstar{\@usstar}{\@usnostar}}
\newcommand{\@dsstar}[1]{{\da}\left(#1\right)}
\newcommand{\@dsnostar}[1]{{\da}(#1)}
\newcommand{\ds}{\@ifstar{\@dsstar}{\@dsnostar}}
\newcommand{\@udsstar}[1]{{\uda}\left(#1\right)}
\newcommand{\@udsnostar}[1]{{\uda}(#1)}
\newcommand{\uds}{\@ifstar{\@udsstar}{\@udsnostar}}
\newcommand{\@Usstar}[1]{{\Ua}\left(#1\right)}
\newcommand{\@Usnostar}[1]{{\Ua}(#1)}
\newcommand{\Us}{\@ifstar{\@Usstar}{\@Usnostar}}
\newcommand{\@Dsstar}[1]{{\Da}\left(#1\right)}
\newcommand{\@Dsnostar}[1]{{\Da}(#1)}
\newcommand{\Ds}{\@ifstar{\@Dsstar}{\@Dsnostar}}
\newcommand{\@Udsstar}[1]{{\Uda}\left(#1\right)}
\newcommand{\@Udsnostar}[1]{{\Uda}(#1)}
\newcommand{\Uds}{\@ifstar{\@Udsstar}{\@Udsnostar}}
\makeatother

\newcommand{\ep}{\epsilon}

\newcommand{\Lam}{\Lambda}
\newcommand{\Sig}{\Sigma}
\newcommand{\sig}{\sigma}
\newcommand{\vd}{\vdash}
\newcommand{\vD}{\vDash}

\newcommand{\nvd}{\nvdash}
\newcommand{\nvD}{\nvDash}

\newcommand{\bx}{\square}
\newcommand{\la}{\leftarrow}
\newcommand{\ra}{\rightarrow}

\newcommand{\Ra}{\Rightarrow}
\newcommand{\Lra}{\Leftrightarrow}

\newcommand{\rsa}{\rightsquigarrow}
\newcommand{\cra}{\circrightarrow}
\newcommand{\ua}{\uparrow}
\newcommand{\da}{\downarrow}
\newcommand{\uda}{\updownarrow}
\newcommand{\Ua}{\Uparrow}
\newcommand{\Da}{\Downarrow}
\newcommand{\Uda}{\Updownarrow}
\newcommand{\es}{\varnothing}
\newcommand{\sse}{\subseteq}

\newcommand{\wt}{\widetilde}
\newcommand{\wh}{\widehat}

\newcommand{\defeq}{\vcentcolon=}

\newcommand{\id}{\ensuremath{\mathsf{id}}}

\newcommand{\dom}{\ensuremath{\mathrm{dom}}}

\newcommand{\last}{\ensuremath{\mathsf{last}}}

\DeclarePairedDelimiter{\ab}{\langle}{\rangle}

\newcommand{\abs}[1]{\mathopen|#1\mathclose|}

\newcommand{\contradiction}{\noindent
	\begin{tikzpicture}[x=0.4ex,y=0.4ex]
		\draw[line width=.15ex] (0,0) -- (1,2) -- (0,2) -- (1,4)
		(0.95,0.32) -- (0,0) -- (-0.32,0.95);
	\end{tikzpicture}\hspace*{0.2em}
}
\newcommand{\circrightarrow}{\mathrel{{\circ}\mkern-3mu{\rightarrow}}}

\renewcommand{\leq}{\leqslant}
\renewcommand{\geq}{\geqslant}

\renewcommand{\preceq}{\preccurlyeq}

	{\newline\vspace{\abovedisplayskip}\hbox to \textwidth\bgroup\hss$\displaystyle}
	{$\hss\egroup\vspace{\belowdisplayskip}}

\newcommand{\FScott}{
	\raisebox{-1.1ex}{\resizebox{!}{3ex}{
		\begin{tikzpicture}
			\node[world] (z) at (0,0) {};
			\node[world] (a1) at (-0.7,1) {};
			\node[world] (a2) at (-0.7,2) {};
			\node[world] (b1) at (0.7,1) {};
			\path[draw] (z) -- (a1) -- (a2);
			\path[draw] (z) -- (b1);
		\end{tikzpicture}
	}}
}

\DeclareFontFamily{U} {MnSymbolC}{}

\DeclareFontShape{U}{MnSymbolC}{m}{n}{
	<-6> MnSymbolC5
	<6-7> MnSymbolC6
	<7-8> MnSymbolC7
	<8-9> MnSymbolC8
	<9-10> MnSymbolC9
	<10-12> MnSymbolC10
	<12-> MnSymbolC12}{}
\DeclareFontShape{U}{MnSymbolC}{b}{n}{
	<-6> MnSymbolC-Bold5
	<6-7> MnSymbolC-Bold6
	<7-8> MnSymbolC-Bold7
	<8-9> MnSymbolC-Bold8
	<9-10> MnSymbolC-Bold9
	<10-12> MnSymbolC-Bold10
	<12-> MnSymbolC-Bold12}{}

\DeclareSymbolFont{MnSyC} {U} {MnSymbolC}{m}{n}
\DeclareMathSymbol{\meddiamond}{\mathbin}{MnSyC}{110}

\usepackage{marginnote}
\makeatletter
\newcommand{\@rednotestar}[1]{\marginnote{\footnotesize\color{red}#1}}
\newcommand{\@rednotenostar}[1]{\textsuperscript{\color{red}†}\@rednotestar{#1}}
\newcommand{\rednote}{\@ifstar{\@rednotestar}{\@rednotenostar}}
\makeatother

\usepackage{xcolor}
\newcommand\bluecolour{black}
\newcommand\blue[1]{\textcolor{\bluecolour}{#1}}
\newenvironment{blueblock}{\color{\bluecolour}}{}

\begin{document}

	\maketitle

	\abstract{
		We investigate a recently-devised polyhedral semantics for intermediate logics, in which formulas are interpreted in $n$-dimensional polyhedra. An intermediate logic is \emph{polyhedrally complete} if it is complete with respect to some class of polyhedra. The first main result of this paper is a necessary and sufficient condition for the polyhedral-completeness of a logic. This condition, which we call the Nerve Criterion, is expressed in terms of Alexandrov's notion of the nerve of a poset. It affords a purely combinatorial characterisation of polyhedrally-complete logics.

		Using the Nerve Criterion we show, easily, that there are continuum many intermediate logics that are not polyhedrally-complete \blue{but which have the finite model property}. We also provide, at considerable combinatorial labour,  a countably infinite class of logics axiomatised by the Jankov-Fine formulas of `starlike trees' all of which are polyhedrally-complete. The polyhedral completeness theorem for these `starlike logics' is the second main result of this paper.}




\section{Introduction}

The genesis of many connections between logic and geometry is rooted in the discovery of topological semantics for intuitionistic and modal logic, as pioneered by Marshall Stone \cite{Stone1938}, Tang Tsao-Chen \cite{tsao-chen1938}, Alfred Tarski \cite{Tarski1939} and John C. C. McKinsey \cite{mckinsey1941}. This semantics is now well-known. In short, one starts with a topological space $X$, and interprets intuitionistic formulas inside the Heyting algebra of open sets of $X$, and modal formulas inside the modal algebra of subsets of $X$ with $\bx$ interpreted as the topological interior operator. A celebrated result due to Tarski \cite{Tarski1939} states that this provides a complete semantics for intuitionistic propositional logic ($\IPC$) on the one hand, and the modal logic \Sf\ on the other. Moreover, one can even obtain completeness with respect to certain individual spaces. Specifically, McKinsey and Tarski showed \cite{mckinseytarski44} that for any separable metric space $X$ without isolated points, if $\IPC \nvd \phi$, then $\phi$ has a countermodel based on $X$, and similarly with \Sf\ in place of \IPC. Later, Helena Rasiowa and Roman Sikorski showed that one can do without the assumption of separability \cite{rasiowasikorski1963}.

This result traces out an elegant interplay between topology and logic; however, it simultaneously establishes limits on the expressive power of this kind of interpretation. Indeed, examples of separable metric spaces without isolated points are the $n$-dimensional Euclidean space $\R^n$ and the Cantor space $2^\omega$. What McKinsey and Tarski's result shows, then, is that  these spaces have the same logic, namely \IPC\ (or \Sf). The upshot is that topological semantics does not allow logic to capture much of the geometric content of a space.

A natural idea is that, if we want to remedy the situation and allow for the capture of more information about a space, then we need an algebra finer than the Heyting algebra of open sets, or the modal algebra of arbitrary subsets with the interior operator. For example, Marco Aiello, Johan van Benthem, Guram Bezhanishvili and Mai Gehrke consider the modal logic of \emph{chequered} subsets of $\R^n$: finite unions of sets of the form $\prod_{i=1}^n C_i$, where each $C_i \sse \R$ is convex (\cite{AvBB03} and \cite{vBBG03}; see also \cite{vBB07}).  In \cite{tarski-polyhedra}, \cite{Gabelaia2017} and 
\cite{planar-polygons},  this algebra-refinement idea is taken one step further.  
To be able to capture some of the geometric content of a space, one may restrict attention to topological spaces and subsets which are \emph{polyhedra} (of arbitrary dimension). Indeed, the set $\Subo(P)$ of open subpolyhedra of $P$ is a Heyting algebra under $\sse$ (and a similar result holds in the modal case). 
This allows for an interpretation of intuitionistic and modal formulas in $\Subo(P)$. 
The main result of \cite{tarski-polyhedra} is that more is true. A polyhedral analogue of Tarski's theorem holds: these polyhedral semantics are complete for \IPC\ and \SfGrz. Furthermore, this approach delivers 
that logic can capture the dimension of the polyhedron in which it is interpreted, via the bounded depth formulas ${\sf bd}_n$ \cite[Sec.~2.4]{chagrovzakharyaschev1997}. In particular, the polyhedron $P$ is $n$ dimensional if, and only if, $P$ validates ${\sf bd}_{n+1}$ and does not validate ${\sf bd}_{n+2}$ for $n\in \omega$ \cite{tarski-polyhedra}. 

In this paper we make further advances in the study of polyhedral semantics. We introduce and study polyhedral completeness for intermediate logics. We say that  an intermediate logic $L$ is polyhedrally complete if there is a class $\mathcal{C}$
of polyhedra such that $L$ is the logic of $\mathcal{C}$. It follows from  \cite{tarski-polyhedra} that $\IPC$ and the logic $\BD_n$ of bounded depth $n$, for each $n$, are  
polyhedrally-complete. We construct 
infinitely many polyhedrally-complete logics, and  show that there are continuum many polyhedrally incomplete ones \blue{all of which have the finite model property}.

To this end we employ a time-honoured tool from combinatorial and polyhedral geometry, the \emph{nerve} of a poset (=partially ordered set).  The nerve will be our key concept  relating logic with polyhedral geometry. In detail, the nerve $\N(F)$ of the poset $F$ is the collection of finite non-empty chains in $F$ ordered by inclusion. As was already noted in \cite{tarski-polyhedra}, given a polyhedron $P$,  a triangulation of $P$ corresponds to  a validity-preserving map from $P$ onto the poset $F$ of the faces of the triangulation. Through Esakia duality, in turn, this validity-preserving map corresponds to an embedding of the Heyting algebra of upsets of $F$ into the Heyting algebra of open subpolyhedra of $P$. Nerves are closely related to \emph{barycentric subdivisions} of triangulations. Indeed, if a finite poset $F$ is the face poset of some triangulation $\Sigma$ of a polyhedron $P$, then $\N(F)$ corresponds to a barycentric subdivision of $\Sigma$. 
	
	Applying methods and results from rational polyhedral geometry we present a proof of our first main result, the \emph{Nerve Criterion for polyhedral completeness} (\cref{thm:nerve criterion}): A logic $L$ is complete with respect to some class of polyhedra if and only if it is the logic of a class of finite posets closed under taking nerves. Thus, we obtain that the logic of any given polyhedron is the logic of the iterated nerves of any one of its triangulations. The criterion yields many negative results, showing in particular that there are continuum-many non-polyhedrally-complete logics with the finite model property (\cref{thm:stable logics poly incomplete}).

As to positive results, we consider logics defined using \emph{starlike trees} --- trees which only branch at the root --- as forbidden configurations. \emph{Starlike logics} are then those defined by the Jankov-Fine formulas of a collection of starlike trees. Exploiting the Nerve Criterion, and a result by Zakharyaschev \cite{zakharyaschev93} that all these logics have the finite model property, we prove our second main result (\cref{thm:starlike completeness}): Every starlike  logic is polyhedrally-complete. This yields a countably infinite class of polyhedrally-complete logics of each finite height and of infinite height. (For instance, Scott's well-known logic \SL\ is in this class.) As forbidden configurations, starlike trees have a natural geometric meaning, expressing connectedness properties of polyhedral spaces.

The paper is organised as follows. In \cref{sec:preliminaries}, we give the required background on intermediate logics and polyhedral geometry. \Cref{sec:suboP} presents the polyhedral semantics first defined in \cite{tarski-polyhedra}. In \cref{sec:nerve criterion}, we present and prove the Nerve Criterion for polyhedral-completeness (\cref{thm:nerve criterion}). Making use of this criterion, \cref{sec:stable logics} establishes that all stable logics (as defined in \cite{bezhbezh2009}) of height at least $2$ are polyhedrally-incomplete. Then in \cref{sec:starlike completeness}, we define the class of starlike logics, and prove that each one is polyhedrally-complete. The techniques \blue{in these} two sections are entirely combinatorial.

Finally, let us briefly comment on further research. One major  problem already mentioned in \cite{tarski-polyhedra} is to characterise the logic of piecewise-linear manifolds of a fixed dimension. Here we announce significant progress on this question; the results will appear in a forthcoming paper. A second relevant goal would be a complete classification of polyhedrally-complete logics. At the the time of writing, we do not know how to attain this goal. One might wonder if our results on starlike logics extend to arbitrary trees, or even to a wider class of posets. As to the latter, some negative results are obtained in  \cite[Corollary~4.12]{sam-thesis}. For the former, the situation is rather obscure to us at the time of writing; cf.\ the discussion on `general trees' in \cite[p.~61]{sam-thesis}. Identifying further classes of polyhedrally complete logics beyond the starlike ones introduced in this paper would be the next immediate task in the direction of obtaining a classification  of polyhedrally-complete logics.\footnote{\blue{This paper is partly based on the fist-named author's M.Sc. thesis \cite{sam-thesis}.}}

\section{Preliminaries}
\label{sec:preliminaries}
 In this section we remind the reader of the relational and algebraic semantics for intermediate logics, and survey the definitions and results which will play their part in the forthcoming. As a main reference we use \cite{chagrovzakharyaschev1997}. We assume rather less familiarity with polyhedral geometry, and thus present in more detail the material we need.

\subsection{Posets as Kripke frames}

	A \emph{\textup{(}Kripke\textup{)} frame} for intuitionistic logic is simply a poset. We thus use the term `frame' in this paper as a synonym of `poset'. The validity relation $\vD$ between frames and formulas is defined in the usual way, see, e.g., \cite[Ch.~2]{chagrovzakharyaschev1997}. Given a class of frames $\CC$, its \emph{logic} is:
	\begin{equation*}
		\Logic(\CC) \defeq \{\phi \text{ a formula } \mid \forall F \in \CC \colon F \vD \phi \}
	\end{equation*}
	Conversely, given a logic $\Lo$, define:
	\begin{gather*}
		\Frames(\Lo) \defeq \{F\text{ a Kripke frame} \mid F \vD \Lo\} \\
		\FramesFin(\Lo) \defeq \{F\text{ a finite Kripke frame} \mid F \vD \Lo\}
	\end{gather*}
	A logic $\Lo$ has the \emph{finite model property} (fmp) if it is the logic of a class of finite frames. Equivalently, if $\Lo = \Logic(\FramesFin(\Lo))$.

Fix a poset $F$. For any $x \in F$, its \emph{upset}, \emph{downset}, \emph{strict upset} and \emph{strict downset} are defined, respectively, as follows.
	\begin{gather*}
		\us x \defeq \{y \in F \mid y \geq x\} \\
		\ds x \defeq \{y \in F \mid y \leq x\} \\
		\Us x \defeq \{y \in F \mid y > x\} \\
		\Ds x \defeq \{y \in F \mid y < x\}
	\end{gather*}
	For any set $S \sse F$, its \emph{upset} and \emph{downset} are defined, respectively, as follows.
	\begin{gather*}
		\uset U \defeq \bigcup_{x \in U} \us x\\
		\dset U \defeq \bigcup_{x \in U} \ds x
	\end{gather*}
	A \emph{subframe} is a subposet. A subframe $U \sse F$ is \emph{upwards-closed} or a \emph{generated subframe} if $U = \uset U$, and it is \emph{downwards-closed} if $\dset U = U$. The \emph{Alexandrov topology} on $F$ is the set $\Up F$ of its upwards-closed subsets. This constitutes a topology on $F$. In the sequel, we will freely switch between thinking of $F$ as a poset and as a topological space. Note that the closed sets in this topology correspond to downwards-closed sets.

	A \emph{chain} in $F$ is $X \sse F$ which as a subposet is linearly-ordered. The \emph{length} of the chain $X$ is $\abs X$. A chain $X \sse F$ is \emph{maximal} if there is no chain $Y \sse F$ such that $X \subset Y$ (i.e. such that $X$ is a proper subset of $Y$). The \emph{height} of $F$ is the element of $\NN \cup \{\infty\}$ defined by:
	\begin{equation*}
		\height(F) \defeq \sup\{\abs X-1 \mid X \sse F\text{ is a chain}\}
	\end{equation*}
	For notational uniformity, say that this value is also the \emph{depth} of $F$, $\depth(F)$. For any $x \in F$, define its \emph{height} and \emph{depth} as follows.
	\begin{gather*}
		\height(x) \defeq \height(\ds x) \\
		\depth(x) \defeq \depth(\us x)
	\end{gather*}
	The \emph{height} of a logic $\Lo$ is the element of $\NN \cup \{\infty\}$ given by:
	\begin{equation*}
		\height(\Lo) \defeq \sup\{\height(F) \mid F \in \Frames(\Lo)\}
	\end{equation*}
	A \emph{top element} of $F$ is $t \in F$ such that $\depth(t)=0$. For any $x,y \in F$, say that $x$ is an \emph{immediate predecessor} of $y$, and that $y$ is an \emph{immediate successor} of $x$, if $x < y$ and there is no $z \in F$ such that $x < z < y$. Write $\Succ(x)$ for the collection of immediate successors of $x$.

	The poset $F$ is \emph{rooted} if it has a minimum element, which is called the \emph{root}, and is usually denoted by $\bot$. Define:
	\begin{gather*}
		\FramesRoot(\Lo) \defeq \{F \in \Frames(\Lo) \mid F\text{ is rooted}\} \\
		\FramesFinRoot(\Lo) \defeq \{F \in \FramesFin(\Lo) \mid F\text{ is rooted}\}
	\end{gather*}

	An \emph{antichain} in $F$ is a subset $Z \sse F$ in which no two elements are comparable. The \emph{width}, notation $\width(F)$, of $F$ is the cardinality of the largest antichain in $F$. 

	A function $f \colon F \to G$ is a \emph{p-morphism} if for every $x \in F$ we have:
	\begin{equation*}
		f(\us x) = \us{f(x)}
	\end{equation*}
	Equivalently, $f$ should satisfy the following conditions.
	\begin{gather*}
		\forall x,y \in F \colon (x \leq y \Ra f(x) \leq f(y)) \tag{Forth}\label{p:p-morphism forth} \\
		\forall x \in F \colon \forall z \in G \colon (f(x) \leq z \Ra \exists y \in F \colon (x \leq y \wedge f(y)=z)) \tag{Back}\label{p:p-morphism back}
	\end{gather*}
	An \emph{up-reduction} from $F$ to $G$ is a surjective p-morphism $f$ from an upwards-closed set $U \sse F$ to $G$. Write $f \colon F \cra G$.

	\begin{proposition}\label{prop:up-reduction logic containment}
		If there is an up-reduction $F \cra G$ then $\Logic(F) \sse \Logic(G)$. In other words, if $G \nvD \phi$ then $F \nvD \phi$.
	\end{proposition}

	\begin{proof}
		See \cite[Corollary~2.8, p.~30 and Corollary~2.17, p.~32]{chagrovzakharyaschev1997}.
	\end{proof}

	\begin{corollary}\label{cor:logic of frames logic of rooted frames}
		If $\CC$ is any collection of frames and $\Lo = \Logic(\CC)$, then:
		\begin{equation*}
			\Lo = \Logic(\FramesRoot(\Lo))
		\end{equation*}
	\end{corollary}

	\begin{proof}
		First, $\Lo \sse \Logic(\FramesRoot(\Lo))$. Conversely, suppose $\Lo \nvd \phi$. Then there exists $F \in \CC$ such that $F \nvD \phi$, hence there is $x \in F$ such that $x \nvD \phi$ (for some valuation on $F$), meaning that $\us x \nvD \phi$. Now, $\us x$ is upwards-closed in $F$, hence $\id_{\us x}$ is an up-reduction $F \cra \us x$. Then by \cref{prop:up-reduction logic containment}, we get that $\us x \vD \Lo$, so that $\us x \in \FramesRoot(\Lo)$.
	\end{proof}

\subsection{Heyting algebras, topological semantics}

	A \emph{Heyting algebra} is a bounded lattice equipped with a \emph{Heyting implication} $\ra$ that satisfies:
	\begin{equation*}
		c \leq a \ra b \quad\Lra\quad c \wedge a \leq b
	\end{equation*}
	The validity relation $\vD$ between Heyting algebras and formulas is defined in the usual way; the  notation $\Logic(-)$ is extended appropriately. The logic of a Heyting algebra is exactly the logic of its finitely generated subalgebras. Say that $A$ is \emph{locally-finite} if for every $S \sse A$ finite, the algebra $\ab S$ generated by $S$ is finite. If $F$ is any poset, the bounded distributive lattice $\Up F$ is a Heyting algebra, and:
	\begin{proposition}\label{prop:finite esakia duality logic-preserving}
		If $F$ is \blue{a poset}, $\Logic(F) = \Logic(\Up F)$
		\end{proposition}

	\begin{proof}
		See \cite[Corollary~8.5, p.~238]{chagrovzakharyaschev1997}.
	\end{proof}

	Co-Heyting algebras are the order-duals of Heyting algebras. Specifically, a \emph{co-Heyting algebra}  is a bounded lattice equipped with a \emph{co-Heyting implication} $\la$ that satisfies:
	\begin{equation*}
		a \la b \leq c \quad\Lra\quad a \leq b \vee c
	\end{equation*}
	For more information on co-Heyting algebras, the reader is referred to \cite[\S1]{mckinseytarski1946} and \cite{rauszer1974}, where they are called `Brouwerian algebras'.

	Given a topological space $X$, we regard the collection of open sets $\Opens(X)$ of $X$ as a Heyting algebra in the standard manner, cf.\ \cite[Proposition~8.31, p.~247]{chagrovzakharyaschev1997}. (Recall that
	\begin{equation*}
		U \ra V = \Int(\comp U \cup V)
	\end{equation*}
	where $\Int$ denotes the topological interior operator, and $\comp -$ is set-theoretic complement.) We can thus  interpret formulas in topological spaces. Write $X \vD \phi$ for $\Opens(X) \vD \phi$, and extend the remaining notations accordingly. 
	
%
%

	The topological space $X$ also comes with a co-Heyting algebra, namely its collection of closed sets $\Closeds(X)$. The co-Heyting implication on $\Closeds(X)$ satisfies:
	\begin{equation*}
		C \la D \defeq \Cl(C \setminus D)
	\end{equation*}
	where $\Cl$ denotes the topological closure operator. If a Heyting algebra $A$  is regarded as a poset category $(A, \leq)$, then its opposite category $A^\mathrm{op} = (A, \geq)$ is a co-Heyting algebra. In the case of the Heyting algebra $\Opens(X)$ of open sets of $X$, $\Opens(X)^\mathrm{op}$ is isomorphic to the   co-Heyting algebra $\Closeds(X)$ of closed subsets of $X$.

		
\subsection{Jankov-Fine formulas as forbidden configurations}

	To every finite, rooted frame $Q$, we associate a formula $\chi(Q)$, the \emph{Jankov-Fine} formula of $Q$ (also called its \emph{Jankov-De Jongh formula}). The precise definition of $\chi(Q)$ is somewhat involved, but the exact details of this syntactical form are not relevant for our considerations. What matters to us is its notable semantic property.

	\begin{theorem}\label{thm:Jankov-Fine up-reductions}
		For any frame $F$, we have that $F \vD \chi(Q)$ if and only if $F$ does not up-reduce to $Q$.
	\end{theorem}

	\begin{proof}
		See \cite[\S9.4, p.~310]{chagrovzakharyaschev1997} for a treatment in which Jankov-Fine formulas are considered as specific instances of more general `canonical formulas'. An alternative proof can be found in \cite[\S3.3, p.~56]{bezhanishvili2006}, which gives a complete definition of $\chi(Q)$. See also \cite{bezhbezh2009} for an algebraic version of this result.
	\end{proof}

	Jankov-Fine formulas formalise the intuition of `forbidden configurations'. The formula $\chi(Q)$ `forbids' the configuration $Q$ from its frames.

	The following consequence of \cref{thm:Jankov-Fine up-reductions} will come in handy later on.

	\begin{corollary}\label{cor:Logic C every finite rooted frame up-reduction}
		Let $\Lo = \Logic(\CC)$ where $\CC$ is a class of frames. Then:
		\begin{equation*}
			\FramesFinRoot(\Lo) = \{F \text{ finite, rooted frame} \mid \exists G \in \CC \colon G \cra F\}
		\end{equation*}
	\end{corollary}

	\begin{proof}
		First, if $F$ is a finite, rooted frame such that there is $G \in \CC$ and an up-reduction $G \cra F$, then by \cref{prop:up-reduction logic containment} we have that $F \in \FramesFinRoot(\Lo)$. Conversely take $F$ finite and rooted, and assume that there is no $G \in \CC$ with $G \cra F$. Then by \cref{thm:Jankov-Fine up-reductions}, $G \vD \chi(F)$ for every $G \in \CC$; whence $\Lo \vd \chi(F)$. By \cref{thm:Jankov-Fine up-reductions}, $F \nvD \chi(F)$ implying $F \nvD \Lo$. This yields $F \notin \FramesFinRoot(\Lo)$.
	\end{proof}

\subsection{Intermediate logics}

	 The logic $\IPC$ is  intuitionistic propositional logic. An \emph{intermediate logic} is any consistent logic extending $\IPC$. Classical logic, $\CPC$, is the largest intermediate logic.

	\begin{proposition}\label{prop:ipc Kripke complete and has fmp}
		\IPC\ is the logic of the class of all finite frames, i.e.\ \IPC\   has the fmp.
	\end{proposition}

	\begin{proof}
		See \cite[Theorem~2.57, p.~49]{chagrovzakharyaschev1997}.
	\end{proof}

	For every $n \in \NN$, let $\BD_n$ be the logic of all finite frames of height at most $n$. This has the following axiomatisation in terms of Jankov-Fine formulas. 

	\begin{proposition}\label{prop:BDn iff no p-morphism to chain n}
		$\BD_n$ is the logic axiomatised by $\IPC$ plus the Jankov-Fine formula of the chain (linear order) on $n+1$ elements.
	\end{proposition}

	\begin{proof}
		See \cite[Table~9.7, p.~317, and \S9]{chagrovzakharyaschev1997}.
	\end{proof}

	Scott's Logic, $\SL$, is usually axiomatised by the Scott sentence:
	\begin{equation*}
		\SL = \IPC + \IPC + ((\neg\neg p \ra p) \ra p \vee \neg p) \ra \neg p \vee \neg\neg p
	\end{equation*}
	This logic can also be axiomatised using a forbidden configuration, as follows.

	\begin{proposition}
		$\SL = \IPC + \chi(\FScott)$.
	\end{proposition}

	\begin{proof}
		See \cite[Table~9.7, p.~317, and \S9]{chagrovzakharyaschev1997}.
	\end{proof}

\subsection{Polytopes, polyhedra, and simplices}

	Polyhedra are certain subsets of finite-dimensional real affine spaces $\R^n$. An \emph{affine combination} of $x_0, \ldots, x_d \in \R^n$ is a point $r_0 x_0 + \cdots + r_d x_d$, where  $r_0, \ldots, r_d \in \R$ are such that $r_0 + \cdots + r_d = 1$. A \emph{convex combination} is an affine combination in which additionally each $r_i \geq 0$. Given a set $S \sse \R^n$, its \emph{convex hull}, notation $\Conv S$, is the collection of convex combinations of its elements. (We stress that each convex combination involves, by definition, a finite subset of $S$ only.) A subspace $S \sse \R^n$ is \emph{convex} if $\Conv S = S$. A \emph{polytope} is the convex hull of a finite set. A \emph{polyhedron} in $\R^n$ is a set which can be expressed as the finite union of polytopes. A \emph{subpolyhedron} of a polyhedron $P$ in $\R^n$ is a subset of $P$ which is itself a polyhedron. Note that every polyhedron is closed and bounded, hence compact, in the canonical (Euclidean) topology carried by the real affine space $\R^n$. All topological notions pertaining to polyhedra in the following refer to this topology.

	A set of points $x_0, \ldots, x_d$ is \emph{affinely independent} if whenever:
	\begin{equation*}
		r_0 x_0 + \cdots + r_d x_d = \mathbf 0 \quad\text{and}\quad r_0 + \cdots + r_d = 0
	\end{equation*}
	we must have that $r_0, \ldots, r_d = 0$. This is equivalent to saying that the vectors:
	\begin{equation*}
		x_1 - x_0, \ldots, x_d - x_0
	\end{equation*}
	are linearly independent. A  \emph{$d$-simplex} is the convex hull $\sig$ of $d+1$ affinely independent points $x_0, \ldots, x_d$, which we call its \emph{vertices}. Write $\sig = x_0\cdots x_d$; the \emph{dimension} of $\sigma$ is $d$.

	\begin{proposition}
		Every simplex determines its vertex set: two simplices coincide if and only if they share the same vertex set.
	\end{proposition}

	\begin{proof}
		See \cite[Proposition 2.3.3, p.~32]{maunder1980algebraic}.
	\end{proof}

	\noindent A \emph{face} of $\sig$ is the convex hull $\tau$ of some non-empty subset of $\{x_0, \ldots, x_d\}$ (note that $\tau$ is then a simplex too). Write $\tau \preceq \sig$, and $\tau \prec \sig$ if $\tau \neq \sig$.

	Since $x_0, \ldots, x_d$ are affinely independent, every point $x \in \sig$ can be expressed uniquely as a convex combination $x = r_0 x_0 + \cdots + r_d x_d$ with $r_0, \ldots, r_d \geq 0$ and $r_0 + \cdots + r_d = 1$. Call the tuple $(r_0, \ldots, r_d)$ the \emph{barycentric coordinates} of $x$ in $\sig$. The \emph{barycentre} $\wh\sig$ of $\sig$ is the special point whose barycentric coordinates are $(\frac{1}{d+1}, \ldots, \frac{1}{d+1})$. The \emph{relative interior} of $\sig$ is defined:
	\begin{equation*}
		\Relint \sig \defeq \{r_0 x_0 + \cdots + r_d x_d \in \sig \mid r_0, \ldots, r_d >0\}
	\end{equation*}
	The relative interior of $\sig$ is `$\sig$ without its boundary' in the following sense. The \emph{affine subspace spanned by $\sig$} is the set of all affine combinations of $x_0, \ldots, x_d$. Then the relative interior of $\sig$ coincides with the topological interior of $\sig$ inside this affine subspace, the latter being equipped with the subspace topology  it inherits from $\R^n$. Note that $\Cl\Relint\sig = \sig$, the closure being taken in the ambient space $\R^n$.

\subsection{Triangulations}

	A \emph{simplicial complex} in $\R^n$ is a finite set $\Sig$ of simplices satisfying the following conditions.
	\begin{enumerate}[label=(\alph*)]
		\item \label{item:downwards-closed; defn:simplicial complex}
		$\Sig$ is $\prec$-downwards-closed: whenever $\sig \in \Sig$ and $\tau \prec \sig$ we have $\tau \in \Sig$.
		\item \label{item:intersection; defn:simplicial complex}
		If $\sig,\tau \in \Sig$, then $\sig \cap \tau$ is either empty or a common face of $\sig$ and $\tau$.
	\end{enumerate}
	The \emph{support} of $\Sig$ is the set $\abs\Sig \defeq \bigcup \Sig$. Note that by definition this set is automatically a polyhedron. We say that $\Sig$ is a \emph{triangulation} of the polyhedron $\abs\Sig$. The set $\Sig$ is a poset under $\prec$, called the \emph{face poset} of the triangulation. A \emph{subcomplex} of $\Sig$ is subset which is itself a simplicial complex. Note that a subcomplex, as a poset, is precisely a downwards-closed set. Given $\sig \in \Sig$, its \emph{open star} is defined:
	\phantomsection\label{defn:open star}
	\begin{equation*}
		\ostar(\sig) \defeq \bigcup \{\Relint(\tau) \mid \tau \in \Sig \text{ and }\sig \sse \tau\}
	\end{equation*}

	\begin{proposition}\label{prop:relints partition support}
		The relative interiors of the simplices in a simplicial complex $\Sig$ partition $\abs\Sig$. That is, for every $x \in \abs\Sig$, there is exactly one $\sig \in \Sig$ such that $x \in \Relint\sig$.
	\end{proposition}

	\begin{proof}
		See \cite[Proposition~2.3.6, p.~33]{maunder1980algebraic}.
	\end{proof}

	In light of \cref{prop:relints partition support}, for any $x \in \abs\Sig$ let us write $\sig^x$ for the unique $\sig \in \Sig$ such that $x \in \Relint\sig$. The simplex $\sig^x$ is known as the \emph{carrier} of $x$.

	\begin{proposition}\label{prop:simplex relint excludes proper faces}
		Let $\Sig$ be a simplicial complex, take $\tau \in \Sig$ and $x \in \Relint\tau$. Then no proper face $\sig \prec \tau$ contains $x$. This means that $\sig^x$ is the inclusion-smallest simplex containing $x$.
	\end{proposition}

	\begin{proof}
		See \cite[Lemma~3.1]{tarski-polyhedra}.
	\end{proof}

	The next result is a basic fact of polyhedral geometry, and is of fundamental importance in its connection with logic. For $\Sig$ a triangulation and $S$ a subset of the ambient space $\R^n$, define:
	\begin{equation*}
		\Sig_S \defeq \{\sig \in \Sig \mid \sig \sse S\}
	\end{equation*}
	This, being a downwards-closed subset of $\Sig$, is a subcomplex of $\Sig$.

	\begin{lemma}[Triangulation Lemma]\label{lem:triangulation lemma}
		Any polyhedron admits a triangulation which simultaneously triangulates each of any fixed finite set of subpolyhedra. That is, for a collection of polyhedra $P, Q_1, \ldots, Q_m$ such that each $Q_i \sse P$, there is a triangulation $\Sig$ of $P$ such that $\Sig_{Q_i}$ triangulates $Q_i$ for each $i$.
	\end{lemma}

	\begin{proof}
		See \cite[Theorem~2.11 and Addendum~2.12, p.~16]{rourkesanderson1972}.
	\end{proof}

\subsection{Stellar and barycentric subdivisions}\label{ss:stellar}
For  $\Sig$ and $\Delta$  simplicial complexes, $\Delta$ is a \emph{subdivision} or \emph{refinement} of $\Sig$, notation $\Delta \lhd \Sig$, if $\abs\Sig=\abs\Delta$ and every simplex of $\Delta$ is contained in a simplex of $\Sig$.

	\begin{lemma}\label{lem:subdivision simplex union of simplicies}
		If $\Delta \lhd \Sig$ then for every $\sig \in \Sig$ we have:
		\begin{equation*}
			\sig = \bigcup\{\tau \in \Delta \mid \tau \sse \sig\}
		\end{equation*}	
	\end{lemma}

\begin{figure}
		\begin{equation*}
			\begin{tikzpicture}[scale=0.45, draw=black]
				\begin{scope}[yshift=100,every node/.style=point]
					\coordinate (x) at (90:2);
					\coordinate (y) at (210:2);
					\coordinate (z) at (330:2);
					\filldraw[draw=black,facefill] (x) -- (y) -- (z) -- cycle;
	 				\node at (x) {};
	 				\node at (y) {};
					\node at (z) {};
					\draw[-Latex] (2.5,0.5) -- (3.5,0.5);
				\end{scope}
				\begin{scope}[xshift=170,yshift=100,every node/.style=point]
					\coordinate (x) at (90:2);
					\coordinate (y) at (210:2);
					\coordinate (z) at (330:2);
					\coordinate (xy) at (bc cs:x,y);
					\filldraw[draw=black,facefill] (x) -- (y) -- (z) -- cycle;
	 				\node at (x) {};
	 				\node at (y) {};
					\node at (z) {};
					\node at (xy) {};
					\draw (z) -- (xy);
					\draw[-Latex] (2.5,0.5) -- (3.5,0.5);
				\end{scope}
				\begin{scope}[xshift=340,yshift=100,every node/.style=point]
					\coordinate (x) at (90:2);
					\coordinate (y) at (210:2);
					\coordinate (z) at (330:2);
					\coordinate (xy) at (bc cs:x,y);
					\coordinate (xyz) at (bc cs:xy,z);
					\filldraw[draw=black,facefill] (x) -- (y) -- (z) -- cycle;
	 				\node at (x) {};
	 				\node at (y) {};
					\node at (z) {};
					\node at (xy) {};
					\node at (xyz) {};
					\draw (z) -- (xy);
					\draw (x) -- (xyz) -- (y);
				\end{scope}
				\begin{scope}[every node/.style=point]
					\fill[facefill] (0,0) -- (2,1) -- (3,-1) -- (2.6,-2.4) -- (0,-2.6) -- (0,-1) -- cycle;
					\draw (0,0) node {} -- (2,1) node {} -- (3,-1) node {} -- (2.6,-2.4) node {} -- (0,-2.6) node {} -- (0,-1.3) node {} -- cycle;
					\draw (3,-1) -- (5,0.1) node {} -- (3.8,1.2) node {};
					\draw (3,-1) -- (0,0) -- (2.6,-2.4) -- (0,-1.3);
					\draw[-Latex] (5.5,-1) -- (6.5,-1);
				\end{scope}
				\begin{scope}[xshift=230, every node/.style=point]
					\fill[facefill] (0,0) -- (2,1) -- (3,-1) -- (2.6,-2.4) -- (0,-2.6) -- (0,-1) -- cycle;
					\draw (0,0) node (x) {} -- (2,1) node {} -- (3,-1) node (z) {} -- (2.6,-2.4) node (y) {} -- (0,-2.6) node {} -- (0,-1.3) node (w) {} -- cycle;
					\draw (3,-1) -- (5,0.1) node {} -- (3.8,1.2) node {};
					\draw (3,-1) -- (0,0) -- (2.6,-2.4) -- (0,-1.3);
					\node (xy) at (bc cs:x,y) {};
					\draw (z) -- (xy) -- (w);
				\end{scope}
				\tdplotsetmaincoords{47}{95}
				\begin{scope}[xshift=15, yshift=-130, tdplot_main_coords, scale=2]
					\coordinate (x0) at (0,0,0);
					\coordinate (x1) at (2,0,0);
					\coordinate (x2) at (1,1.732,0);
					\coordinate (x3) at (1,0.577,1.633);
					\fill[facefill] (x0) -- (x3) -- (x2);
					\fill[facefill] (x0) -- (x1) -- (x2);
					\draw (x0) -- (x2);
					\filldraw[facefill] (x1) -- (x0) -- (x3);
					\filldraw[facefill] (x1) -- (x2) -- (x3) -- cycle;
					\node[point] at (0,0,0) {};
					\node[point] at (2,0,0) {};
					\node[point] at (1,1.732,0) {};
					\node[point] at (1,0.577,1.633) {};
					\draw[scale=0.5,tdplot_screen_coords,-Latex] (5,-1) -- (6,-1);
				\end{scope}
				\begin{scope}[xshift=250, yshift=-130, tdplot_main_coords, scale=2]
					\coordinate (x0) at (0,0,0);
					\coordinate (x1) at (2,0,0);
					\coordinate (x2) at (1,1.732,0);
					\coordinate (x3) at (1,0.577,1.633);
					\coordinate (x12) at (bc cs:x1,x2);
					\fill[facefill] (x0) -- (x3) -- (x2);
					\fill[facefill] (x0) -- (x12) -- (x3);
					\fill[facefill] (x0) -- (x1) -- (x2);
					\draw (x0) -- (x2);
					\draw (x0) -- (x12);
					\filldraw[facefill] (x1) -- (x0) -- (x3);
					\filldraw[facefill] (x1) -- (x2) -- (x3) -- cycle;
					\draw (x3) -- (x12);
					\node[point] at (0,0,0) {};
					\node[point] at (2,0,0) {};
					\node[point] at (1,1.732,0) {};
					\node[point] at (1,0.577,1.633) {};
					\node[point] at (x12) {};
				\end{scope}
			\end{tikzpicture}
		\end{equation*}
	 	\caption{Examples of elementary stellar subdivisions}
	 	\label{fig:elementary subdivision examples}
	 \end{figure}
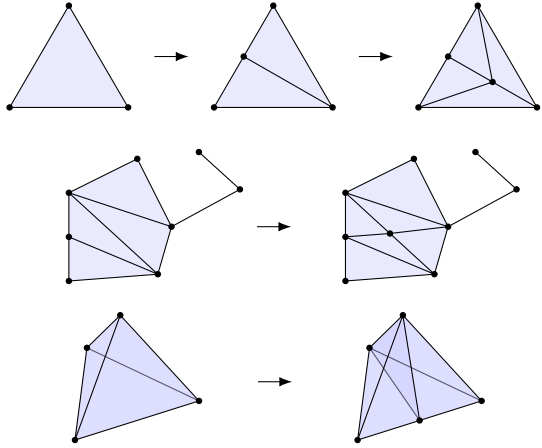

	\begin{proof}
		Let $S \defeq \{\tau \in \Delta \mid \tau \sse \sig\}$. Clearly $\bigcup S \sse \sig$. Conversely, for $x \in \sig$, let $\tau^x \in \Delta$ be such that $x \in \Relint\tau^x$. Since $\Delta$ refines $\Sig$, there is some $\rho \in \Sig$ such that $\tau^x \sse \rho$; assume that $\rho$ is inclusion-minimal with this property. It follows from \cite[\S3, Lemma~3, p.~121]{spanier1966} that $\Relint\tau^x \sse \Relint\rho$, meaning that $x \in \sig \cap \Relint \rho$. By condition \ref{item:intersection; defn:simplicial complex} in the definition of a simplicial complex, we have that $\sig \cap \rho$ is face of $\rho$. But then by \cref{prop:simplex relint excludes proper faces}, $\rho \preceq \sig$, since otherwise $\sig \cap \rho$ would be a proper face of $\rho$ containing $x \in \Relint \rho$. Therefore $\tau^x \sse \rho \sse \sig$ so that $x \in \bigcup S$.
	\end{proof}
We now introduce a special class of subdivisions, for which  the original source \cite{alexander1930} remains a fundamental reference. Let $\Sigma$ be a simplicial complex, and let $c\in |\Sigma|$. The \emph{elementary stellar subdivision of $\Sigma$ at $c$} is the set of simplices $\Delta$ obtained from $\Sigma$ via the following transformation: Replace each simplex $\sigma\in\Sigma$ that contains $c$ by the set of all simplices $\Conv{\{\tau\cup \{c\}\}}$, where $\tau$ ranges over all faces of $\sigma$ that do not contain $c$. It can then be proved that $\Delta$ is again a triangulation, and a subdivision of $\Sigma$. The equality $\Sigma=\Delta$ holds precisely when  the chosen $c$ is a vertex of $\Sigma$. If $\Delta$ is a subdivision of $\Sigma$ that is obtained via a finite number of successive elementary subdivisions of $\Sigma$, then $\Delta$ is a \emph{stellar subdivision} of $\Sigma$.  See \cref{fig:elementary subdivision examples}.

	If $\Delta$ is obtained from $\Sigma$ via an elementary stellar subdivision at $c\in |\Sigma|$, and $c$ is moreover the barycentre of the vertices of its carrier simplex $\sigma^c\in \Sigma$ (see \cref{prop:relints partition support} and the comments following it), then $\Delta$ is an \emph{elementary barycentric subdivision of $\Sigma$} (at the barycentre $c$).

	\begin{figure}
		\begin{equation*}
			\begin{tikzpicture}[scale=0.45, draw=black]
				\begin{scope}[yshift=100,every node/.style=point]
					\coordinate (x) at (90:2);
					\coordinate (y) at (210:2);
					\coordinate (z) at (330:2);
					\filldraw[draw=black,facefill] (x) -- (y) -- (z) -- cycle;
	 				\node at (x) {};
	 				\node at (y) {};
					\node at (z) {};
					\draw[-Latex] (2.5,0.5) -- (3.5,0.5);
				\end{scope}
				\begin{scope}[xshift=170,yshift=100,every node/.style=point]
					\coordinate (x) at (90:2);
					\coordinate (y) at (210:2);
					\coordinate (z) at (330:2);
					\filldraw[draw=black,facefill] (x) -- (y) -- (z) -- cycle;
	 				\node at (x) {};
	 				\node at (y) {};
					\node at (z) {};
					\node (xy) at (bc cs:x,y) {};
					\node (xz) at (bc cs:x,z) {};
					\node (yz) at (bc cs:y,z) {};
					\node (xyz) at (bc cs:x,y,z) {};
					\graph[use existing nodes]{
						xyz -- {x,y,z,xy,xz,yz};
					};
					\draw[-Latex] (2.5,0.5) -- (3.5,0.5);
				\end{scope}
				\begin{scope}[xshift=340,yshift=100,every node/.style=point]
					\coordinate (x) at (90:2);
					\coordinate (y) at (210:2);
					\coordinate (z) at (330:2);
					\filldraw[draw=black,facefill] (x) -- (y) -- (z) -- cycle;
	 				\node at (x) {};
	 				\node at (y) {};
					\node at (z) {};
					\node (xy) at (bc cs:x,y) {};
					\node (xz) at (bc cs:x,z) {};
					\node (yz) at (bc cs:y,z) {};
					\node (xyz) at (bc cs:x,y,z) {};
					\node (xxy) at (bc cs:x,xy) {};
					\node (xxz) at (bc cs:x,xz) {};
					\node (yxy) at (bc cs:y,xy) {};
					\node (yyz) at (bc cs:y,yz) {};
					\node (zxz) at (bc cs:z,xz) {};
					\node (zyz) at (bc cs:z,yz) {};
					\node (xxyz) at (bc cs:x,xyz) {};
					\node (yxyz) at (bc cs:y,xyz) {};
					\node (zxyz) at (bc cs:z,xyz) {};
					\node (xyxyz) at (bc cs:xy,xyz) {};
					\node (xzxyz) at (bc cs:xz,xyz) {};
					\node (yzxyz) at (bc cs:yz,xyz) {};
					\node (xxyxyz) at (bc cs:x,xy,xyz) {};
					\node (xxzxyz) at (bc cs:x,xz,xyz) {};
					\node (yxyxyz) at (bc cs:y,xy,xyz) {};
					\node (yyzxyz) at (bc cs:y,yz,xyz) {};
					\node (zxzxyz) at (bc cs:z,xz,xyz) {};
					\node (zyzxyz) at (bc cs:z,yz,xyz) {};
					\graph[use existing nodes]{
						xyz -- {x,y,z,xy,xz,yz};
						xxyxyz -- {x,xy,xyz,xxy,xxyz,xyxyz};
						xxzxyz -- {x,xz,xyz,xxz,xxyz,xzxyz};
						yxyxyz -- {y,xy,xyz,yxy,yxyz,xyxyz};
						yyzxyz -- {y,yz,xyz,yyz,yxyz,yzxyz};
						zxzxyz -- {z,xz,xyz,zxz,zxyz,xzxyz};
						zyzxyz -- {z,yz,xyz,zyz,zxyz,yzxyz};
					};
				\end{scope}
				\begin{scope}[every node/.style=point]
					\fill[facefill] (0,0) -- (2,1) -- (3,-1) -- (2.6,-2.4) -- (0,-2.6) -- (0,-1) -- cycle;
					\draw (0,0) node {} -- (2,1) node {} -- (3,-1) node {} -- (2.6,-2.4) node {} -- (0,-2.6) node {} -- (0,-1.3) node {} -- cycle;
					\draw (3,-1) -- (5,0.1) node {} -- (3.8,1.2) node {};
					\draw (3,-1) -- (0,0) -- (2.6,-2.4) -- (0,-1.3);
					\draw[-Latex] (5.5,-1) -- (6.5,-1);
				\end{scope}
				\begin{scope}[xshift=230, every node/.style=point]
					\fill[facefill] (0,0) -- (2,1) -- (3,-1) -- (2.6,-2.4) -- (0,-2.6) -- (0,-1.3) -- cycle;
					\draw (0,0) node (a) {} -- (2,1) node (b) {} -- (3,-1) node (c) {} -- (2.6,-2.4) node (d) {} -- (0,-2.6) node (e) {} -- (0,-1.3) node (f) {} -- cycle;
					\draw (c) -- (5,0.1) node (g) {} -- (3.8,1.2) node (h) {};
					\draw (c) -- (0,0) -- (2.6,-2.4) -- (0,-1.3);
					\node at (bc cs:a,b) (ab) {};
					\node at (bc cs:a,c) (ac) {};
					\node at (bc cs:b,c) (bc) {};
					\node at (bc cs:a,c) (ac) {};
					\node at (bc cs:a,d) (ad) {};
					\node at (bc cs:c,d) (cd) {};
					\node at (bc cs:a,d) (ad) {};
					\node at (bc cs:a,f) (af) {};
					\node at (bc cs:d,f) (df) {};
					\node at (bc cs:d,e) (de) {};
					\node at (bc cs:d,f) (df) {};
					\node at (bc cs:e,f) (ef) {};
					\node at (bc cs:a,b,c) (abc) {};
					\node at (bc cs:a,c,d) (acd) {};
					\node at (bc cs:a,d,f) (adf) {};
					\node at (bc cs:d,e,f) (def) {};
					\node at (bc cs:c,g) (cg) {};
					\node at (bc cs:g,h) (gh) {};
					\draw (abc)
						edge (a)
						edge (b)
						edge (c)
						edge (ab)
						edge (ac)
						edge (bc);
					\draw (acd)
						edge (a)
						edge (c)
						edge (d)
						edge (ac)
						edge (ad)
						edge (cd);
					\draw (adf)
						edge (a)
						edge (d)
						edge (f)
						edge (ad)
						edge (af)
						edge (df);
					\draw (def)
						edge (d)
						edge (e)
						edge (f)
						edge (de)
						edge (df)
						edge (ef);
				\end{scope}
				\tdplotsetmaincoords{47}{95}
				\begin{scope}[xshift=15, yshift=-130, tdplot_main_coords, scale=2]
					\coordinate (x0) at (0,0,0);
					\coordinate (x1) at (2,0,0);
					\coordinate (x2) at (1,1.732,0);
					\coordinate (x3) at (1,0.577,1.633);
					\fill[facefill] (x0) -- (x3) -- (x2);
					\fill[facefill] (x0) -- (x1) -- (x2);
					\draw (x0) -- (x2);
					\filldraw[facefill] (x1) -- (x0) -- (x3);
					\filldraw[facefill] (x1) -- (x2) -- (x3) -- cycle;
					\node[point] at (0,0,0) {};
					\node[point] at (2,0,0) {};
					\node[point] at (1,1.732,0) {};
					\node[point] at (1,0.577,1.633) {};
					\draw[scale=0.5,tdplot_screen_coords,-Latex] (5,-1) -- (6,-1);
				\end{scope}
				\begin{scope}[xshift=250, yshift=-130, tdplot_main_coords, scale=2]
					\node[point] (x0) at (0,0,0) {};
					\node[point] (x1) at (2,0,0) {};
					\node[point] (x2) at (1,1.732,0) {};
					\node[point] (x3) at (1,0.577,1.633) {};
					\coordinate (x01) at (1,0,0);
					\coordinate (x12) at (1.5,0.866,0)	;
					\coordinate (x23) at (1,1.155,0.816) ;
					\coordinate (x03) at (0.5,0.289,0.816) ;
					\coordinate (x02) at (0.5,0.866,0);
					\coordinate (x13) at (1.5,0.289,0.816) ;
					\coordinate	(x012) at (1,0.577,0);
					\coordinate (x013) at (1,0.192,0.544);
					\coordinate (x023) at (0.667,0.770,0.544);
					\coordinate (x123) at (1.333,0.770,0.544);
					\node[point] (x0123) at (1,0.577,0.544) {};
					\graph[use existing nodes] {
						x01 -- {x0, x1};
						x02 -- {x0, x2};
						x03 -- {x0, x3};
						x12 -- {x1, x2};
						x13 -- {x1, x3};
						x23 -- {x2, x3};
						x012 -- {x0, x1, x2, x01, x02, x12};
						x013 -- {x0, x1, x3, x01, x03, x13};
						x023 -- {x0, x2, x3, x02, x03, x23};
						x123 -- {x1, x2, x3, x12, x13, x23};
						x0123 -- {x0, x1, x2, x3, x01, x02, x03, x12, x13, x23, x012, x013, x023, x123}
					};
				\end{scope}
			\end{tikzpicture}
		\end{equation*}
		\caption[Examples of barycentric subdivision]{Examples of barycentric subdivision. \blue{Each simplex in the simplicial complex is divided at its barycentre, proceeding in decreasing  order of dimension. The bottom right tetrahedron is drawn without filled-in faces to aid clarity.}}
		\label{fig:barycentric subdivision examples}
	\end{figure}
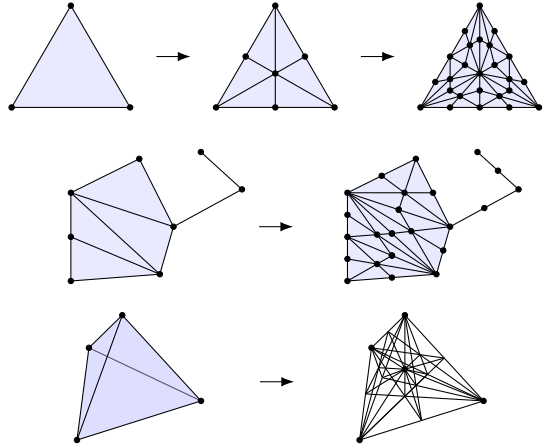

	The \emph{barycentric subdivision} $\Sd\Sig$ of $\Sig$ is then defined as the refinement of $\Sigma$ obtained by successively applying   elementary barycentric subdivisions at each simplex of $\Sigma$, proceeding in decreasing  order of dimension. It can be proved that $\Sd\Sig$ does not depend on the chosen ordering of the simplices of $\Sigma$. See the examples in \cref{fig:barycentric subdivision examples}. In the literature,   $\Sd\Sig$ is also often called the \emph{first derived} subdivision of $\Sigma$; cf. e.g. \cite{rourkesanderson1972}. We inductively define, for each $k\in\mathbb{N}$, the \emph{$k^{\rm th}$ derived subdivision of $\Sig$}: $\Sig^{(0)}\coloneqq\Sig$; and $\Sig^{(k)} =\Sd \Sig^{(k-1)}$.


\section{The algebra of open subpolyhedra}
\label{sec:suboP}

	With the preliminaries in place, we relate  intuitionistic logic and polyhedra. For further details please see \cite{tarski-polyhedra}.

\subsection{Polyhedral semantics} 

	Given a polyhedron $P$, let $\Sub P$ denote the set of its subpolyhedra.

	\begin{theorem}\label{thm:subP locally-finite co-Heyting algebra}
		$\Sub P$ is a co-Heyting algebra, and a subalgebra of \blue{$\Closeds(P)$}.
	\end{theorem}

	\begin{proof}
		See \blue{\cite[Corollary~3.4]{tarski-polyhedra}}. The proof makes fundamental use of the Triangulation Lemma.
	\end{proof}

	By an \emph{open subpolyhedron} of a polyhedron $P$ in this paper we mean the complement (in $P$) of a  subpolyhedron of $P$. Denote by $\Subo P$ the collection of open subpolyhedra in $P$. Evidently, there is an isomorphism $\Subo  P\cong(\Sub P)^{\rm op}$, and  \cref{thm:subP locally-finite co-Heyting algebra} yields the following. 

	\begin{theorem}\label{thm:suboP locally-finite Heyting algebra}For any polyhedron $P$,
		$\Subo P$ is a Heyting algebra, and a subalgebra of $\Opens(P)$.
	\end{theorem}
 For any formula $\phi$ and polyhedron $P$, say that $P \vD \phi$ if and only if $\Subo P \vD \phi$ as a Heyting algebra. Call an intermediate logic \emph{polyhedrally-complete} if it is the logic of some class of polyhedra.  In \cite{tarski-polyhedra}, it is shown that $\IPC$ is polyhedrally-complete, being the logic of all polyhedra, while $\BD_n$ is the logic of all polyhedra of dimension at most $n$. It is also shown that all polyhedrally-complete logics must have the finite model property; cf.\  \cref{thm:logic of poly logic of tris} below.

\subsection{Triangulation subalgebras}
Let $\Sig$ be a triangulation of the polyhedron $P$. Then $\Sig \sse \Sub P$. Let $\Pc(\Sig)$ be the sublattice of $\Sub P$ generated by $\Sig$.

	\begin{lemma}\label{lem:subalgebra}
		$\Pc(\Sig)$ is a co-Heyting subalgebra of $\Sub P$.
	\end{lemma}

	\begin{proof}
		See \cite[Lemma 3.6]{tarski-polyhedra}.
	\end{proof}

	Call any algebra of the form $\Pc(\Sig)$ a \emph{triangulation subalgebra}.
	\begin{lemma}\label{lem:fin-gen subalgebras in triangulation subalgebras}
		Every finitely-generated subalgebra of $\Sub P$ is contained in some triangulation algebra.
	\end{lemma}

	\begin{proof}
		See \cite[Lemma~3.2]{tarski-polyhedra}. Essentially, this is the content of the Triangulation Lemma \ref{lem:triangulation lemma}.
	\end{proof}

	Turning now to the dual, every triangulation $\Sig$ of a polyhedron $P$ gives rise to a Heyting subalgebra $\Po(\Sig)$ of $\Subo P$, which we also call a \emph{triangulation subalgebra}, generated by the complements of the simplices in $\Sig$. 
	\begin{corollary}\label{cor:SuboP locally-finite}For any polyhedron $P$,
		$\Subo P$ is a locally-finite Heyting algebra.
	\end{corollary}

	\begin{proof}
		This follows from the dual of \cref{lem:fin-gen subalgebras in triangulation subalgebras} and the fact that triangulation subalgebras are finite.
	\end{proof}

	The algebra $\Po(\Sig)$, though not necessarily easy to visualise geometrically, is in fact precisely the algebra of upsets of the poset  $\Sigma$.

	\begin{lemma}\label{lem:UpP PoSig duality}
		The map:
		\begin{align*}
			\gamma^\ua \colon \Up \Sig &\to \Po(\Sig) \\
			U &\mapsto \bigcup_{\sig \in U}\Relint(\sig)
		\end{align*}
		is an isomorphism of Heyting algebras.
	\end{lemma}

	\begin{proof}
		See \cite[Lemma~4.3]{tarski-polyhedra}.
	\end{proof}

As a consequence of the preceding results, we have:

	\begin{theorem}\label{thm:logic of poly logic of tris}
		The logic of a polyhedron is the logic of its triangulations.
	\end{theorem}

	The following additional facts about triangulation algebras will be useful later on.

	\begin{lemma}\label{lem:properties of triangulation subalgebras}
		\begin{enumerate}[label=(\arabic*)]
			\item\label{item:determine; lem:properties of triangulation subalgebras}
				Triangulation algebras determine their corresponding triangulations. That is, for any two triangulations $\Sig$ and $\Delta$, if $\Po(\Sig)=\Po(\Delta)$ then $\Sig=\Delta$.
			\item\label{item:isomorphism; lem:properties of triangulation subalgebras}
				If $\Sig$ and $\Delta$ are triangulations which are isomorphic as posets then $\Po(\Sig) \cong \Po(\Delta)$.
			\item\label{item:refinement; lem:properties of triangulation subalgebras}
				If $\Delta$ refines $\Sig$, then $\Po(\Sig)$ is a subalgebra of $\Po(\Delta)$.
		\end{enumerate}
	\end{lemma}

	\begin{proof}
		\begin{enumerate}[label=(\arabic*)]
			\item 
				It follows from conditions \ref{item:downwards-closed; defn:simplicial complex} and \ref{item:intersection; defn:simplicial complex} on simplicial complexes that $\Pc(\Sig)$ consists exactly of the unions of elements of $\Sig$, and similarly for $\Delta$. Assume $\Po(\Sig) = \Po(\Delta)$, so that $\Pc(\Sig) = \Pc(\Delta)$, and take $\sig \in \Sig$. Then $\sig \in \Pc(\Delta)$, so $\sig = \bigcup S$ for some $S \sse \Delta$, and similarly each $\tau \in S$ is $\tau = \bigcup T_\tau$ for some $T_\tau \sse \Sig$. Hence:
				\begin{equation*}
					\sig = \bigcup\bigcup_{\tau \in S} T_\tau
				\end{equation*}
				But then by condition \ref{item:intersection; defn:simplicial complex} on $\Sig$, every $\rho \in \bigcup_{\tau \in S} T_\tau$ must either be equal to $\sig$ or be a proper face of $\sig$. Since $\Relint \sig$ contains no proper face of $\sig$, we must have $\sig \in T_\tau$ for some $\tau \in S$. But then $\sig \sse \tau \sse \sig$, and so $\sig \in \Delta$. Applying this argument also in the other direction, we get that $\Sig = \Delta$.
			\item
				This follows from \cref{lem:UpP PoSig duality}.
			\item 
				By \cref{lem:subdivision simplex union of simplicies}, every $\sig \in \Sig$ is the union of simplices in $\Delta$. Whence $\Sig \sse \Pc(\Delta)$. Therefore, by definition $\Pc(\Sig) \sse \Pc(\Delta)$. From this is follows that $\Po(\Sig) \sse \Po(\Delta)$.\qedhere
		\end{enumerate}
	\end{proof}

\subsection{PL homeomorphisms}
\label{ssec:polyhedral maps}

Let $P\subseteq \R^m$ and $Q\subseteq \R^n$ be polyhedra. A continuous function $f\colon P\to Q$ is  \emph{piecewise-linear}, or is a \emph{PL map}, if the graph of $f$ is a polyhedron in the product space $\R^m\times \R^n$. A  \emph{PL homeomorphism} is a PL map that is a homeomorphism. 	
\begin{proposition}\label{p:PLinv}
		The inverse of a PL homeomorphism is a PL homeomorphism.
	\end{proposition}
\begin{proof}
		See \cite[p.~6]{rourkesanderson1972}.
\end{proof}

\begin{proposition}\label{prop:PL homeomorphism HA isomorphism}	
If $P$ and $Q$ are PL homeomorphic then $\Subo{Q}$ and  $\Subo{P}$ are isomorphic Heyting algebras, and $\Logic(P)=\Logic(Q)$.
\end{proposition}

	\begin{proof}It is obvious that any homeomorphism between $P$ and $Q$ induces an isomorphism of their open-set lattices by taking inverse images. Since
		the inverse image of a subpolyhedron under a PL homeomorphism is again a subpolyhedron \cite[Corollary~2.5, p.~13]{rourkesanderson1972}, and in light of \cref{p:PLinv}, we see that when the homeomorphism is PL this isomorphism of open-set lattices descends to an isomorphism of distributive lattices between $\Subo P$ and $\Subo Q$. This implies that  $\Subo P$ and $\Subo Q$ are isomorphic as Heyting algebras, too, because the Heyting implication is uniquely determined by the underlying lattice structure, and the  proof is complete.\end{proof}

%
%
%
%

\section{The Nerve Criterion}
\label{sec:nerve criterion}

Given a poset $F$, its \emph{nerve}, $\N(F)$, is the collection of finite non-empty chains in $F$ ordered by inclusion. The following theorem is the first main contribution of the paper: 

\begin{theorem}[The Nerve Criterion]\label{thm:nerve criterion}
	A logic is polyhedrally-complete if, and only if, it is the logic of a class of finite posets closed under the nerve operator $\N$.
\end{theorem}

The utility of the Nerve Criterion is that it transforms logico-geometric questions into questions about finite posets, to which finite combinatorial methods are applicable.  

\begin{remark}
	We cannot strengthen the left-to-right direction to the following. ``If a logic $\Lo$ is polyhedrally-complete then $\FramesFin(\Lo)$ is closed under the nerve operator $\N$''. By \cref{cor:SL ploy-complete} below Scott's Logic $\SL$ is polyhedrally-complete. However $\FramesFin(\SL)$ contains the frame $F$ given in \cref{fig:SL frames not closed under N}. As can be seen there, the nerve $\N(F)$ does not validate $\SL$, since there is an up-reduction $\N(F) \cra \FScott$. Using the terminology introduced in \cref{sec:starlike completeness}, the problem is that while $F$ is $(2 \cdot 1)$-connected, it is not $(2 \cdot 1)$-diamond-connected.
\end{remark}

\begin{figure}
	\begin{equation*}
		\begin{tikzpicture}
			\begin{scope}
				\graph[poset] {
					r[x=0.5] -- {
						a1 -- a2,
						b[y=0.5]
					} -- t[x=0.5]
				};
				\node at (0.5,-1) {$F$};
			\end{scope}
			\begin{scope}[xshift=100]
				\graph[poset] {
					r[x=1.5] -- {
						{
							a1 -- {b1, b2},
							a2[x=-1] -- {b1, b3},
							a3[x=-1] -- {b2, b3}
						} -- c[x=1],
						{a3, a4[x=-1]} -- b4[x=-1]
					}
				};
				\node at (1.5,-1) {$\N(F)$};
			\end{scope}
		\end{tikzpicture}
	\end{equation*}
	\caption{An example showing that the $\FramesFin(\SL)$ is not closed under $\N$, even though \SL\ is polyhedrally-complete.}
	\label{fig:SL frames not closed under N}
\end{figure}

Achieving a proof of the Nerve Criterion will require considerable work with rational triangulations and their subdivisions. We next state the key intermediate result to be obtained. Let $A$ be a triangulation subalgebra of $\Subo P$, for some polyhedron $P$. By \cref{lem:properties of triangulation subalgebras} \ref{item:determine; lem:properties of triangulation subalgebras}, there is a unique triangulation $\Sig$ of $P$ such that $A = \Po(\Sig)$. For any $k \in \NN$, let $A^{(k)} \defeq \Po(\Sig^{(k)})$, where $\Sig^{(k)}$ is the $k$-th derived subdivision of $\Sigma$ (see \cref{ss:stellar}).

	\begin{theorem}\label{thm:kth derived subalgebra includes all finite subalgebras}
		Let $P$ be a polyhedron, and let $A$ be any triangulation subalgebra of $\Subo P$. For any finitely-generated subalgebra $B$ of $\Subo{P}$, there is $k \in \NN$ such that $B$ is isomorphic to a subalgebra of $A^{(k)}$.
	\end{theorem}

Sections~\ref{sec:rational polyhedra}–\ref{sec:putting-all-together} will be devoted to proving  Theorem \ref{thm:kth derived subalgebra includes all finite subalgebras}. The proof of the Nerve Criterion is completed in \cref{ssec:proof of the criterion}.

\subsection{Rational polyhedra and unimodular triangulations}\label{sec:rational polyhedra}

	The geometric intuition behind \cref{thm:kth derived subalgebra includes all finite subalgebras} is that any triangulation can be approximated from any other by taking iterated barycentric subdivisions.  One difficulty  with spelling out such an intuition is that if we start with a triangulation $\Sig$ on vertices with irrational coordinates, and try to approximate it using the iterated barycentric subdivisions of a triangulation on  vertices with rational coordinates, the approximations can never quite capture (a refinement of) $\Sig$. The approach taken here is effectively to show that it suffices to restrict attention to the rational case. In order to make this idea precise, we need  tools on rational triangulations that go beyond the standard polyhedral topology handbooks, which typically deal with the real case only. For these tools we  mainly use \cite{mundici2011advanced} as a background reference.

	A polytope in $\R^n$ is \emph{rational} if it may be written as the convex hull of finitely many points in $\Q^n \sse \R^n$. A polyhedron in $\R^n$ is \emph{rational} if it may be written as a union of a finite collection of rational polytopes. A simplicial complex $\Sig$ is \emph{rational} if it consists of rational simplices. Note that when this is the case, $\abs\Sig$ is a rational polyhedron.

	For any $x \in \Q^n \sse \R^n$, there is a unique way to write out $x$ in coordinates as $x = (\frac{p_1}{q_1},\ldots,\frac{p_n}{q_n})$ such that for each $i$, we have $p_i,q_i \in \Z$ coprime. The \emph{denominator} of $x$ is defined:
	\begin{equation*}
		\Den(x) \defeq \lcm\{q_1, \ldots, q_n\}
	\end{equation*}
	Thus, $\Den(x)=1$ if and only if $x$ has integer coordinates. Letting $q = \Den(x)$, the \emph{homogeneous correspondent} of $x$ is defined to be the integer vector:
	\begin{equation*}
		\wt x \defeq \left(\frac{qp_1}{q_1},\ldots,\frac{qp_n}{q_n}, q\right)
	\end{equation*}

	A rational $d$-simplex $\sig = x_0 \cdots x_d$ is \emph{unimodular} if there is an $(n+1)\times(n+1)$ matrix with integer entries whose first $d+1$ columns are $\wt{x}_0,\ldots,\wt{x}_d$, and whose determinant is $\pm 1$. This is equivalent to requiring that the set $\{\wt{x}_0,\ldots,\wt{x}_d\}$ can be completed to a $\Z$-module basis of $\Z^{d+1}$. A simplicial complex is \emph{unimodular} if each one of its simplices is unimodular.

\subsection{Farey subdivisions}
	\begin{proposition}\label{prop:mediant}
		Given a rational simplex $\sigma$ with vertices $x_0, \ldots, x_d \in \Q^n\subseteq\R^n$, there is a unique $m \in \Q^n$ such that $\wt m = \sum_{i=0}^d\wt{x}_i$. Moreover,  $m\in\Relint{\sigma}$.
	\end{proposition}

	\begin{proof}
		Let $H_{n+1} \sse \R^{n+1}$ be the hyperplane specified by:
		\begin{equation*}
		 	H_{n+1} \defeq \{(x_1, \ldots, x_{n+1}) \in \R^{n+1} \mid x_{n+1} =1\}
		\end{equation*} 
		Identify $\Q^n$ with the set of rational points of $H_{n+1}$ via the map $(q_1,\ldots,q_n)\mapsto(q_1,\ldots,q_n,1)$. Under this identification, $\wt m$ lies in the rational cone:
		\begin{equation*}
			\left\{\sum_{i=0}^d c_i\wt x_i\mid c_i\in\R, c_i \geq 0\right\}
		\end{equation*}
		A routine computation then proves the geometrically evident fact that $m$ is the point of intersection of the line spanned in $\R^{n+1}$ by the vector $\wt m$, with the hyperplane $H_{n+1}$; from which the result follows.
	\end{proof}

	\noindent The element $m \in \Q^n$  in Proposition \ref{prop:mediant} is called the \emph{Farey mediant} of (the vertices of) the simplex $\sigma$. Note that when $d=0$, i.e.\ when $\sigma$ is a vertex of $\Sigma$, then $m$ coincides with the vertex $\sigma$. Also observe that the Farey mediant and the barycentre of $\sigma$   are in general distinct, though both lie in $\Relint{\sigma}$.

	We can now define a specific type of stellar subdivision based on Farey mediants, cf.\ \cite[\S5.1, p.~55]{mundici2011advanced}.
 Let $\Sigma$ be a simplicial complex,  let $\sigma\in\Sigma$, and let $m$ be the Farey mediant of $\sigma$. The \emph{elementary Farey subdivision of $\Sigma$ at $m$}   is the elementary stellar subdivision of $\Sigma$ at $m$. In general, the triangulation $\Delta$ is a \emph{Farey subdivision} of $\Sigma$ if it is obtained from the latter via finitely many successive elementary Farey subdivisions.

At the combinatorial level, Farey and barycentric subdivisions are indiscernible:

	\begin{lemma}\label{lem:elementary Farey and barycentric subdivisions isomorphic}
		Let $\Sig, \Delta$ be simplicial complexes with $\Sig$ rational, assume that $\gamma \colon \Sig \to \Delta$ is an isomorphism of $\Sig$ and $\Delta$ as posets, let $\sig \in \Sig$, and let $m$ be the Farey mediant of $\sig$. Then the elementary Farey subdivision of $\Sig$ at  $m$ and the elementary barycentric subdivision of $\Delta$ at the barycentre of $\gamma(\sig)$ are isomorphic as posets.
	\end{lemma}

	\begin{proof}
		Indeed, at the level of posets, elementary Farey subdivision and elementary bary\-centric subdivision are the same operation, as direct inspection of the definitions confirms. For further details see also \cite[\S III]{alexander1930}.
	\end{proof}

However, going beyond the combinatorial level, the construction of universal approximations of arbitrary rational polyhedra does require Farey subdivisions and cannot be done with barycentric ones. This is made precise in the following fundamental fact of rational polyhedral geometry.
	\begin{lemma}[The De Concini-Procesi Lemma]\label{lem:De Concini-Procesi lemma}
		Let $P$ be a rational polyhedron, and let $\Sig$ be a unimodular triangulation of $P$. There exists a sequence $(\Sig_i)_{i \in \NN}$ of unimodular triangulations of $P$ with $\Sig_0=\Sig$ such that:
		\begin{enumerate}[label=(\alph*)]
			\item For each $i\in \NN$, $\Sig_{i+1}$ is an elementary Farey subdivision of $\Sig_i$, and
			\item For any rational polyhedron $Q\sse P$, there is $i\in\NN$ such that $\Sig_i$ triangulates $Q$.
		\end{enumerate}
	\end{lemma}

	\begin{proof}
		See \cite[Theorem~5.3, p.~57]{mundici2011advanced}.
	\end{proof}

\subsection{From \texorpdfstring{$\R$}{R} to \texorpdfstring{$\Q$}{Q}}

To deploy the power of \cref{lem:De Concini-Procesi lemma}, we  need to relate general polyhedra to rational polyhedra, and general triangulations to unimodular ones.

	\begin{lemma}\label{lem:general to rational polyhedra and simplicial complexes}
		Let $P$ be a polyhedron, and let $\Sigma$ be a triangulation of $P$. There exist an integer $n\in\NN$, a rational polyhedron $Q\subseteq\R^n$, and a unimodular triangulation $\Delta$ of $Q$ such that $P$ and $Q$ are PL-homeomorphic via a map that induces an isomorphism of $\Sigma$ and $\Delta$ as posets.
	\end{lemma}

	\begin{proof}
		This is a standard argument. Fix a bijection $\beta$ from the vertices of $\Sigma$ to the standard basis of $\R^n$, where $n$ is the number of vertices in $\Sigma$. Take a simplex $\sigma = x_0 \cdots x_d$ in $\Sigma$. Note that the points $\beta(x_0), \ldots, \beta(x_d)$ are affinely independent; let $\alpha(\sig)$ be the $d$-simplex spanned by their convex hull: $\alpha(\sig) \defeq \Conv\{\beta(x_0), \ldots, \beta(x_d)\}$. Since the vertices of $\alpha(\sig)$ are standard basis elements, $\alpha(\sig)$ is a unimodular simplex by definition. Let $f_\sigma \colon \sigma \to \alpha(\sig)$ be the linear map determined by $f_\sigma(x_i) = \beta(x_i)$ for each $i$, and let $g_\sigma \colon \alpha(\sig) \to \sigma$ be its inverse, determined by $g_\sigma(\beta(x_i)) = x_i$.

		Now, let $Q \defeq \bigcup_{\sigma \in \Sigma}\alpha(\sig)$. For any simplices $\sigma \preceq \tau$, the map $f_\sigma$ agrees with $f_\tau$ on $\sigma$. Hence we may glue these maps together to form a map $f \colon P \to Q$, i.e. $f(x)=f_\sigma(x)$, where $\sigma$ is any simplex of $\Sigma$ containing $x$. Similarly, we may glue together the maps $g_\sigma$ for $\sigma \in \Sigma$ to form an inverse to $f$. By definition $f$ is a PL homeomorphism. Finally, note that $\Delta \defeq \{\alpha(\sig) \mid \sigma \in \Sigma\}$ is a triangulation of $Q$, and that $f$ induces the poset isomorphism $\sigma \mapsto \alpha(\sig)$ between $\Sigma$ and $\Delta$.
	\end{proof}

	\begin{lemma}\label{lem:farey subdivision refined by barycentric}
		Let $\Sig$ be a unimodular triangulation of the rational polyhedron $P$, and suppose $\Sig'$ is a Farey subdivision of $\Sig$. There is a triangulation $\Delta$ of $P$ which is isomorphic as a poset to $\Sig'$, and $k \in \NN$ such that $\Sig^{(k)}$ refines $\Delta$.
	\end{lemma}

	\begin{proof}
		The proof works by replacing each elementary Farey subdivision by an elementary barycentric subdivision. We induct on the number $m \in \NNp$ of elementary Farey subdivisions needed to obtain $\Sig'$ from $\Sig$. If $m=1$, let $\sigma$ be the simplex of $\Sig$ being subdivided at its Farey mediant. Then the first barycentric subdivision $\Sig^{(1)}$ of $\Sig$ refines the elementary barycentric subdivision $\Sig^*$ of $\Sig$ at the barycentre of $\sigma$. By \cref{lem:elementary Farey and barycentric subdivisions isomorphic}, $\Sig^*$ and $\Sig'$ are isomorphic.

		For the induction step, suppose $m>1$, and write $(\Sig_i)_{i=0}^m$ for the finite sequence of triangulations connecting $\Sig=\Sig_0$ to $\Sig'=\Sig_m$ through elementary Farey subdivisions. By the induction hypothesis, there is $k\in\NN$ such that $\Sig^{(k)}$ refines a triangulation $\Delta$ isomorphic to $\Sig_{m-1}$; let us fix one such isomorphism $\gamma$. Let $\sigma$ be the $d$-simplex of $\Sig_{m-1}$ that must be subdivided through its Farey mediant in order to obtain $\Sig_m$. Let further $\delta$ be the simplex of $\Delta$ that corresponds to $\sigma$ through the isomorphism $\gamma$. Since the $d$-simplices are exactly the height-$d$ elements of $\Delta$, we get that $\delta$ is a $d$-simplex. Then $\Sig^{(k+1)}$ refines $\Delta^*$, the latter denoting the elementary barycentric subdivision of $\Delta$ at the barycentre of $\delta$. But $\Delta$ is isomorphic to $\Sig_{m-1}$, and therefore by \cref{lem:elementary Farey and barycentric subdivisions isomorphic}, $\Delta^*$ is isomorphic to $\Sig_{m}$.
	\end{proof}

Finally, we shall need the  non-trivial fact that  arbitrary triangulations of a rational polyhedron realise no more combinatorial types than its rational triangulations; this is due to Meurig Beynon:
	\begin{lemma}[Beynon's Lemma]\label{lem:Beynons Lemma}
		Let $P$ be a rational polyhedron, and let $\Sigma$ be a triangulation of $P$. There exists a rational triangulation of $P$ which is isomorphic as a poset to $\Sigma$.
	\end{lemma}

	\begin{proof}
		This is the main result of \cite{beynon1977}.
	\end{proof}

\subsection{End of proof of \texorpdfstring{\cref{thm:kth derived subalgebra includes all finite subalgebras}}{Theorem \ref{thm:kth derived subalgebra includes all finite subalgebras}}}\label{sec:putting-all-together}

	\begin{proof}[Proof of \cref{thm:kth derived subalgebra includes all finite subalgebras}]
		Let $\Sigma$ be the triangulation of $P$ such that $A=\Po{(\Sigma)}$. Using \cref{lem:general to rational polyhedra and simplicial complexes}, \cref{lem:properties of triangulation subalgebras} \ref{item:isomorphism; lem:properties of triangulation subalgebras} and \cref{prop:PL homeomorphism HA isomorphism} we may assume without loss of generality that $P$ is rational and $\Sigma$ is unimodular. By \cref{lem:fin-gen subalgebras in triangulation subalgebras}, there is a triangulation $\Delta$ of $P$ such that $B$ is isomorphic to a subalgebra of $\Po{(\Delta)}$. By Beynon's Lemma \ref{lem:Beynons Lemma} and \cref{lem:properties of triangulation subalgebras} \ref{item:isomorphism; lem:properties of triangulation subalgebras}, we may assume that $\Delta$ is rational (and hence each member of $B$ is, too). By the De Concini-Procesi Lemma \ref{lem:De Concini-Procesi lemma}, there is a Farey subdivision $\Sigma'$ of $\Sigma$ that refines $\Delta$. Therefore by \cref{lem:properties of triangulation subalgebras} \ref{item:refinement; lem:properties of triangulation subalgebras}, $B$ is isomorphic to a subalgebra of $\Po{(\Sigma')}$. By \cref{lem:farey subdivision refined by barycentric}, there is $k\in\NN$ such that $\Sigma^{(k)}$ refines $\Sigma'$ up to isomorphism. Hence by \cref{lem:properties of triangulation subalgebras} \ref{item:refinement; lem:properties of triangulation subalgebras} again, $A^{(k)}$ contains a subalgebra isomorphic to $\Po{(\Sigma')}$, and therefore also a subalgebra isomorphic to $B$. This completes the proof.
	\end{proof}

\subsection{Nerves,  subdivisions, and geometric realisations}\label{ssec:nerves and deriveds}

	The reason that \cref{thm:kth derived subalgebra includes all finite subalgebras} is relevant to the  Nerve Criterion is the following classical connection between nerves of posets and derived subdivisions, which is  deeply rooted in the work of Pavel S. Alexandrov.
	
		\begin{proposition}\label{prop:Sd cong N}
		Let $\Sig$ be a simplicial complex, regarded as a poset under inclusion of faces. Then the barycentric subdivision of $\Sig$ is isomorphic as a poset to the nerve of $\Sig$:
		\begin{equation*}
			\Sd \Sig \cong \N(\Sig)
		\end{equation*}
	\end{proposition}

	\begin{proof}The proof flows readily from the definitions, and is in any case available from multiple sources; see e.g. \cite[Ch.\ IV, \S 2.2]{alexandrov1998} (Alexandrov's own textbook treatment), or \cite[Proposition~2.5.10, p.~51]{maunder1980algebraic}, or \cite[\S3]{ranickiweiss2012}. Details are left to the reader.
	\end{proof}

	\begin{corollary}\label{cor:logic of P logic of iterated nerves of Sig}
		For $P$ a polyhedron and $\Sig$ a triangulation of $P$ we have:
		\begin{equation*}
			\Logic(P) = \Logic(\N^k(\Sig) \mid k \in \NN)
		\end{equation*}
	\end{corollary}

	\begin{proof}
		Indeed:
		\begin{align*}
			\Logic(P)
				&= \Logic(\Subo P) \\
				&= \Logic (A \mid A\text{ finitely-generated subalgebra of }\Subo P) \tag{\cite[Ch.~7]{chagrovzakharyaschev1997}}\\
				&= \Logic(\Po(\Sig^{(k)}) \mid k \in \NN) \tag{\cref{thm:kth derived subalgebra includes all finite subalgebras}} \\
				&= \Logic(\Sig^{(k)} \mid k \in \NN) \tag{as above} \\
				&= \Logic(\N^k(\Sig) \mid k \in \NN) \tag{\cref{prop:Sd cong N}}
		\end{align*}
	\end{proof}

	In order to complete a proof of the Nerve Criterion, we will also need to use geometric realisations of finite posets via nerves, another classical tool. Let $F = \{x_1, \ldots, x_m\}$ be a finite poset, and let $e_1, \ldots, e_m$ be the standard basis vectors of $\R^m$. The set:
		\begin{equation*}
		\nabla F \defeq \{\Conv\{e_{i_1}, \ldots, e_{i_k}\} \mid \{x_{i_1}, \ldots, x_{i_k}\} \in \N(F)\}
	\end{equation*}
	can be proved to be a triangulation by elementary arguments; its underlying polyhedron $\abs{\nabla F}$ is the \emph{geometric realisation} of $F$. For us, the key fact about geometric realisations is:
	\begin{lemma}\label{l:geomrea}Let $F$ be a finite poset. The map $\max \colon \N(F) \to F$, which sends a chain to is maximum element, is a p-morphism, and $\Logic(\abs{\nabla F}) \sse \Logic(F)$.
	\end{lemma}
	\begin{proof}The first statement is easy to verify by direct inspection, and a detailed proof was already given in \cite[p.\ 389]{tarski-polyhedra}. For the second statement, observe first that $\nabla F$ and $\N(F)$ are isomorphic as posets (under inclusion), by their definitions. Thus, by \cref{prop:up-reduction logic containment}, the surjective p-morphism $\nabla F \to F$ yields  $\Logic(\nabla F)\sse\Logic(F)$. But $\Up(\nabla F)$ is a subalgebra of $\Subo(\abs{\nabla(F)})$ by \cref{lem:UpP PoSig duality} together with \cref{lem:subalgebra}, so that   $\Logic(\abs{\nabla F}) \sse \Logic(\nabla F) \sse \Logic(F)$, as was to be shown.
	\end{proof}


\subsection{End of proof of the Nerve Criterion}\label{ssec:proof of the criterion}
	\begin{proof}[Proof of \cref{thm:nerve criterion}, the Nerve Criterion]
		Assume that $\Lo$ is the logic of a class $\CC$ of polyhedra. For each $P \in \CC$ fix a triangulation $\Sig_P$, and let:
		\begin{equation*}
			\CC^* \defeq \{\N^k(\Sig_P) \mid P \in \CC\text{ and }k \in \NN\}
		\end{equation*}
		Then:
		\begin{align*}
			\Logic(\CC^*)
				&= \bigcap_{P \in \CC}\Logic(\N^k(\Sig_P) \mid k \in \NN) \\
				&= \bigcap_{p \in \CC}\Logic(P) \tag{\cref{cor:logic of P logic of iterated nerves of Sig}} \\
				&= \Logic(\CC) = \Lo
		\end{align*}
		Conversely, assume that $\Lo = \Logic(\D)$, where $\D$ is a class of finite frames closed under $\N$. Let:
		\begin{equation*}
			\D_* \defeq \{\abs{\nabla(F)} \colon F \in \D\}
		\end{equation*}	
	where $\abs{\nabla(F)}$ is the geometric realisation of $F$ as in \cref{ssec:nerves and deriveds}.
		We will show that $\Lo = \Logic(\D_*)$. First suppose that $\Lo \nvd \phi$, so that $F \nvD \phi$ for some $F \in \D$. Then we have that $\abs{\nabla(F)} \nvD \phi$, so that $\Logic(\D_*) \nvd \phi$. Conversely, suppose that $\Logic(\D_*) \nvd \phi$, so that $\abs{\nabla(F)} \nvD \phi$ for some $F \in \D$. By definition $\nabla(F)$ is a triangulation of $\abs{\nabla(F)}$, hence by \cref{cor:logic of P logic of iterated nerves of Sig} there is $k \in \NN$ such that $\nabla(F)^{(k)} \nvD \phi$. But $\nabla(F) \cong \N(F)$ by definition, and so by \cref{prop:Sd cong N} we get $\N^{k+1}(F) \cong\nabla(F)^{(k)}$. Thus, as $\D$ is closed under $\N$, we get that $\Lo \nvd \phi$.
	\end{proof}

\section{Polyhedrally-incomplete logics}
\label{sec:stable logics}

	In this section, we apply the Nerve Criterion to show that every stable logic \blue{other than \IPC\ is polyhedrally-incomplete. A logic $\Lo$ is \emph{stable} if $\FramesRoot(\Lo)$ is closed under monotone images.
	(We point out that the original definition of \cite[Def.~6.6]{BB17} used Esakia spaces. However, it can be shown that these definitions are equivalent \cite[Theorem 3.3.17]{julia-thesis}.)}

	\begin{proposition}
		The following well-known logics\footnote{For more information on these logics see \cite[Table~4.1, p.~112]{chagrovzakharyaschev1997}.} are all stable.
		\begin{enumerate}[label=(\roman*)]
			\item The logic of weak excluded middle, $\KC = \IPC + (\neg p \vee \neg\neg p)$.
			\item Gödel-Dummett logic, $\LC = \IPC + (p \ra q) \vee (q \ra p)$.
			\item $\LC_n = \LC + \BD_n$.
			\item The logic of bounded width $n$, $\BW_n = \IPC + \bigvee_{i=0}^n(p_i \ra \bigvee_{j \neq i} p_j)$.
			\item The logic of bounded top width $n$, defined:
			\begin{equation*}
				\BTW_n \defeq \bigwedge_{0 \leq i < j \leq n} \neg(\neg p_i \wedge \neg p_j) \ra \bigvee_{i=0}^n( \neg p_i \ra \bigvee_{j \neq i} \neg p_j)
			\end{equation*}
			\item The logic of bounded cardinality $n$, defined:
			\begin{equation*}
				\BC_n \defeq p_0 \vee (p_0 \ra p_1) \vee ((p_0 \wedge p_1) \ra p_2) \vee \cdots \vee ((p_0 \wedge \cdots \wedge p_{n-1}) \ra p_n)
			\end{equation*}
		\end{enumerate}
	\end{proposition}

	\begin{proof}
		See \cite[Theorem~7.3]{BB17}.
	\end{proof}

	\noindent In fact:

	\begin{theorem}
		There are continuum-many stable logics.
	\end{theorem}

	\begin{proof}
		See \cite[Theorem 6.13]{BB17}.
	\end{proof}

	\begin{theorem}\label{t:stablefmp}
		Every stable logic has the finite model property.
	\end{theorem}

	\begin{proof}
		See \cite[Theorem 6.8]{BB17}.
	\end{proof}

	However,  \cref{t:stablefmp} notwithstanding:

	\begin{theorem}\label{thm:stable logics poly incomplete}
		If $\Lo$ is a stable logic other than $\IPC$, and $\Frames(\Lo)$ contains a frame of height at least $2$, then $\Lo$ is not polyhedrally-complete.
	\end{theorem}

	\begin{proof}
		Let $\Lo$ be a polyhedrally-complete stable logic of height at least 2. We show that $\Lo = \IPC$. 
		
		By the Nerve Criterion \ref{thm:nerve criterion}, there is a class $\CC$ of finite frames closed under $\N$ such that $\Lo = \Logic(\CC)$. Since $\Frames(\Lo)$ contains a frame of height at least $2$, we must have $\Lo \nvd \BD_1$. Since $\Lo = \Logic(\CC)$, there is therefore $F \in \CC$ such that $\height(F) \geq 2$. This means there are $x_0, x_1, x_2 \in F$ with $x_0 < x_1 < x_2$. Without loss of generality, we may assume that $x_2$ is a top element and that $x_1$ is an immediate predecessor of $x_2$ and $x_0$ an immediate predecessor of $x_1$. Now, by assumption $\N^k(F) \in \CC$ for every $k \in \NN$. Let us examine the structure of these frames a little. Note that $\{x_0,x_1,x_2\}$ is a chain. Let $X$ be a maximal chain in $\Ds{x_0}$. We have the following relations occurring in $\N(F)$. 
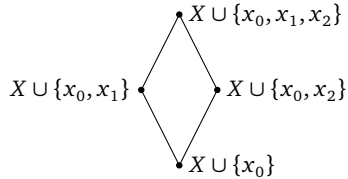
		\begin{equation*}
			\begin{tikzpicture}
				\graph[poset] {
					z[x=0.5,label=right:$X \cup \{x_0\}$] -- 
						{
							a[label=left:{$X \cup \{x_0,x_1\}$}],
							b[label=right:{$X \cup \{x_0,x_2\}$}]
						} -- c[x=0.5,label=right:{$X \cup \{x_0,x_1,x_2\}$}];
				};
			\end{tikzpicture}
		\end{equation*}
		Moreover, by assumptions on $x_0, x_1, x_2$ and $X$, we have that $X \cup \{x_0,x_1,x_2\}$ is a top element of $\N(F)$, with $X \cup \{x_0,x_1\}$ and $X \cup \{x_0,x_2\}$ immediate predecessors, and $X \cup \{x_0\}$ an immediate predecessor of those. So, we may apply this argument once more, to obtain the following structure sitting at the top of $\N^2(F)$.
		\begin{equation*}
			\begin{tikzpicture}
				\graph[poset] {
					z[x=1] -- {
						a1 -- b1[x=0.5],
						a2 -- {b1, b2[x=0.5]},
						a3 -- b2
					};
				};
			\end{tikzpicture}
		\end{equation*}
		Iterating, we see that at the top of $\N^k(F)$ we have the following structure.
		
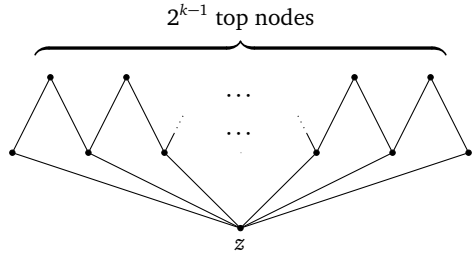
\begin{equation*}
			\begin{tikzpicture}
				\graph[poset] {
					a[x=3cm,label=below:$z$] -- {
						b0 -- c0[x=0.5cm],
						b1 -- {c0, c1[x=0.5cm]},
						b2 -- c1,
						-!- e0[minimum size=0],
						b3 -- c2[x=0.5cm],
						b4 -- {c2, c3[x=0.5cm]},
						b5 -- c3
					}
				};
				\draw (b2) -- ++(0.125,0.25) edge[dotted] +(0.125,0.25);
				\draw (b3) -- ++(-0.125,0.25) edge[dotted] +(-0.125,0.25);
				\node[above=1 of a] {$\cdots$};
				\node[above=1.5 of a] {$\cdots$};
				\draw[decoration={calligraphic brace,amplitude=5pt}, decorate, line width=1.25pt] 
					(0.3,2.3) -- (5.7,2.3);
				\node at (3,2.8) {$2^{k-1}$ top nodes};
			\end{tikzpicture}
		\end{equation*}
		Let $z$ be the base element of this structure, as indicated. Now, take $k \in \NN$ and let $\{t_1, \ldots, t_m\}$ be the top nodes of $\N^k(F)$ produced by this construction, where $m = 2^{k-1}$. By \cref{prop:up-reduction logic containment}, $\us z \in \FramesRoot(\Lo)$. 
		
		Let now $G$ be an arbitrary poset with up to $m$ elements $\{y_1, \ldots, y_m\}$ (possibly with duplicates) plus a root $\bot$. Define $f \colon \us z \to G$ as follows.
		\begin{equation*}
			x \mapsto \left\{
			\begin{array}{ll}
				y_i & \text{if }x=t_i, \\
				\bot & \text{otherwise}.
			\end{array}
			\right.
		\end{equation*}
		Then $f$ is monotonic. Since $\Lo$ is stable, this means that $G \in \FramesRoot(\Lo)$. Thus (since, by \cref{prop:ipc Kripke complete and has fmp} and \cref{cor:logic of frames logic of rooted frames}, $\IPC$ is the logic of finite rooted frames) we get that $\Lo = \IPC$. 
	\end{proof}

\section{Polyhedrally-complete logics: starlike completeness}
\label{sec:starlike completeness}

	In this section, we use the Nerve Criterion to establish a class of logics which are polyhedrally-complete. This constitutes the second main result of the paper.

\subsection{Starlike trees}

	A finite poset $T$ is a \emph{tree} if it has a root $\bot$, and every other $x \in T \setminus \{\bot\}$ has exactly one immediate predecessor. A \emph{branch} in $T$ is a maximal chain. Say that $T$ is a \emph{starlike tree} if every $x \in T \setminus \{\bot\}$ has at most one immediate successor. (The terminology `starlike' comes from graph theory, see \cite{watanabe_schwenk_1979}.) A starlike tree is determined by the multiset of its branch heights, which motivates the following notation.

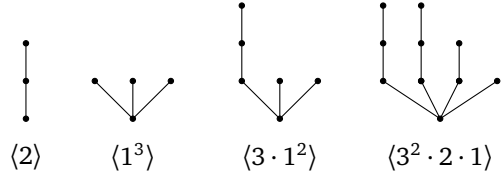
\begin{figure}
		\begin{equation*}
			\begin{tikzpicture}[ scale=0.5]
				\begin{scope}
					\draw[every node/.style=point] 
						(0,0) node {}
						-- (0,1) node {}
						-- (0,2) node {};
					\node at (0,-1) {$\starlike{2}$};
				\end{scope}
				\begin{scope}[xshift=80]
					\draw[every node/.style=point] 
						(0,0) node {} -- (-1,1) node {}
						(0,0) -- (0,1) node {}
						(0,0) -- (1,1) node {};
					\node at (0,-1) {$\starlike{1^3}$};
				\end{scope}
				\begin{scope}[xshift=190]
					\draw[every node/.style=point] 
						(0,0) node {} 
						-- (-1,1) node {} 
						-- (-1,2) node {}
						-- (-1,3) node {}
						(0,0) -- (0,1) node {}
						(0,0) -- (1,1) node {};
					\node at (0,-1) {$\starlike{3 \cdot 1^2}$};
				\end{scope}
				\begin{scope}[xshift=310]
					\draw[every node/.style=point] 
						(0,0) node {} 
						-- (-1.5,1) node {} 
						-- (-1.5,2) node {}
						-- (-1.5,3) node {}
						(0,0) -- (-0.5,1) node {}
						-- (-0.5,1) node {} 
						-- (-0.5,2) node {}
						-- (-0.5,3) node {}
						(0,0) -- (0.5,1) node {}
						-- (0.5,1) node {} 
						-- (0.5,2) node {}
						(0,0) -- (1.5,1) node {};
					\node at (0,-1) {$\starlike{3^2 \cdot 2 \cdot 1}$};
				\end{scope}
			\end{tikzpicture}
		\end{equation*}
		\caption{Some examples of starlike trees}
		\label{fig:starlike tree examples}
	\end{figure}

	Let $n_1, \ldots, n_k, m_1, \ldots, m_k \in \NNp$, with $n_1, \ldots, n_k$ distinct. Then let us define $T=\starlike{n_1^{m_1} \cdots n_k^{m_k}}$ as the starlike tree, uniquely determined to within an isomorphism, with the property that if we remove the root $\bot$ we are left with exactly, for each $i$, $m_i$ chains of length $n_i$. Let $\starlike \ep = \bullet$, the singleton poset. Call $\alpha = n_1^{m_1} \cdots n_k^{m_k}$ (or $\ep$) the \emph{signature} of $T$. We will always assume that $n_1 > n_2 > \cdots > n_k$.  See \cref{fig:starlike tree examples} for some examples of starlike trees together with their signatures. We will sometimes write $1^0$ for $\ep$.

 The \emph{length} of a signature $\alpha = n_1^{m_1} \cdots n_k^{m_k}$ is defined as $\abs{\alpha} \defeq m_1 + \cdots + m_k$. Let $\abs\ep \defeq 0$. For $j \leq \abs\alpha$, the \emph{$j$th height}, $\alpha(j)$, is $n_i$, where:
	\begin{equation*}
		m_1 + \cdots + m_{i-1} \leq j < m_1 + \cdots + m_i
	\end{equation*}

	Let $\alpha$ and $\beta$ be signatures. Say that $\alpha \leq \beta$ if $\abs\alpha \leq \abs\beta$ and for every $j \leq \abs\alpha$ we have $\alpha(j) \leq \beta(j)$. Considering the examples in \cref{fig:starlike tree examples}, we have the following relations:
	\begin{equation*}
		1^3 < 3 \cdot 1^2 < 3^2 \cdot 2 \cdot 1, \quad 2 < 3 \cdot 1^2
	\end{equation*}
	Note that if $\alpha = n_1^{m_1} \cdots n_k^{m_k}$, we have $\alpha \leq \beta$ if and only if $\abs\alpha \leq \abs\beta$ and for every $i \leq k$, we have:
	\begin{equation*}
		\beta(m_1 + \cdots + m_i) \geq n_i
	\end{equation*}

	\begin{proposition}\label{prop:signature ordering induces p-morphisms}
		If $\alpha \leq \beta$ then there is a p-morphism $\starlike\beta \to \starlike\alpha$.
	\end{proposition}

	\begin{proof}
		We can realise $\starlike\alpha$ as a downwards-closed subset of $\starlike\beta$. The p-morphism $f \colon \starlike\alpha \to \starlike\beta$ is then defined as follows. First, $f$ is the identity on $\starlike\alpha$. Second, for any branch of $\starlike\beta$ which contains a branch of $\starlike\alpha$, we let $f$ send any remaining elements to the maximum of the branch of $\starlike\alpha$. Finally, any remaining elements of $\starlike\beta$ are mapped to the maximum element of some fixed branch in $\starlike\alpha$. A routine calculation shows that $f$ is a p-morphism.
	\end{proof}

	Note that the starlike tree $\starlike k$ is the chain on $k+1$ elements; we will use this  notation for chains from now on. \blue{We will write the signature as $k^1$, to disambiguate it from $k$ as a number.} For $k \in \NNp$, the \emph{$k$-fork} is the starlike tree $\starlike{1^k}$.

\subsection{Starlike logics}

	We are now in a position to define the principal class of logics that will be investigated in this section. Let $\Starlikes \defeq \{\alpha \text{ signature} \mid \alpha \neq 1^2\}$. Take $\Lam \sse \Starlikes$ (possibly infinite). The \emph{starlike logic} $\SFL(\Lam)$ based on $\Lam$ is the logic axiomatised by $\IPC$ plus $\jsl\alpha$ for each $\alpha \in \Lam$. Write $\SFL(\alpha_1, \ldots, \alpha_k)$ for $\SFL(\{\alpha_1, \ldots, \alpha_k\})$.

	\begin{remark}
		For an explanation as to why the difork $\starlike{1^2}$ is omitted, see \cref{prop:poly-complete extension of difork implies CPC} below and the preceding discussion.
	\end{remark}

	\begin{proposition}\label{prop:SL starlike logic}
		$\SL = \SFL(2 \cdot 1)$. So Scott's Logic is a starlike logic.
	\end{proposition}

	\begin{proof}
		See \cite[\S9 and Table~9.7, p.~317]{chagrovzakharyaschev1997}.
	\end{proof}

	Let us examine what $\SFL(\Lam)$ `means' in terms of its class of frames. The formula $\jsl\alpha$ turns out to express a kind of connectedness property. We make this precise using the following definitions.

	Let $F$ be a finite poset. A \emph{path} in $F$ is a sequence $p=x_0\cdots x_k$ of elements of $F$ such that for each $i$ we have $x_i < x_{i+1}$ or $x_i > x_{i+1}$. Write $p \colon x_0 \rsa x_k$. The path $p$ is \emph{closed} if $x_0=x_k$. The poset $F$ is \emph{path-connected} if between any two points there is a path.

	\begin{lemma}\label{lem:finite poset path-connected iff connected}
		A poset is path-connected if, and only if, it is connected as a topological space.
	\end{lemma}

	\begin{proof}
		See \cite[Lemma~3.4]{Bezhanishviligabelaia2011}.
	\end{proof}

	\noindent A \emph{connected component} of $F$ is a subposet $U \sse F$ which is connected as a topological subspace and is such that there is no connected $V$ with $U \subset V$.

	\begin{lemma}\label{lem:properties of connectedness} Let $F$ be a poset.
		\begin{enumerate}[label=(\arabic*)]
			\item The connected components of $F$ partition $F$.
			\item The connected components of $F$ are downwards-closed and upwards-closed.
		\end{enumerate}
	\end{lemma}

	\begin{proof}
		\blue{These results follow straightforwardly from the fact that by \cref{lem:finite poset path-connected iff connected} the connected components are exactly the equivalence classes under the relation `there is a path from $x$ to $y$'.}
	\end{proof}

	Define $\concomps(F)$ to be the set of connected components of $F$. The \emph{connectedness type} $\contype(F)$ of $F$ is the signature $n_1^{m_1} \cdots n_k^{m_k}$ such that $\concomps(F)$ contains for each $i$ exactly $m_i$ sets of height $n_i-1$, and nothing else. Let $\contype(\es) \defeq \ep$.

	\begin{remark}
		Note that when $F$ is connected, $\contype(F)=n+1$, where $n = \height(F)$.
	\end{remark}

	Let $\alpha > \ep$ be a signature. An \emph{$\alpha$-partition} of $F$ is an open partition in which the number and heights of the connected components are specified by $\alpha$. In other words, it is a partition:
	\begin{equation*}
		F = C_1 \sqcup \cdots \sqcup C_{\abs{\alpha}}
	\end{equation*}
	into open sets such that $C_j$ has height at least $\alpha(j)-1$. For notational uniformity, say that $F$ has an $\ep$-partition if $F=\es$. The following lemma is a straightforward consequence of the definitions.

	\begin{lemma}\label{lem:alpha partition iff alpha leq contype}
		A finite poset $F$ has an $\alpha$-partition if and only if $\alpha \leq \contype(F)$.
	\end{lemma}

	\begin{corollary}\label{cor:connected alpha partition iff alpha is 1 k}
		When $F$ is connected, $F$ has an $\alpha$-partition if and only if \blue{$\alpha = k^1$}, where $k \leq \height(F)+1$.\blue{\footnote{\blue{Recall that $k^1$ is the signature of length $1$ which contains the single value $k$. The starlike tree $\starlike{k^1}$ is the chain on $k+1$ elements.}}}
	\end{corollary}

	Let $F$ be a poset and $\alpha$ be a signature. $F$ is \emph{$\alpha$-connected} if there is no $x \in F$ such that there is an $\alpha$-partition of $\Us x$. By \cref{lem:alpha partition iff alpha leq contype}, this is equivalent to requiring that $\alpha \not\leq \contype(\Us x)$ for each $x \in F$. 

	We can now express the meaning of $\jsl\alpha$ on frames.

	\begin{theorem}\label{thm:Jankov-Fine of starlike iff alpha-connected}
		For $F$ a finite poset and $\alpha$ any signature, $F \vD \jsl\alpha$ if and only if $F$ is $\alpha$-connected.
	\end{theorem}

	To prove this result, we make use of the following slight strengthening of \cref{thm:Jankov-Fine up-reductions}. Let $F$ and $Q$ be finite posets, and assume that $Q$ has root $\bot$. An up-reduction $f \colon F \cra Q$ is \emph{pointed} with \emph{apex} $x \in F$ if we have $\dom(f) = \us x$ and $\inv f\{\bot\} = \{x\}$.

	\begin{lemma}\label{lem:up-reduction to pointed up-reduction}
		If there is an up-reduction $F \cra Q$ then there is a pointed up-reduction $F \cra Q$.
	\end{lemma}

	\begin{proof}
		Take $f \colon F \cra Q$, and choose $x \in \inv f\{\bot\}$ maximal. Then $f |_{\us x}$ is still a p-morphism, and is moreover a pointed up-reduction $F \cra Q$.
	\end{proof}

	\begin{corollary}\label{cor:Jankov-Fine pointed up-reductions}
		Let $F, Q$ be finite posets, with $Q$ rooted. Then $F \vD \chi(Q)$ if and only if there is no pointed up-reduction $F \cra Q$.
	\end{corollary}

	\begin{proof}[Proof of \cref{thm:Jankov-Fine of starlike iff alpha-connected}]
		Assume that $F \nvD \jsl\alpha$. Then by \cref{cor:Jankov-Fine pointed up-reductions} there is a pointed up-reduction $f \colon F \to \starlike\alpha$ with apex $x$. This means that $\inv f[\starlike\alpha \setminus \{\bot\}] = \Us x$. Let $C_j$ be the preimage of the $j$th branch of $\starlike\alpha \setminus \{\bot\}$ under $f$, for each $j \leq \abs\alpha$. Since $f$ is a p-morphism, $C_j$ is upwards-closed. Note that the $C_j$'s are disjoint and hence they form an open partition of $\Us x$. Now, since $C_j$ is the preimage of a chain of length $\alpha(j)$, we can find a chain of the same length inside $C_j$. From this it follows that $C_j$ has height at least $\alpha(j)-1$. But then $(C_j \mid j \leq \abs\alpha)$ is an $\alpha$-partition of $\Us x$, meaning that $F$ is not $\alpha$-connected.

		Conversely, assume that $F$ is not $\alpha$-connected, so that there is $x \in F$ and an $\alpha$-partition $(C_j \mid j \leq k)$ of $\Us x$. For each $C_j$, we have, by definition, that $\height(C_j) \geq \alpha(j)-1$. Hence by \cref{prop:BDn iff no p-morphism to chain n} there is a p-morphism $f_j \colon C_j \to \starlike{\alpha(j)-1}$. Define $f \colon \us x \to \starlike\alpha$ as follows.
		\begin{equation*}
			y \mapsto \left\{
			\begin{array}{ll}
				\bot & \text{if }y = x, \\
				f_j(y) & \text{if }y \in C_j
			\end{array}
			\right.
		\end{equation*}
		Then $f$ is a p-morphism, so an up-reduction $F \cra \starlike\alpha$.
	\end{proof}

	\begin{remark}\label{cor:BDn Jankov-Fine of starlike}
		In particular it follows that $\BD_n = \IPC + \chi(\starlike{n+1})$. This is just \cref{prop:BDn iff no p-morphism to chain n} of course.
	\end{remark}

	The last matter to resolve before moving on to consider the completeness of starlike logics is their number. For this we make use of Higman's Lemma. A \emph{quasi-well-order} is a preorder which is well-founded and has no infinite antichain. Given a preorder $I$, let $I^{<\omega}$ be the set of finite sequences of elements of $I$ ordered by $(x_1, \ldots, x_n) \leq (y_1, \ldots, y_m)$ if and only if there is $f \colon \{1, \ldots n\} \to \{1,\ldots,m\}$ injective such that for each $k \leq n$ we have $x_k \leq y_{f(k)}$.

	\begin{lemma}[Higman's Lemma, \cite{higman52}]\label{lem:Higman}
		If $I$ is a quasi-well-order then so is $I^{<\omega}$.
	\end{lemma}

	\begin{blueblock}
		\begin{proposition}\label{prop:number of starlike logics}\
			\begin{enumerate}[label=(\arabic*)]
				\item\label{item:finite axiomatizable; prop:number of starlike logics}
					Every starlike logic is finitely axiomatizable.
				\item\label{item:number; prop:number of starlike logics}
					There are exactly countably-many starlike logics.
			\end{enumerate}
		\end{proposition}
	\end{blueblock}

	\begin{proof}
		\begin{enumerate}[label=(\arabic*)]
			\item \blue{As every starlike logic is axiomatizable by Jankov formulas of starlike trees, it suffices to show that} there is no infinite antichain of starlike trees with respect to p-morphic reduction. In light of \cref{prop:signature ordering induces p-morphisms}, it therefore suffices to show that there is no infinite antichain of signatures with respect to the ordering defined on them. Now, we can recast signatures as (monotonic decreasing) finite sequences of integers. Indeed, the signature $\alpha$ is determined by the sequence $(\alpha(1), \ldots, \alpha(\abs\alpha))$. In this way, the set of signatures is seen to be a suborder of $\omega^{<\omega}$. Now, $(\omega, \leq)$ is clearly a quasi-well-order, and hence by Higman's Lemma \ref{lem:Higman}, so is $\omega^{<\omega}$. Thus there is no infinite antichain of signatures, as required.
			\item \blue{The result follows from \ref{item:finite axiomatizable; prop:number of starlike logics} as there are only countably many finitely axiomatizable logics.}
		\end{enumerate}
	\end{proof}

\subsection{Starlike completeness}

	The main theorem to be proved in this section is the following.

	\begin{theorem}\label{thm:starlike completeness}
		Every starlike logic is polyhedrally-complete.
	\end{theorem}

	As an immediate consequence, we obtain:

	\begin{corollary}\label{cor:SL ploy-complete}
		Scott's Logic is polyhedrally-complete.
	\end{corollary}

	\begin{remark}
		The starlike logic $\SFL(2 \cdot 1, 1^3)$ is particularly important geometrically. In \cite{convex}, it is shown that this is the logic of all convex polyhedra.
	\end{remark}

	In order to prove \cref{thm:starlike completeness}, we introduce the following new validity concept on frames. Let $F$ be a poset and $\phi$ be a formula. $F$ \emph{nerve-validates} $\phi$, notation $F \vDN \phi$, if for every $k \in \NN$ we have $\N^k(F) \vD \phi$.

	\begin{remark}
		Since for every $G$ we  have the p-morphism $\max \colon \N(G) \to G$ (see \cref{l:geomrea}), by \cref{prop:up-reduction logic containment} this is equivalent to requiring that $\N^k(F) \vD \phi$ for infinitely-many $k \in \NN$.
	\end{remark}

	\begin{lemma}\label{lem:poly-complete iff fmp and p-morphic of nerve-validated}
		A logic $\Lo$ is polyhedrally-complete if and only if it has the finite model property and every rooted finite frame of $\Lo$ is the up-reduction of a poset which nerve-validates $\Lo$.
	\end{lemma}

	\begin{proof}
		Assume that $\Lo$ is polyhedrally-complete. Then by the Nerve Criterion \ref{thm:nerve criterion} it is the logic of a class $\CC$ of finite frames which is closed under $\N$, and so has the fmp. Then by \cref{cor:Logic C every finite rooted frame up-reduction}, every finite rooted frame $F$ of $\Lo$ is the up-reduction of some $F' \in \CC$. Since $\CC \sse \Frames(\Lo)$ and is closed under $\N$, such an $F'$ nerve-validates $\Lo$.

		Conversely, let $\CC$ be the class of all finite rooted frames which nerve-validate $\Lo$. Note that $\CC$ is closed under $\N$. Further, clearly $\Lo \sse \Logic(\CC)$. To see the reverse inclusion, suppose that $\Lo \nvd \phi$. Since $\Lo$ has the fmp, there is $F \in \FramesFinRoot(\Lo)$ such that $F \nvD \phi$. By assumption, $F$ is the up-reduction of $F' \in \CC$. Then by \cref{prop:up-reduction logic containment}, $F' \nvD \phi$, meaning that $\Logic(\CC) \nvd \phi$.
	\end{proof}

	\begin{lemma}\label{lem:starlike logics fmp}
		Every starlike logic has the finite model property.
	\end{lemma}

	\begin{proof}
		In \cite[Corollary~0.11]{zakharyaschev93}, Zakharyaschev shows that every logic axiomatised by the Jankov-Fine formulas of trees has the finite model property.
	\end{proof}

	\begin{blueblock}
		Now, as every finitely axiomatizable logic with the finite model property is decidable we obtain from \cref{prop:number of starlike logics}\ref{item:finite axiomatizable; prop:number of starlike logics} and \cref{lem:starlike logics fmp} the following.
	
		\begin{corollary}
			Every starlike logic is decidable.
		\end{corollary}
	\end{blueblock}

	With \cref{lem:starlike logics fmp}, we can now use \cref{lem:poly-complete iff fmp and p-morphic of nerve-validated} to produce a proof of \cref{thm:starlike completeness}. Given a rooted finite frame $F$ of $\SFL(\Lam)$, we proceed as follows.
	\begin{enumerate}[label=(\arabic*)]
		\item We examine what it means for a frame to nerve-validate $\jsl\alpha$.
		\item We see that it can be assumed that $F$ is \emph{graded} (a structural property of posets defined below).
		\item Using this additional structure, we construct a frame $F'$ and the p-morphism $F' \to F$, with the property that $F' \vDN \SFL(\Lam)$
	\end{enumerate}

		The reader will have noticed that the difork $\starlike{1^2}$ is omitted from the definition of a starlike logic, and consequently from  \cref{thm:starlike completeness}. In fact, polyhedral semantics is quite fond of this tree: when we take it as a forbidden configuration, the resulting landscape of polyhedrally-complete logics is as sparse as possible, as is shown below. 	

	\begin{proposition}\label{prop:poly-complete extension of difork implies CPC}
		Let $\Lo$ be a polyhedrally-complete logic containing $\SFL(1^2)$. Then $\Lo = \CPC$, the maximum logic.
	\end{proposition}

	\begin{proof}
		Suppose for a contradiction that $\Lo$ is a polyhedrally-complete logic containing $\SFL(1^2)$ other than $\CPC$. By the Nerve Criterion \ref{thm:nerve criterion}, $\Lo = \Logic(\CC)$ where $\CC$ is a class of finite posets closed under $\N$. Since $\Lo \neq \CPC$, there must be $F \in \CC$ with $\height(F) \geq 1$. This means that $F$ has a chain $x_0 < x_1$. As in the proof of \cref{thm:stable logics poly incomplete}, we may assume that $x_1$ is a top element of $F$ and that $x_0$ is an immediate predecessor of $x_1$. Take $X$ a maximal chain in $\Ds{x_0}$. Then, as in that proof, we obtain the following structure lying at the top of $\N(F)$.
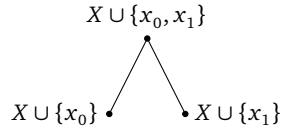
		\begin{equation*}
			\begin{tikzpicture}
				\graph[poset] {
					a[label=left:$X \cup \{x_0\}$] -- b[x=0.5,label=above:{$X \cup \{x_0,x_1\}$}];
					c[label=right:$X \cup \{x_1\}$] -- b
				};
			\end{tikzpicture}
		\end{equation*}
		Applying the nerve once more, we obtain the following structure at the top of $\N^2(F)$.
		\begin{equation*}
			\begin{tikzpicture}
				\graph[poset] {
					a1 -- b1[x=0.5];
					a2[label=below:$Z$] -- {b1,b2[x=0.5]};
					a3 -- b2;
				};
			\end{tikzpicture}
		\end{equation*}
		Since $\CC$ is closed under $\N$, we get that $\N^2(F) \in \Frames(\Lo)$. But $\us Z$ maps p-morph\-ically onto $\starlike{1^2}$, contradicting that $\Lo \vd \chi(\starlike{1^2})$. \contradiction
	\end{proof}

	We now proceed with the proof of \cref{thm:starlike completeness}.

\subsection{Nerve-validation}

	While validating $\jsl\alpha$ corresponds to $\alpha$-connectedness (as shown in \cref{thm:Jankov-Fine of starlike iff alpha-connected}), \emph{nerve}-validating $\jsl\alpha$ corresponds to $\alpha$-\emph{nerve}-connectedness. Let $F$ be a poset and $x<y$ in $F$. The \emph{diamond} and \emph{strict diamond} of $x$ and $y$ are defined, respectively:
	\begin{gather*}
		\uds{x,y} \defeq \us x \cap \ds y \\
		\Uds{x,y} \defeq \uds{x,y} \setminus \{x,y\}
	\end{gather*}

	A poset $F$ is \emph{$\alpha$-diamond-connected} if there are no $x < y$ in $F$ such that there is an $\alpha$-partition of $\Uds{x,y}$. The poset $F$ is \emph{$\alpha$-nerve-connected} if it is $\alpha$-connected and $\alpha$-diamond-connected.

	With a slight conceptual change, $\alpha$-connectedness and $\alpha$-diamond-connectedness can be harmonised as follows. For any poset $F$, we take a new element $\infty$, and let $\check F \defeq F \cup \{\infty\}$, where $\infty$ lies above every element of $F$. Then $F$ is $\alpha$-nerve-connected if and only if there are no $x < y$ in $\check F$ for which there is an $\alpha$-partition of $\Uds{x,y}$.

	The following result shows that $\alpha$-nerve-connectedness is exactly the notion we want.

	\begin{theorem}\label{thm:nerve-validates Jankov of starlike iff alpha-nerve-connected}
		Let $F$ be a finite poset and take $\alpha \in \Starlikes$. Then $F \vDN \jsl\alpha$ if and only if $F$ is $\alpha$-nerve-connected.
	\end{theorem}

	\begin{proof}
		Assume that $F$ is not $\alpha$-nerve-connected with the aim of showing $F \nvDN \chi(\starlike \alpha)$. Choose $x<y$ in $\check F$ such that $\Uds{x,y}$ has an $\alpha$-partition. That is, there is an open partition $(C_j \mid j \leq \abs\alpha)$ of $\Uds {x,y}$ such that $\height(C_j) = \alpha(j)$. 
		\begin{blueblock}
			Choose a chain $X \sse F$ such that:
			\begin{enumerate}[label=(\roman*)]
				\item $x,y \in X \cup \{\infty\}$, and
				\item $X \cap \Uds{x,y} = \es$,
			\end{enumerate}
			which is moreover maximal with respect to these requirements.
		\end{blueblock} 
		We will show that $\Us X^{\N(F)}$ has an $\alpha$-partition. Note that by maximality of $X$, elements $Y \in \Us X^{\N(F)}$ are determined by their intersection $Y \cap \Uds{x,y}$. For $j \leq \abs\alpha$, let:
		\begin{equation*}
			\wh C_j \defeq \{Y \in \Us X^{\N(F)} \mid Y \cap C_j \neq \es\}
		\end{equation*}
		Take $j,l \leq \abs\alpha$ distinct. Since both $C_j$ and $C_l$ are upwards- and downwards-closed in $\Uds{x,y}$, there is no chain $Y \in \Us X^{\N(F)}$ such that $Y \cap C_j \neq \es$ and $Y \cap C_l \neq \es$. This means that:
		\begin{enumerate}[label=(\arabic*)]
			\item $\wh C_j$ and $\wh C_l$ are disjoint.
			\item For any $Y \in \Us X^{\N(F)}$ we have $Y \in \wh C_j$ if and only if $Y \cap \Uds{x,y} \sse C_j$. Hence each $\wh C_j$ is upwards- and downwards-closed in $\Us X^{\N(F)}$.
		\end{enumerate}
		Furthermore, since $(C_j \mid j \leq \abs\alpha)$ covers $\Uds {x,y}$, we get that $(\wh C_j \mid j \leq \abs\alpha)$ covers $\Us X^{\N(F)}$. Finally, any maximal chain in $\wh C_j$ is a sequence of chains $Y_0 \subset \cdots \subset Y_l$ such that $\abs{Y_{i+1} \setminus Y_i} = 1$; this then corresponds to a maximal chain in $C_j$. Therefore:
		\begin{equation*}
			\height(\wh C_j) = \height(C_j)
		\end{equation*}
		Ergo $(\wh C_j \mid j \leq \abs\alpha)$ is an $\alpha$-partition of $\Us X^{\N(F)}$, meaning that $\N(F)$ is not $\alpha$-connected. Then, by \cref{thm:Jankov-Fine of starlike iff alpha-connected}, $\N(F) \nvD \jsl\alpha$, hence by definition $F \nvDN \jsl\alpha$.

		For the converse direction, we will show that if $F$ is $\alpha$-nerve-connected, then so is $\N(F)$, which will give the result by induction (note that $\alpha$-nerve-connectedness is stronger than $\alpha$-connectedness, and hence by \cref{thm:Jankov-Fine of starlike iff alpha-connected} if $\N^k(F)$ is $\alpha$-nerve-connected then $\N^k(F) \vD \jsl\alpha$). So assume that $F$ is $\alpha$-nerve-connected. We will first prove $\alpha$-connectedness. Take $X \in \N(F)$ with the aim of showing that $\Us X^{\N(F)}$ has no $\alpha$-partition. 

		Firstly, assume that $X$ has more than one `gap'; that is, there are distinct $w_1,w_2 \in F \setminus X$ such that $X \cup \{w_1\}$ and $X \cup \{w_2\}$ are still chains, but such that there exists $z \in X$ with $w_1 < z < w_2$. Take $Y,Z \in \Us X^{\N(F)}$. We will use the two gaps to juggle elements between the two sets so as to provide a path $Y \rsa Z$ which never touches $X$ (i.e. lies in $\Us X^{\N(F)}$). For $i \in \{1,2\}$, let $u_i \in X \cap \Ds {w_i}$ be greatest and $v_i \in X \cap \Us {w_i}$ be least. See \cref{fig:X more than one gap; thm:nerve-validates Jankov of starlike iff alpha-nerve-connected} for a representation of the situation. Now, without loss of generality, we may assume that $Y \cap \Uds{u_1,v_1} \neq \es$ (we may add $w_1$ to $Y$, noting that $w_1 \in \Uds{u_1,v_1}$). Similarly, we may assume that $Y \cap \Uds{u_2,v_2} \neq \es$, and likewise for $Z$. We then have the following path in $\Us X^{\N(F)}$ (note that some of the sets along the path may be equal, but in all cases the path is still there):
		\begin{equation*}
			\begin{tikzpicture}[xscale=1.65]
				\node (a) at (0,1) {$Y$};
				\node (b) at (1,0) {$Y \setminus \Uds {u_1,v_1}$};
				\node (c) at (2,1) {$(Y \setminus \Uds {u_1,v_1}) \cup \{w_1\}$};
				\node (d) at (3,0) {$X \cup \{w_1\}$};
				\node (e) at (4,1) {$(Z \setminus \Uds {u_1,v_1}) \cup \{w_1\}$};
				\node (f) at (5,0) {$Z \setminus \Uds {u_1,v_1}$};
				\node (g) at (6,1) {$Z$};
				\graph[use existing nodes] {
					a -- b -- c -- d -- e -- f -- g;
				};
			\end{tikzpicture}
		\end{equation*}
		Here, the gap $\Uds{u_2,v_2}$ is used to ensure that $Y \setminus \Uds {u_1,v_1}$ and $Z \setminus \Uds {u_1,v_1}$ are not equal to $X$, and the fact that we have $v_1 \leq z \leq u_2$ ensures that all these sets are indeed in $\N(F)$. Hence, $\Us X^{\N(F)}$ is path-connected so connected. Therefore, by \cref{cor:connected alpha partition iff alpha is 1 k}, it suffices to show that $\height(\Us X^{\N(F)}) < \height(F)$. But this is immediate from the definition of $\N$.

		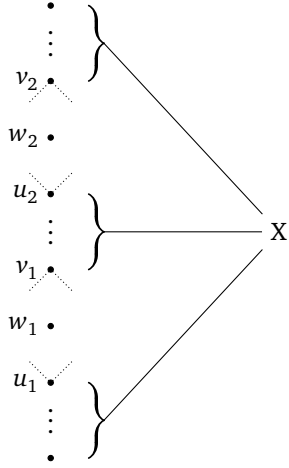
\begin{figure}
			\begin{equation*}
				\begin{tikzpicture}
					\path
						(0,0) node[point] {}
						-- ++(0,0.6) node {$\vdots$}
						-- ++(0,0.4) node[point] [label=left:$u_1$] {}
						edge[draw,densely dotted] +(-0.3,0.3)
						edge[draw,densely dotted] +(0.3,0.3)
						-- +(0,0.75) node[point] [label=left:$w_1$] {}
						-- ++(0,1.5) node[point] [label=left:$v_1$] {}
						edge[draw,densely dotted] +(-0.3,-0.3)
						edge[draw,densely dotted] +(0.3,-0.3)
						-- ++(0,0.6) node {$\vdots$}
						-- ++(0,0.4) node[point] [label=left:$u_2$] {}
						edge[draw,densely dotted] +(-0.3,0.3)
						edge[draw,densely dotted] +(0.3,0.3)
						-- +(0,0.75) node[point] [label=left:$w_2$] {}
						-- ++(0,1.5) node[point] [label=left:$v_2$] {}
						edge[draw,densely dotted] +(-0.3,-0.3)
						edge[draw,densely dotted] +(0.3,-0.3)
						-- ++(0,0.6) node {$\vdots$}
						-- ++(0,0.4) node[point] {};
						\draw[decoration={calligraphic brace,amplitude=0.2cm}, decorate, line width=1.25pt] (0.5,1) -- (0.5,0);
						\draw[decoration={calligraphic brace,amplitude=0.2cm}, decorate, line width=1.25pt] (0.5,3.5) -- (0.5,2.5);
						\draw[decoration={calligraphic brace,amplitude=0.2cm}, decorate, line width=1.25pt] (0.5,6) -- (0.5,5);
						\node at (3,3) (X) {X};
						\draw (X)
							edge (0.7,0.5)
							edge (0.7,3)
							edge (0.7,5.5);
				\end{tikzpicture}
			\end{equation*}
			\caption{The set-up when $X$ has more than one gap}
			\label{fig:X more than one gap; thm:nerve-validates Jankov of starlike iff alpha-nerve-connected}
		\end{figure}

		Hence we may assume that $X$ has exactly one gap (when $X$ has no gaps, $\Us X^{\N(F)} = \es$). This means that there are $x,y \in X$ with $x<y$ such that $X \cap \Uds{x,y} = \es$ and $X$ is maximal outside of $\Uds{x,y}$. As before then, elements $Y \in \Us X^{\N(F)}$ are determined by their intersection $Y \cap \Uds{x,y}$. Suppose that $\Us X^{\N(F)}$ has an $\alpha$-partition $(\wh C_j \mid j \leq \abs\alpha)$. For each $j \leq \abs\alpha$, let:
		\begin{equation*}
			C_j \defeq \bigcup \wh C_j \cap \Uds{x,y}
		\end{equation*}
		Note that $\bigcup_{j \leq \abs\alpha}C_j = \Uds{x,y}$. For each $j \leq \abs\alpha$, since $\wh C_j$ is downwards-closed, we have that, for $z \in \Uds{x,y}$:
		\begin{equation*}
			z \in C_j \quad\Lra\quad \exists Y \in \wh C_j \colon z \in Y \quad\Lra\quad X \cup \{z\} \in \wh C_j
		\end{equation*}
		This means in particular that the $C_j$'s are pairwise disjoint. Further, if $z \in C_j$ and $w \in \Uds{x,y}$ with $w < z$, then $X \cup \{w,z\}$ is a chain, and so as $\wh C_j$ is upwards-closed, we have $X \cup \{w,z\} \in \wh C_j$, meaning that $w \in C_j$; similarly when $w > z$. Whence each $C_j$ is upwards- and downwards-closed. Finally, as above, maximal chains in $\wh C_j$ correspond to maximal chains in $C_j$ of the same length, whence:
		\begin{equation*}
			\height(\wh C_j) = \height(C_j)
		\end{equation*}
		But then $(C_j \mid j \leq \abs\alpha)$ is an $\alpha$-partition of $\Uds{x,y}$, contradicting the fact that $F$ is $\alpha$-nerve-connected. \contradiction

		This shows that $\N(F)$ is $\alpha$-connected. What about $\alpha$-diamond-connectedness? In fact we can show this without using any assumptions on $F$. Take $X,Y \in \N(F)$ with $X \subset Y$. We will show that $\Uds{X,Y}^{\N(F)}$ has no $\alpha$-partition. We may assume that $\abs{Y \setminus X} \geq 2$, otherwise $\Uds{X,Y}^{\N(F)} = \es$. Note that this means in particular that $\alpha > 1$, since $F$ is $\alpha$-connected. If $\abs{Y \setminus X} = 2$, then $\Uds{X,Y}^{\N(F)}$ is the antichain on two elements, which, since $\alpha \neq 1^2$ by assumption, has no $\alpha$-partition. So assume that $\abs{Y \setminus X} \geq 3$; we will show that in fact $\Uds{X,Y}^{\N(F)}$ is connected. Take distinct $Z,W \in \Uds{X,Y}^{\N(F)}$. Choose $z \in Z \setminus X$ and $w \in W \setminus X$. Since $\abs{Y \setminus X} \geq 3$, we have that $X \cup \{z,w\} \in \Uds{X,Y}^{\N(F)}$. Hence the following is a path in $\Uds{X,Y}^{\N(F)}$:
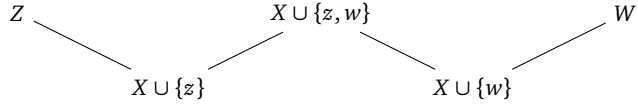
		\begin{equation*}
			\begin{tikzpicture}[xscale=2]
				\node (a) at (0,1) {$Z$};
				\node (b) at (1,0) {$X \cup \{z\}$};
				\node (c) at (2,1) {$X \cup \{z,w\}$};
				\node (d) at (3,0) {$X \cup \{w\}$};
				\node (e) at (4,1) {$W$};
				\graph[use existing nodes] {
					a -- b -- c -- d -- e;
				};
			\end{tikzpicture}
		\end{equation*}
		Therefore, $\Uds{X,Y}^{\N(F)}$ is connected. Finally, note that:
		\begin{equation*}
			\height(\Uds{X,Y}^{\N(F)}) \leq \height(\N(F)) = \height(F) \qedhere
		\end{equation*}
	\end{proof}

	\begin{remark}\label{rem:vDN jsl alpha iff NF vD jsl alpha}
		Note that the proof shows an interesting property of the formulas $\jsl\alpha$: we have $F \vDN \jsl\alpha$ if and only if $\N(F) \vD \jsl\alpha$. This is not true in general. For example, formulas expressing bounded width can take many iterations of the nerve construction to become falsified. 
	\end{remark}

\subsection{Graded posets}

	The next step is to show that we can put $F \in \FramesFinRoot(\SFL(\Lam))$ into a special form. The following definition comes from combinatorics (see e.g. \cite[p.~99]{stanley1997}).

	\begin{definition}[Graded poset]
		A \emph{rank function} on a poset $F$ is a map $\rho \colon F \to \NN$ such that:
		\begin{enumerate}[label=(\roman*)]
			\item whenever $x$ is minimal in $F$, we have $\rho(x)=0$,
			\item whenever $y$ is the immediate successor of $x$, we have $\rho(y)=\rho(x)+1$.
		\end{enumerate}
		If $F$ is non-empty and has a rank function, then it is \emph{graded}.
	\end{definition}

	The notion of gradedness has a strong visual connection. When a poset is graded, we can draw it out in well-defined layers such that any element's immediate successors lie entirely in the next layer up. 

	\begin{lemma}\label{lem:properties of gradedness}
		Let $F$ be a finite poset.
		\begin{enumerate}[label=(\arabic*)]
			\item\label{item:maximal chains; lem:properties of gradedness}
				$F$ is graded if and only if for every $x \in F$, all maximal chains in $\ds x$ have the same length.
			\item\label{item:height; lem:properties of gradedness}
				When $F$ is graded, $\rho(x) = \height(x)$ for every $x \in F$, and $\height(F) = \max\rho[F]$.
			\item\label{item:unique; lem:properties of gradedness}
				Rank functions, when they exist, are unique.
		\end{enumerate}
	\end{lemma}

	\begin{proof}
		\begin{enumerate}[label=(\arabic*)]
			\item See \cite[p.~99]{stanley1997}. Assume that $F$ is graded, and take $X$ a maximal chain in $\ds x$ for some $x \in F$. Let $k = \rho(x)$. We will show that $\abs X = k+1$. Since $X$ is a chain, the ranks of each of its elements are distinct. Since $X$ is maximal, $x \in X$. Suppose for a contradiction that there is $j < k$ such that there is no $x \in X$ of rank $j$. We may assume that $j$ is minimal with this property. We can't have $j=0$, since otherwise $X$ wouldn't contain any minimal element, so wouldn't be a maximal chain. Hence, there is $y \in X$ with $\rho(y) = j-1$. Let $z$ be next in $X$ after $y$. Then $y$ has an immediate successor $w$ such that $w \leq z$. By definition, $\rho(w) = j$, so $w \notin X$. But $X \cup \{w\}$ is a chain, contradicting the maximality of $X$. \contradiction Therefore, $\abs X = k+1$.

			Conversely, define $\rho \colon F \to \NN$ by:
			\begin{equation*}
				x \mapsto \height(x)
			\end{equation*}
			Let us check that $\rho$ is a rank function.
			\begin{enumerate*}[label=(\roman*),mode=unboxed]
				\item Clearly, when $x$ is minimal, $\rho(x) = 0$.
				\item Suppose for a contradiction that there are $x, y \in F$, with $y$ an immediate successor of $x$, such that $\rho(y) \neq \rho(x) + 1$. First, by definition, $\rho(y) > \rho(x)$, so we must have $\rho(y) > \rho(x) + 1$. Choose maximal chains $X \sse \ds x$, $Y \sse \ds y$. Note that by assumption:
				\begin{equation*}
					\abs Y > \abs X + 1
				\end{equation*}
				But now, since $y$ is an immediate successor of $x$, both $X \cup \{y\}$ and $Y$ are maximal chains in $\ds y$ of different heights. \contradiction
			\end{enumerate*}

			\item This follows from the proof of \ref{item:maximal chains; lem:properties of gradedness}.

			\item This follows from \ref{item:height; lem:properties of gradedness}.\qedhere 
		\end{enumerate}
	\end{proof}

	\begin{corollary}\label{cor:trees and nerves graded}
		\begin{enumerate}[label=(\arabic*)]
			\item\label{item:trees; cor:trees and nerves graded}
				Every tree is graded.
			\item\label{item:nerves; cor:trees and nerves graded}
				For any finite poset $F$, its nerve $\N(F)$ is graded, with rank function given by $\rho(X) = \abs{X}-1$.
		\end{enumerate}
	\end{corollary}

	\begin{proof}
		For \ref{item:nerves; cor:trees and nerves graded}, note that for any $X \in \N(F)$ we have $\height(X) = \abs{X}-1$.
	\end{proof}

	\begin{blueblock}
		What we will show in the proceeding two subsections is that any frame $\SFL(\Lam)$ can be assumed to be graded. In other words, we prove the following `gradification' theorem.
	
		\begin{theorem}\label{thm:gradification}
			Take $\Lam \sse \Starlikes$ and let $F$ be a finite rooted poset such that $F \vD \SFL(\Lam)$. Then there is a finite graded rooted poset $F'$ and a p-morphism $f \colon F' \to F$ such that $F' \vD \SFL(\Lam)$.
		\end{theorem}
	
		The proof of the theorem works differently depending on whether we have Scott's tree $\starlike{2\cdot 1}$ present. \Cref{thm:gradification with Scott} deals with the case $2 \cdot 1 \in \Lam$, while \cref{thm:gradification without Scott} deals with the case $2 \cdot 1 \notin \Lam$.
	\end{blueblock}

\subsection{Gradification in the presence of Scott's tree}

	\begin{blueblock}
		Let us first consider the gradification theorem in the case $2 \cdot 1 \in \Lam$.

		\begin{theorem}\label{thm:gradification with Scott}
			Let $\Lam \sse \Starlikes$ be such that $2 \cdot 1 \in \Lam$. Let $F$ be a finite rooted poset such that $F \vD \SFL(\Lam)$. Then there is a finite graded rooted poset $F'$ and a p-morphism $f \colon F' \to F$ such that $F' \vD \SFL(\Lam)$.
		\end{theorem}

		To begin with, the following lemmas show us that this case is not too complicated.
	\end{blueblock}


	\begin{lemma}\label{lem:SFL with Scott equal SFL with minimal 1 k}
		Take $\Lam \sse \Starlikes$ such that $2\cdot1 \in \Lam$ but $n \notin \Lam$ for any $n \in \NN$.
		\begin{enumerate}[label=(\arabic*)]
			\item\label{item:no 1; lem:SFL with Scott equal SFL with minimal 1}
				If there is no $k \in \NNp$ such that $1^k \in \Lam$, then $\SFL(\Lam) = \SFL(2 \cdot 1)$.
			\item\label{item:1; lem:SFL with Scott equal SFL with minimal 1}
				Otherwise, let $k \in \NNp$ be minimal such that $1^k \in \Lam$. Then $\SFL(\Lam) = \SFL(2 \cdot 1,1^k)$.
		\end{enumerate}
	\end{lemma}

	\begin{proof}
		\begin{enumerate}[label=(\arabic*)]
			\item Take $\alpha \in \Lam$. Then by assumption $\alpha(1) \geq {2}$, hence, as $\alpha \neq n$, we have ${2 \cdot 1} \leq \alpha$. Then by \cref{prop:signature ordering induces p-morphisms} there is a p-morphism $\starlike\alpha \to \starlike {2 \cdot 1}$. Hence by the semantic meaning of Jankov-Fine formulas, \cref{thm:Jankov-Fine up-reductions}, we have that any frame validating $\chi(\starlike {2 \cdot 1})$ will also validate $\chi(\starlike\alpha)$. This means that $\SFL(\Lam) \sse \SFL(2 \cdot 1)$. The converse direction is immediate.
			\item Take $\alpha \in \Lam$. If $\alpha(1) \geq 2$ then by \cref{prop:signature ordering induces p-morphisms} there is a p-morphism $\starlike\alpha \to \starlike {2 \cdot 1}$. Assume that $\alpha(1) \not\geq 2$. Since $\alpha \neq \ep$, we have $\alpha(1)=1$, meaning that $\alpha=1^l$ for some $l \in \NNp$. By assumption $k \leq l$. But then $1^k \leq \alpha$, giving that there is a p-morphism $\starlike\alpha \to \starlike {1^k}$. It follows that for any $\alpha \in \Lam$, $\starlike\alpha$ up-reduces to either $\starlike {2 \cdot 1}$ or $\starlike {1^k}$. By \cref{thm:Jankov-Fine up-reductions}, any frame validating $\chi(\starlike {2 \cdot 1})$ and $\chi(\starlike {1^k})$ will also validate $\chi(\starlike {\alpha})$. This implies that $\SFL(\Lam) \subseteq \SFL(2 \cdot 1,1^k)$. The converse direction is obvious. \qedhere
		\end{enumerate}
	\end{proof}

	\begin{corollary}\label{cor:SFL with Scott equal SFL with minimal 1 k plus BDn}
		Take $\Lam \sse \Starlikes$ such that $2\cdot1 \in \Lam$ and there is $n \in \NN$ with $n \in \Lam$; assume that $n$ is the minimal such natural number.
		\begin{enumerate}[label=(\arabic*)]
			\item\label{item:no 1; cor:SFL with Scott equal SFL with minimal 1 k plus BDn}
				If there is no $k \in \NNp$ such that $1^k \in \Lam$, then $\SFL(\Lam) = \SFL(n, 2 \cdot 1)$.
			\item\label{item:1; cor:SFL with Scott equal SFL with minimal 1 k plus BDn}
				Otherwise, let $k \in \NNp$ be minimal with $1^k \in \Lam$. Then $\SFL(\Lam) = \SFL(n, 2 \cdot 1,1^k)$.
		\end{enumerate}
	\end{corollary}

	\begin{proof}
		This follows from \cref{lem:SFL with Scott equal SFL with minimal 1 k} and the fact that when $n_1 < n_2$ every frame validating $\jsl{n_1}$ also validates $\jsl{n_2}$.
	\end{proof}

	Using this, the `meaning' of $\SFL(\Lam)$ can be expressed relatively simply. Note that this meaning is expressed in terms of the depth of elements $x \in F$. Up until this point we have mainly been concerned with the height of elements.

	\begin{lemma}\label{lem:SFL meaning in presence of Scott}
		Take $\Lam \sse \Starlikes$ such that $2\cdot1 \in \Lam$, and let $F$ be a finite poset. Let $n \in \NN$ be minimal such that $n \in \Lam$, or $\infty$ if no such signature is present. Similarly, let $k \in \NNp$ be minimal with $1^k \in \Lam$, or $\infty$. Then $F \vD \SFL(\Lam)$ if and only if the following three conditions are satisfied for every $x \in F$.
		\begin{enumerate}[label=(\roman*)]
			\item\label{item:height n; lem:SFL meaning in presence of Scott}
			We have $\height(F) < n$.
			\item\label{item:depth 1; lem:SFL meaning in presence of Scott}
			Whenever $\depth(x) = 1$, we have $\abs{\Us x} < k$.
			\item\label{item:depth gt1; lem:SFL meaning in presence of Scott}
			Whenever $\depth(x) > 1$, the set $\Us x$ is connected.
		\end{enumerate}
	\end{lemma}

	\begin{proof}
		By \cref{cor:SFL with Scott equal SFL with minimal 1 k plus BDn} and the fact that $F \vD \jsl{n}$ if and only if $\height(F) \leq n-1$, it suffices to treat the case $n=\infty$. Now by \cref{lem:SFL with Scott equal SFL with minimal 1 k}, $\SFL(\Lam) = \SFL(2 \cdot 1, 1^k)$ when $k < \infty$, and $\SFL(\Lam) = \SFL(2 \cdot 1)$ otherwise. 

		Assume that $F \vD \SFL(\Lam)$. 
		\begin{itemize*}
			\item[\ref{item:depth 1; lem:SFL meaning in presence of Scott}] In the case $k < \infty$, take $x \in F$ with $\depth(x) = 1$. Note that $\Us x$ is an antichain, so $(\{y\} \mid y \in \Us x)$ is an open partition of $\Us x$. Since $x \vD \chi(\starlike{1^k})$, by \cref{lem:alpha partition iff alpha leq contype} and \cref{thm:Jankov-Fine of starlike iff alpha-connected} we must have $\abs{\Us x} < k$. 
			\item[\ref{item:depth gt1; lem:SFL meaning in presence of Scott}] Now take $x \in F$ with $\depth(x) > 1$, and suppose for a contradiction that $\Us x$ is disconnected. Then we can partition $\Us x$ into disjoint upwards-closed sets $U,V$. Since $\depth(x) > 1$, one of $U$ and $V$ (say $U$) must have height at least $1$. But then $(U,V)$ is a $(2 \cdot 1)$-partition of $\Us x$, contradicting that $F \vD \chi(\starlike{2 \cdot 1})$ by \cref{thm:Jankov-Fine of starlike iff alpha-connected}. \contradiction
		\end{itemize*}

		Conversely, assume that $F \nvD \SFL(\Lam)$ We will show that one of \ref{item:depth 1; lem:SFL meaning in presence of Scott} and \ref{item:depth gt1; lem:SFL meaning in presence of Scott} is violated. If $F \nvD \chi(\starlike{2 \cdot 1})$, then by \cref{thm:Jankov-Fine of starlike iff alpha-connected} there is $x \in F$ and a $(2 \cdot 1)$-partition $(U,V)$ of $\Us x$. But then $\height(U) \geq 1$, meaning that $\depth(x) > 1$, and furthermore $\Us x$ is disconnected, violating \ref{item:depth gt1; lem:SFL meaning in presence of Scott}. So let us assume that $k<\infty$, that $F \vD \chi(\starlike{2 \cdot 1})$ but that $F \nvD \chi(\starlike{1^k})$. Again, we get $x \in F$ and a $1^k$-partition $(C_1, \ldots, C_k)$ of $\Us x$. We must have that $\height(C_1) = 0$, otherwise $(C_1,C_2 \cup \cdots \cup C_k)$ is a $(2 \cdot 1)$-partition of $\Us x$. Similarly $\height(C_i) = 0$ for every $i \leq k$. This means that $\depth(x) = 1$, and that $\abs{\Us x} \geq k$, violating \ref{item:depth 1; lem:SFL meaning in presence of Scott}.
	\end{proof}

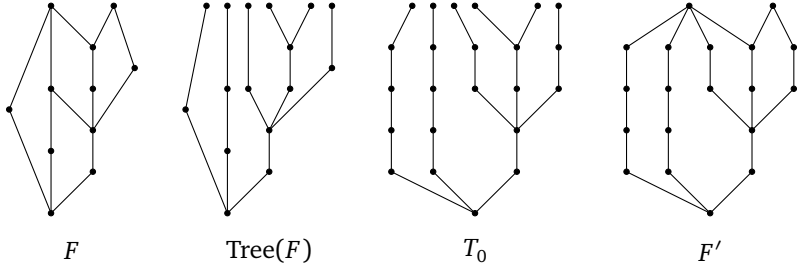
\begin{figure}
		\begin{equation*}
			\begin{tikzpicture}[scale=0.55]
				\begin{scope}
					\graph[poset=0.55] {
						a0[x=1cm] -- {
							b0[y=1.5cm] -- c0[x=1cm,y=3cm],
							b1[y=0.5cm] -- c1[y=1cm] -- c0,
							b2 -- c2 -- {
								c1,
								d1 -- e1 -- {c0, f1[x=0.5cm]},
								d2[y=0.5cm] -- f1
							}
						};
					};
					\node at (1.5,-0.9) {$F$};
				\end{scope}
				\begin{scope}[xshift=120]
					\graph[poset=0.55] {
						a0[x=1cm] -- {
							b0[y=1.5cm] -- c0[x=0.5cm,y=3cm],
							b1[y=0.5cm] -- c1[y=1cm] -- d1[y=2cm],
							b2 -- c2 -- {
								d2[x=-0.5cm] -- e1[x=-0.5cm,y=1cm],
								d3[x=-0.5cm] -- e2[x=-0.5cm] -- {f1[x=-1cm], f2[x=-1cm]},
								d4[x=-1.5cm,y=0.5cm] -- e3[x=-1.5cm,y=1cm]
							}
						};
					};
					\node at (2,-0.9) {$\Tree(F)$};
				\end{scope}
				\begin{scope}[xshift=260]
					\graph[poset=0.55] {
						a0[x=2cm] -- {
							b0 -- c0 -- d0 -- e0 -- f0[x=0.5cm],
							b1 -- c1 -- d1 -- e1 -- f1,
							b2[x=1cm] -- c2[x=1cm] -- {
								d2 -- e2 -- f2[x=-0.5cm],
								d3 -- e3 -- {f3[x=-1cm], f4[x=-0.5cm]},
								d4[x=-1cm] -- e4[x=-1cm] -- f5[x=-1cm]
							}
						};
					};
					\node at (2,-0.9) {$T_0$};
				\end{scope}
				\begin{scope}[xshift=420]
					\graph[poset=0.55] {
						a0[x=2cm] -- {
							b0 -- c0 -- d0 -- e0 -- f0[x=1.5cm],
							b1 -- c1 -- d1 -- e1 -- f0,
							b2[x=1cm] -- c2[x=1cm] -- {
								d2 -- e2 -- f0,
								d3 -- e3 -- {f0, f1[x=0.5cm]},
								d4 -- e4 -- f1
							}
						};
					};
					\node at (2,-0.9) {$F'$};
				\end{scope}
			\end{tikzpicture}
		\end{equation*}
		\caption{An example of gradification in the presence of Scott's tree}
		\label{fig:gradification with Scott example}
	\end{figure}

	\blue{Let us turn now to the proof of \cref{thm:gradification with Scott}. We first} outline the construction before coming to the full proof.
	\begin{itemize}
		\item We first split $F$ up into its tree unravelling $\Tree(F)$ (defined below).
		\item We then lengthen branches so that every top element has the same height.
		\item Lastly, we join top nodes of this tree in order to recover any $\alpha$-connectedness that we lost.
	\end{itemize}
	See \cref{fig:gradification with Scott example} for an example of this process.

	Given any finite, rooted poset $F$, its \emph{tree unravelling} $\Tree(F)$ is the set of chains $X$ in $F$ such that $X$ is maximal in $\ds{\max (X)}$, ordered by subset inclusion. Define the function $\last \colon \Tree(F) \to F$ by:
	\begin{equation*}
		X \mapsto \max (X)
	\end{equation*}
	Then $\Tree(F)$ is a tree and $\last$ is a p-morphism (see \cite[Theorem~2.19, p.~32]{chagrovzakharyaschev1997}).

	\blue{We make use of the following abbreviations. For any poset $F$, the set of top elements (i.e.\@ elements of depth $0$) in $F$ is denoted by $\Top(F)$; let $\Trunk(F) \defeq F \setminus \Top(F)$.}

	\begin{proof}[Proof of \cref{thm:gradification with Scott}]
		\newcommand{\newel}[2]{#1^*(#2)}

		Let $n \defeq \height(F)$. We may assume $\ep \notin \Lam$. If $2 \in \Lam$, then by \cref{cor:BDn Jankov-Fine of starlike}, $n \leq 1$, so $F$ is already graded. So assume that $2 \notin \Lam$.

		Start with the tree unravelling $T=\Tree(F)$ of $F$. Form a new tree $T_0$ by replacing each top node $t \in \Top(T)$ with a chain of new elements $\newel{t}{0}, \ldots, \newel{t}{m_t}$, where $m_t = n-\height(t)$. The relations between these new elements and the rest of $T$ is as follows:
		\begin{gather*}
			\newel{t}{0} < \cdots < \newel{t}{{m_t}}, \\
			x < \newel{t}{0} \quad \Lra \quad x < t \qquad \forall x \in T
		\end{gather*}
		Note that in $T_0$ all branches have the same length $n+1$. Define the p-morphism $g \colon T_0 \to T$ by:
		\begin{equation*}
			x \mapsto \left\{
			\begin{array}{ll}
				x & \text{ if }x \in \Trunk(T), \\
				\last(t) & \text{ if }x=\newel{t}{i}\text{ for some }t \in \Top(T)\text{ and }i \leq m_t
			\end{array}
			\right.
		\end{equation*}
		Form $F'$ from $T_0$ by identifying, for top nodes $t,s \in \Top(T)$, the elements $\newel{t}{m_t}$ and $\newel{s}{m_s}$ whenever $\last(t)=\last(s)$. That is, let $F' \defeq T_0/\W$, where: 
		\begin{equation*}
			\W \defeq \{\{\newel{t}{m_t} \mid \last(t)=u\} \mid u \in \Top(F)\}
		\end{equation*}

		Note that we have a p-morphism $f = \last \circ g \circ q_\W \colon F' \to F$. Furthermore, $F$ is clearly finite and rooted. As to gradedness, take $x \in F'$ with the aim of showing that all maximal chains in $\ds x$ are of the same length, utilising \cref{lem:properties of gradedness}. If $x \in \Trunk(F')$, then $\ds x^{F'}$ is a linear order. So assume that $x \in \Top(F')$. Then any maximal chain $X$ in $\ds x$ corresponds to a branch of $T_0$, and therefore has length $n+1$.

		Let us now use \cref{lem:SFL meaning in presence of Scott} to verify that our construction preserves $\alpha$-connected\-ness for $\alpha \in \Lam$ and complete the proof. Let $k \in \NNp$ be minimal such that $1^k \in \Lam$, or $\infty$ if no such signature is present. For $u \in \Top(F)$ let $\wh u$ be the equivalence class of those elements $\newel{t}{m_t}$ such that $\last(t) = u$. Note that by construction, for $x \in \Trunk(T)$ and $u \in \Top(F)$:
		\begin{equation}\label{eq:x lt hat u iff; proof:gradification with Scott}
			x < \wh u \quad\Lra\quad \last(x) < u \tag{$\star$}
		\end{equation}
		We need to check the three conditions of \cref{lem:SFL meaning in presence of Scott}.
		\begin{itemize}
			\item[\ref{item:height n; lem:SFL meaning in presence of Scott}] 
				Note that $\height(F') = \height(F)$.
			\item[\ref{item:depth 1; lem:SFL meaning in presence of Scott}] 
				For any $x \in F'$ with $\depth(x)=1$, either $x \in \Trunk(T)$ or $x=\newel{t}{n_t-1}$ for some top node $t \in T$. In the former case, the fact that $\abs{\Us x} \leq k$ follows from \eqref{eq:x lt hat u iff; proof:gradification with Scott} and the fact that $\abs{\Us {\last(x)}^F} \leq k$. In the latter case we have $\Us x = \left\{\wh{\last(t)}\right\}$.
			\item[\ref{item:depth gt1; lem:SFL meaning in presence of Scott}] 
				Similarly, for any $x \in F'$ with $\depth(x) > 1$, either $x \in \Trunk(T)$ or $x=\newel{t}{r}$ for some top node $t \in T$ and $r < n_t-1$. In the latter case, $\Us x$ is a chain, so connected. For the former case, it suffices to show that any two top elements $\wh u, \wh v \in \Us x$ are connected by a path in $\Us x$. Note that $\depth(\last(x))^F > 1$. Now, since $F \vD \chi(\starlike{2 \cdot 1})$, by \cref{lem:SFL meaning in presence of Scott} there is a path $u \rsa v$ in $\Us{\last(x)}^F$. We may assume that this path is of form given in \cref{fig:form of paths in F and F prime; proof:gradification with Scott} (a), where $w_0, \ldots, w_k$ are top nodes in $F$. Using \eqref{eq:x lt hat u iff; proof:gradification with Scott}, this path then translates into a path $\wh u \rsa \wh v$ as in \cref{fig:form of paths in F and F prime; proof:gradification with Scott} (b), where $y_i \in \last^{-1}\{a_i\} \cap \Us x$ for each $i$. \qedhere
		\end{itemize}

		\begin{figure}
			\begin{equation*}
				\begin{tikzpicture}[scale=0.7]
					\begin{scope}
						\node at (-1.5,-0.5) {(a)};
						\graph[grow right, math nodes] {
							w_0
								-- a_0[y=-1cm]
								-- w_1
								-- a_1[y=-1cm]
								-!- e0[x=0.5cm,y=-0.25cm,as={$\cdots$}]
								-!- e1[empty nodes]
								-!- a_{k-2}[y=-1cm]
								-- w_{k-1}
								-- a_{k-1}[y=-1cm]
								-- w_k
								;
						};
						\draw (a_1) -- ++(0.5,0.5);
						\draw (a_{k-2}) -- ++(-0.5,0.5);
					\end{scope}
					\begin{scope}[yshift=-60]
						\node at (-1.5,-0.5) {(b)};
						\graph[grow right, math nodes] {
							"\wh{w_0}"
								-- y_0[y=-1cm]
								-- "\wh{w_1}"
								-- y_1[y=-1cm]
								-!- e0[x=0.5cm,y=-0.25cm,as={$\cdots$}]
								-!- e1[empty nodes]
								-!- y_{k-2}[y=-1cm]
								-- "\wh{w_{k-1}}"
								-- y_{k-1}[y=-1cm]
								-- "\wh{w_k}"
								;
						};
						\draw (y_1) -- ++(0.5,0.5);
						\draw (y_{k-2}) -- ++(-0.5,0.5);
					\end{scope}
				\end{tikzpicture}
			\end{equation*}
			\caption{The form of the paths in $\protect\Us{\last(x)}^F$ and $\protect\Us x^{F'}$}
			\label{fig:form of paths in F and F prime; proof:gradification with Scott}
		\end{figure}
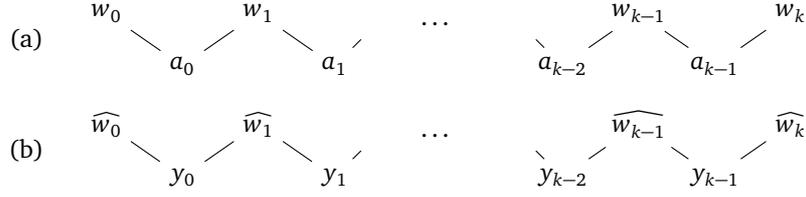
	\end{proof}

\subsection{Gradification without Scott's tree}

	Now that the situation $2 \cdot 1 \in \Lam$ has been dealt with, let us turn to the case $2 \cdot 1 \notin \Lam$.

	\begin{blueblock}
		\begin{theorem}\label{thm:gradification without Scott}
			Let $\Lam \sse \Starlikes$ be such that $2 \cdot 1 \notin \Lam$. Let $F$ be a finite, rooted poset such that $F \vD \SFL(\Lam)$. Then there is a finite, graded, rooted poset $F'$ and a p-morphism $f \colon F' \to F$ such that $F' \vD \SFL(\Lam)$.
		\end{theorem}
	\end{blueblock}
	
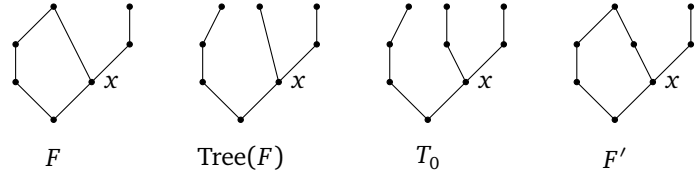
\begin{figure}
		\begin{equation*}
			\begin{tikzpicture}[scale=0.5]
				\begin{scope}
					\node[point] (a) at (1,0) {};
					\node[point] (b1) at (0,1) {};
					\node[point] (b2) at (0,2) {};
					\node[point] (c) at (2,1) [label=right:$x$] {};
					\node[point] (f1) at (3,2) {};
					\node[point] (f2) at (3,3) {};
					\node[point] (d) at (1,3) {};
					\graph[use existing nodes] {
						a -- b1 -- b2 -- d;
						a -- c -- d;
						c -- f1 -- f2;
					};
					\node at (1,-1) {$F$};
				\end{scope}
				\begin{scope}[xshift=140]
					\node[point] (a) at (1,0) {};
					\node[point] (b1) at (0,1) {};
					\node[point] (b2) at (0,2) {};
					\node[point] (c) at (2,1) [label=right:$x$] {};
					\node[point] (f1) at (3,2) {};
					\node[point] (f2) at (3,3) {};
					\node[point] (db) at (0.5,3) {};
					\node[point] (dc) at (1.5,3) {};
					\graph[use existing nodes] {
						a -- b1 -- b2 -- db;
						a -- c -- dc;
						c -- f1 -- f2;
					};
					\node at (1,-1) {$\Tree(F)$};
				\end{scope}
				\begin{scope}[xshift=280]
					\node[point] (a) at (1,0) {};
					\node[point] (b1) at (0,1) {};
					\node[point] (b2) at (0,2) {};
					\node[point] (c) at (2,1) [label=right:$x$] {};
					\node[point] (f1) at (3,2) {};
					\node[point] (f2) at (3,3) {};
					\node[point] (db) at (0.5,3) {};
					\node[point] (dc1) at (1.5,2) {};
					\node[point] (dc2) at (1.5,3) {};
					\graph[use existing nodes] {
						a -- b1 -- b2 -- db;
						a -- c -- dc1 -- dc2;
						c -- f1 -- f2;
					};
					\node at (1,-1) {$T_0$};
				\end{scope}
				\begin{scope}[xshift=420]
					\node[point] (a) at (1,0) {};
					\node[point] (b1) at (0,1) {};
					\node[point] (b2) at (0,2) {};
					\node[point] (c) at (2,1) [label=right:$x$] {};
					\node[point] (f1) at (3,2) {};
					\node[point] (f2) at (3,3) {};
					\node[point] (dc1) at (1.5,2) {};
					\node[point] (dt) at (1,3) {};
					\graph[use existing nodes] {
						a -- b1 -- b2 -- dt;
						a -- c -- dc1 -- dt;
						c -- f1 -- f2;
					};
					\node at (1,-1) {$F'$};
				\end{scope}
			\end{tikzpicture}
		\end{equation*}
		\caption{The technique in the proof of \cref{thm:gradification with Scott} does not work in general}
		\label{fig:gradification with Scott counterexample}
	\end{figure}

	Unfortunately, the proof of \cref{thm:gradification with Scott} crucially relied on the fact that the original frame $F$ was $(2 \cdot 1)$-connected. Consider for instance the frame $F$ given in \cref{fig:gradification with Scott counterexample}, which at $x$ is not $(2 \cdot 1)$-connected. If we apply the construction to $F$, we end up with a frame $F'$ in which $x$ sits below two connected components of height $1$, that is\footnote{Recall that $\concomps(F)$ is the set of connected components of $F$ and that $\contype(F)$ of $F$ is the signature $n_1^{m_1} \cdots n_k^{m_k}$ such that $\concomps(F)$ contains for each $i$ exactly $m_i$ sets of height $n_i-1$, and nothing else.}, $\contype(\Us x^{F'}) = 2^2$. Hence $F'$ is not $2^2$-connected, while $F$ is. Taking $2\cdot1$ away from $\Lam$ is a double-edged sword however, since it allows for more complex constructions in $F'$.

	The following reusable lemma will come in handy a couple of times. 

	\begin{lemma}\label{lem:concomps same condition}
		Let $f \colon F' \to F$ be a surjective p-morphism between finite posets, and take $x \in F'$. Assume that for any $y,z \in \Succ(x)$ there is a path $y \rsa z$ in $\Us x$ whenever there is a path $f(y) \rsa f(z)$ in $\Us{f(x)}$. Then:
		\begin{equation*}
			\concomps(\Us x) = \{\inv f[C] \mid C \in \concomps(\Us{f(x)})\}
		\end{equation*}
		In particular, if $\height(\inv f[C]) = \height(C)$ for any $C \in \concomps(\Us{f(x)})$ then:
		\begin{equation*}
			\contype(\Us x) = \contype(\Us{f(x)}
		\end{equation*}
	\end{lemma}

	\begin{proof}
		Note that, since $f$ is a p-morphism and $F$ and $F'$ are finite, $\{\inv f[C] \mid C \in \concomps(\Us{f(x)})\}$ is a partition of $\Us x$ into upwards- and downwards-closed sets. So it suffices to show that $\inv f[C]$ is connected for every $C \in \concomps(\Us{f(x)})$. Take $y_0, z_0 \in \inv f [C]$. Since $\inv f[C]$ is downwards-closed in $\Us x$, there are $y, z \in \Succ(x) \cap \inv f[C]$ such that $y \leq y_0$ and $z \leq z_0$. Then $f(y), f(z) \in C$, so by assumption there is a path $f(y) \rsa f(z)$ in $\Us {f(x)}$. But then by assumption there is a path $y \rsa z$ in $\Us x$, which lies in $\inv f[C]$ since the latter is upwards- and downwards-closed. This yields a path $y_0 \rsa z_0$.
	\end{proof}

\begin{figure}
		\begin{equation*}
			\begin{tikzpicture}[scale=0.6]
				\begin{scope}
					\graph[poset=0.6] {
						a[x=1] -- {
							b3 -- c3 -- d3 -- e3[x=1],
							b1[x=1,y=1] -- {
								e3,
								e2[y=1,x=2] -- n[y=1,x=2]
							},
						};
					};
					\node[below=0.5 of a] {$F$};
				\end{scope}
				\begin{scope}[xshift=150]
					\graph[poset=0.6] {
						a[x=1] -- {
							b1 -- c1 -- d1 -- e1,
							b2[x=1] -- {
								c2,
								c3[x=1] -- n[x=1]
							},
						};
					};
					\node[below=0.5 of a] {$\Tree(F)$};
				\end{scope}
				\begin{scope}[xshift=300]
					\graph[poset=0.6] {
						a[x=2.5] -- {
							xb7 -- xc7 -- xd7,
							xb6[x=-0.5] -- xc6[x=-0.5] -- xd6[x=-0.5] -- {
								xe7[x=-1] -- xd7, xe5[x=-1]
							},
							xb3[x=-1] -- xc3[x=-1] -- {xd4[x=-1.5] -- xe5, xd2[x=-1.5]},
							xb1[x=-1.5] -- {xc1[x=-2] -- xd2, xa2[x=-2]},
							x1[x=-2] -- {xa2, x2[x=-1] -- x3[x=-1]}
						};
					};
				\end{scope}
				\node[below=0.5 of a] {$F'$};
			\end{tikzpicture}
		\end{equation*}
		\caption{An example of gradification in the absence of Scott’s tree.}
		\label{fig:gradification without Scotts tree example}
	\end{figure}
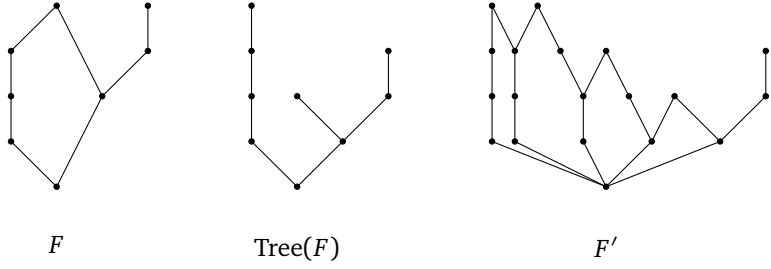

	\blue{Let us turn now to the proof of \cref{thm:gradification without Scott}.} The construction works in two steps as follows (see \cref{fig:gradification without Scotts tree example} for an example).
	\begin{itemize}
		\item Again, we start by splitting $F$ up into its tree unravelling $\Tree(F)$.
		\item Then, in order to connect the frame back up again while ensuring that it remains graded, we construct `zigzag roller-coasters' connecting top nodes of different heights.
	\end{itemize}

	\begin{proof}[Proof of \cref{thm:gradification without Scott}]
		As in the proof of \cref{thm:gradification with Scott}, we may assume that $\ep,1, 2 \notin \Lam$.

		Start with $T= \Tree(F)$. For every two distinct $p,q \in \Top(T)$ such that $\last(p)=\last(q)=t$, we will build a `roller-coaster' structure $Z(p,q)$, which will furnish a bridge between $p$ and $q$. Every such structure $Z(p,q)$ is independent, so that they can all be added to $T$ at the same time. First note that by \cref{cor:trees and nerves graded}, $T$ is graded; let $\rho \colon T \to \NN$ be its rank function.

		Now, take distinct $p,q \in \Top(T)$ such that $\last(p)=\last(q)=t$. Let $l \defeq \rho(q)-\rho(p)$. By swapping $p$ and $q$, we may assume that $l \geq 0$. We need to join $p$ to $q$ with a path which ascends in grade. We do this using a zigzagging path, which consists of lower points $a_0, \ldots, a_l$, upper points $b_0, \ldots, b_{l-1}$ and intermediate points $c_0, \ldots, c_{l-1}$. The relations between these points are as follows (see \cref{fig:zigzag relations; proof:gradification without Scott}).
		\begin{gather*}
			a_i < c_i < b_i, \qquad a_{i+1}<b_i
		\end{gather*}

		\begin{figure}
			\begin{equation*}
				\begin{tikzpicture}
					\node[point] (a0) [label=left:$a_0$] {};
					\node[point] (c0) [label=left:$c_0$,above right=5ex of a0] {};
					\node[point] (b0) [label=left:$b_0$,above right=5ex of c0] {};
					\node[point] (a1) [label=left:$a_1$,below right=5ex of b0] {};
					\node[point] (c1) [label=left:$c_1$,above right=5ex of a1] {};
					\node[point] (b1) [label=left:$b_1$,above right=5ex of c1] {};
					\node[point] (a2) [label=left:$a_2$,below right=5ex of b1] {};
					\node[point] (c2) [label=left:$c_2$,above right=5ex of a2] {};
					\node[point] (b2) [label=left:$b_2$,above right=5ex of c2] {};
					\node[point] (a3) [label=left:$a_3$,below right=5ex of b2] {};
					\path[draw] (a0) -- (c0) -- (b0) -- (a1) -- (c1) -- (b1) -- (a2) -- (c2) -- (b2) -- (a3);
				\end{tikzpicture}
			\end{equation*}
			\caption{The relations between the zigzag points in case $l=3$.}
			\label{fig:zigzag relations; proof:gradification without Scott}
		\end{figure}
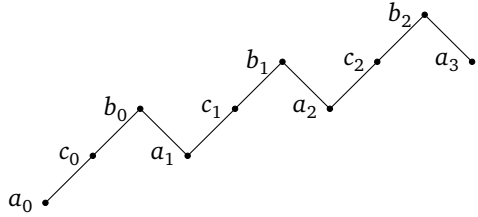

		Consider $p \wedge q$ (i.e. the intersection of $p$ and $q$, regarded as chains), and let $k \defeq \rho(p)-\rho(p \wedge q)-1$. Note that $k \geq 0$ since $p$ and $q$ are incomparable. Moreover, $k \geq 1$ as follows. Suppose for a contradiction that $k = 0$, so that $p$ is an immediate successor of $p \wedge q$. Then $\last(p)$ is an immediate successor of $\last(p \wedge q)$. But $\last(q)=\last(p)$, so we have, as chains:
		\begin{equation*}
			p = (p \wedge q) \cup \{\last(p)\} = (p \wedge q) \cup \{\last(q)\} = q
		\end{equation*}
		contradicting that $p$ and $q$ are distinct. \contradiction

		To ensure that the new poset $F'$ is still graded, we need to dangle some scaffolding down from the zigzag path to $p \wedge q$. Below each lower point $a_i$ we will dangle a chain of $k+i-1$ points $d(i,1), \ldots, d(i,{k+i-1})$. The relations are as follows:
		\begin{equation*}
			d(i,1) < d(i,2) < \cdots < d(i,{k+i-1}) < a_i
		\end{equation*}

		Finally, let $\mathrm Z(p,q)$ denote the whole structure of the zigzag path plus the dangling scaffolding. Attach $\mathrm Z(p,q)$ to $T$ by adding the following relations and closing under transitivity (see \cref{fig:zigzag plus ladders in context; proof:gradification without Scott}).
		\begin{equation*}
			a_0 < p, \qquad a_l < q, \qquad \forall i \colon p \wedge q < d(i,1)
		\end{equation*}

		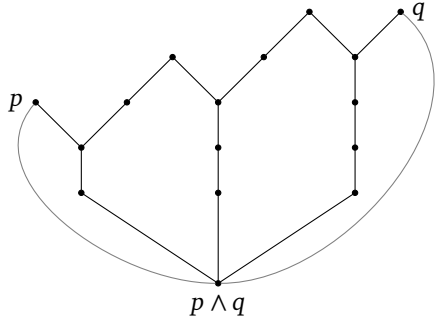
\begin{figure}
			\begin{equation*}
				\begin{tikzpicture}[scale=0.6]
					\node[point] (p) at (0,3) [label=left: $p$] {};
					\node[point] (a0) at (1,2) {};
					\node[point] (c0) at (2,3) {};
					\node[point] (b0) at (3,4) {};
					\node[point] (a1) at (4,3) {};
					\node[point] (c1) at (5,4) {};
					\node[point] (b1) at (6,5) {};
					\node[point] (a2) at (7,4) {};
					\node[point] (q) at (8,5) [label=right:$q$] {};

					\node[point] (d01) at (1,1) {};
					\node[point] (d12) at (4,2) {};
					\node[point] (d11) at (4,1) {};
					\node[point] (d23) at (7,3) {};
					\node[point] (d22) at (7,2) {};
					\node[point] (d21) at (7,1) {};

					\node[point] (r) at (4,-1) [label=below:$p \wedge q$] {};


					\graph[use existing nodes] {
						p -- a0 -- c0 -- b0 -- a1 -- c1 -- b1 -- a2 -- q;
						r -- d01;
						r -- d11;
						r -- d21;
						a0 -- d01;
						a1 -- d12 -- d11;
						a2 -- d23 -- d22 -- d21;
					};
					\draw[gray, out=180, in=230] (r) edge (p);
					\draw[gray, out=0, in=320] (r) edge (q);
				\end{tikzpicture}
			\end{equation*}
			\caption{The zigzag path and the ladder structure in place.}
			\label{fig:zigzag plus ladders in context; proof:gradification without Scott}
		\end{figure}

		Let $F'$ be the result of adding $\mathrm Z(p,q)$ to $T$ for every pair $p,q$, and define the function $f \colon F' \to F$ by:
		\begin{equation*}
			f(x) \defeq \left\{
			\begin{array}{ll}
				\last(x) & \text{if }x \in T \\
				\last(p) & \text{if }x \in Z(p,q)\text{ for some }p,q
			\end{array}
			\right.
		\end{equation*}

		First, let us see that $f$ is a p-morphism. The \eqref{p:p-morphism forth} condition follows from the fact that $\last$ is monotonic, and that:
		\begin{itemize}
			\item if $x \leq y$ with $x \in T$ and $y \in Z(p,q)$, then by construction $x \leq p \wedge q$, meaning that $f(x) = \last(x) \leq \last(p \wedge q) \leq \last(p) = f(y)$, and
			\item if $x \leq y$ with $x \in Z(p,q)$ and $y \in T$, then by construction $y \in \{p,q\}$, so that $f(x) = \last(p) = f(y)$.
		\end{itemize}
		The \eqref{p:p-morphism back} condition follows from the fact that $\last$ is open, and that each $Z(p,q)$ maps to a top node.

		Second, for any pair $p,q$, we can extend the rank function $\rho$ to the new structure $\mathrm Z(p,q)$ as follows (as indicated by the heights of the nodes in \cref{fig:zigzag plus ladders in context; proof:gradification without Scott}):
		\begin{gather*}
			\rho(a_i) = \rho(p) + i - 1 \\
			\rho(b_i) = \rho(p) + i + 1 \\
			\rho(c_i) = \rho(p) + i \\
			\rho(d(i,j)) = \rho(p \wedge q) + j 
		\end{gather*}
		To see that, thus extended, $\rho$ is still a rank function, it suffices to check that the newly-ranked $Z(p,q)$ fits into $T$ as a ranked structure. That is, we need to check the following equations.
		\begin{gather*}
			\rho(p) = \rho(a_0) + 1 \\
			\rho(q) = \rho(a_l) + 1 \\
			\rho(d(i,1)) = \rho(p \wedge q) + 1
		\end{gather*}
		But these follow by definition. In this way we see that $F'$ is graded.

		Finally, it remains to be shown that $F \vD \SFL(\Lam)$. So take $x \in F$. First, whenever $x \in Z(p,q)$ for some $p,q$, by construction $\Us x$ is $\alpha$-connected for every signature other than $\ep$, $1^2$, $2 \cdot 1$ and $k$ where $k \geq \height(F)+1$. Hence we may assume that $x \in T$. Let us use \cref{lem:concomps same condition}. Take $y,z \in \Succ(x)$ such that there is a path $f(y) \rsa f(z)$ in $\Us{\last(x)}$, with the aim of finding a path $y \rsa z$ in $\Us x$.

		First assume that $y \in Z(p,q)$ for some $p,q$. Then since $y \in \Succ(x)$ and $x \in T$, by construction $x = p \wedge q$. All of $Z(p,q)$ is connected in $\Us x$, hence there is a path $y \rsa p$. Let $p' \in T$ be the immediate successor of $x$ which lies below $p$ (this exists since $T$ is a tree). Then we have a path $y \rsa p'$ in $\Us x$. With this case thus dealt with, we may now assume that $y \in T$, and similarly that $z \in T$.

		So, we have a path $\last(y) \rsa \last(z)$. We now proceed in a similar fashion to the proof of \cref{thm:gradification with Scott}. We may assume that the path $\last(y) \rsa \last(z)$ has the form in \cref{fig:form of paths in F and F prime; proof:gradification without Scott} (a), where $t_0, \ldots, t_k$ are top nodes in $F$. Let $u_0 \defeq y$ and $u_k \defeq z$. For each $i \in \{1, \ldots, k-1\}$, choose $u_i \in \inv\last\{a_i\}$. For $i \in \{0, \ldots, k-1\}$, take $p_i,q_i \in \inv\last\{t_i\}$ such that $u_i \leq p_i$ and $u_{i+1} \leq q_i$. For each such $i$, since $\last(p_i)=\last(q_i)$, there is a path $p_i \rsa q_i$ which lies in $Z(p_i,q_i)$, and hence lies in $\Us x$. Compose all these paths as in \cref{fig:form of paths in F and F prime; proof:gradification without Scott} to form a path $y \rsa z$ in $\Us x$ as required.

		\begin{figure}
			\begin{equation*}
				\begin{tikzpicture}[scale=0.75]
					\begin{scope}
						\node at (-1.4,0.5) {(a)};
						\begin{scope}[xshift=30]
							\graph[grow right=0.75, math nodes] {
								a_0
									-- t_0[y=1cm]
									-- a_1
									-- t_1[y=1cm]
									-!- e0[x=0.5cm,y=0.25cm,as={$\cdots$}]
									-!- e1[empty nodes]
									-!- t_{k-2}[y=1cm]
									-- a_{k-1}
									-- t_{k-1}[y=1cm]
									-- a_k
									;
							};
							\draw (t_1) -- ++(0.5,-0.5);
							\draw (t_{k-2}) -- ++(-0.5,-0.5);
						\end{scope}
					\end{scope}
					\begin{scope}[yshift=-130]
						\node at (-1.4,1.4) {(b)};
						\begin{scope}[xscale=1.3,yscale=1.8]
							\draw (0,0) node[point] [label=below:$u_0$] {}
								-- (0.5,1) node[point] [label=above:$p_0$] {}
								.. controls (0.75,1.2) and (1.25,0.8) .. (1.5,1) node[point] [label=above:$q_0$] {}
								-- (2,0) node[point] [label=below:$u_1$] {}
								-- (2.5,1) node[point] [label=above:$p_1$] {}
								.. controls (2.75,1.2) and (3.25,0.8) .. (3.5,1) node[point] [label=above:$q_1$] {}
								-- (3.75,0.5);
							\draw[dashed] (1,1.2) ellipse (0.9 and 0.4);
							\node at (1,1.75) {$t_0$};
							\draw[dashed] (3,1.2) ellipse (0.9 and 0.4);
							\node at (3,1.75) {$t_1$};
							\node at (4.25,0.5) {$\cdots$};
							\draw (4.75,0.5)
								-- (5,1) node[point] [label=above:$p_{k-2}$] {}
								.. controls (5.25,1.2) and (5.75,0.8) .. (6,1) node[point] [label=above:$q_{k-2}$] {}
								-- (6.5,0) node[point] [label=below:$u_{k-1}$] {}
								-- (7,1) node[point] [label=above:$p_{k-1}$] {}
								.. controls (7.25,1.2) and (7.75,0.8) .. (8,1) node[point] [label=above:$q_{k-1}$] {}
							-- (8.5,0) node[point] [label=below:$u_k$] {};
							\draw[dashed] (5.5,1.2) ellipse (0.9 and 0.4);
							\node at (5.5,1.75) {$t_{k-2}$};
							\draw[dashed] (7.5,1.2) ellipse (0.9 and 0.4);
							\node at (7.5,1.75) {$t_{k-1}$};
						\end{scope}
					\end{scope}
				\end{tikzpicture}
			\end{equation*}
			\caption{The form of the paths in $\protect\Us{\last(x)}$ and $\protect\Us x$}
			\label{fig:form of paths in F and F prime; proof:gradification without Scott}
		\end{figure}
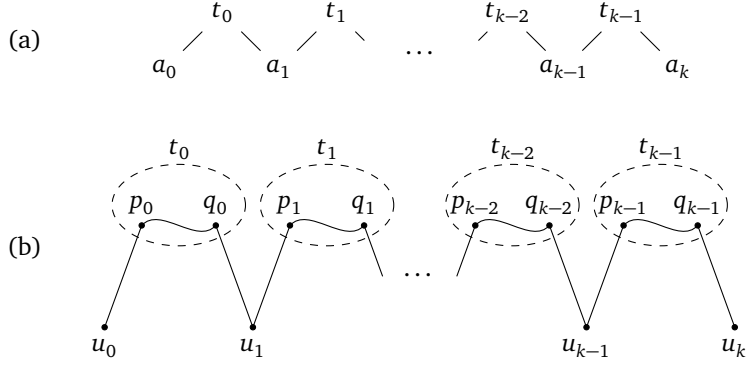

		It now remains to show that if $C \in \concomps(\Us{\last(x)})$, then $\height(\inv f[C]) = \height(C)$. First, since $f$ is a p-morphism, $\height(\inv f[C]) \geq \height(C)$. Conversely, let $X \sse \inv f[C]$ be a maximal chain. Assume $X$ intersects with some $Z(p,q)$. Then we can replace the part $X \cap (Z(p,q) \cup \{p,q\})$ with the unique maximal chain in $\Us{p \wedge q}^T$ containing $q$ (this exists since $T$ is a tree). Then by construction this does not decrease the length of $X$ nor does it move $X$ outside of $\inv f[C]$ (since the latter is upwards- and downwards-closed). Therefore, we may assume that $X \sse T$, so $X$ corresponds to a chain $\last[X]$ of the same length in $C$.

		Therefore, by \cref{lem:concomps same condition} we get that $\contype(\Us x) = \contype(\Us{\last(x)}$. Applying \cref{lem:alpha partition iff alpha leq contype}, we have that $\Us x$ has an $\alpha$-partition if and only if $\Us{\last(x)}$ has an $\alpha$-partition. 
	\end{proof}

\subsection{Nervification}

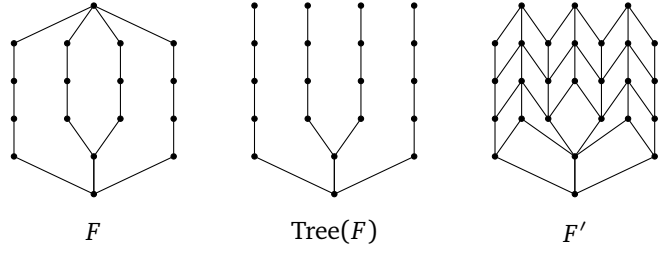
\begin{figure}
		\begin{equation*}
			\begin{tikzpicture}[scale=0.5]
				\begin{scope}
					\node[point] (a1) at (0,-2.5) {};
					\node[point] (b1) at (-2.1,-1.5) {};
					\node[point] (b2) at (-2.1,-0.5) {};
					\node[point] (b3) at (-2.1,0.5) {};
					\node[point] (b4) at (-2.1,1.5) {};
					\node[point] (c1) at (0,-1.5) {};
					\node[point] (d1) at (-0.7,-0.5) {};
					\node[point] (d2) at (-0.7,0.5) {};
					\node[point] (d3) at (-0.7,1.5) {};
					\node[point] (e1) at (0.7,-0.5) {};
					\node[point] (e2) at (0.7,0.5) {};
					\node[point] (e3) at (0.7,1.5) {};
					\node[point] (g1) at (2.1,-1.5) {};
					\node[point] (g2) at (2.1,-0.5) {};
					\node[point] (g3) at (2.1,0.5) {};
					\node[point] (g4) at (2.1,1.5) {};
					\node[point] (x) at (0,2.5) {};
					\graph[use existing nodes] {
						a1 -- b1 -- b2 -- b3 -- b4;
						a1 -- c1 -- d1 -- d2 -- d3;
						a1 -- c1 -- e1 -- e2 -- e3;
						a1 -- g1 -- g2 -- g3 -- g4;
						x -- {b4, d3, e3, g4};
					};
					\node at (0,-3.5) {$F$};
				\end{scope}
				\begin{scope}[xshift=180]
					\node[point] (a1) at (0,-2.5) {};
					\node[point] (b1) at (-2.1,-1.5) {};
					\node[point] (b2) at (-2.1,-0.5) {};
					\node[point] (b3) at (-2.1,0.5) {};
					\node[point] (b4) at (-2.1,1.5) {};
					\node[point] (c1) at (0,-1.5) {};
					\node[point] (d1) at (-0.7,-0.5) {};
					\node[point] (d2) at (-0.7,0.5) {};
					\node[point] (d3) at (-0.7,1.5) {};
					\node[point] (e1) at (0.7,-0.5) {};
					\node[point] (e2) at (0.7,0.5) {};
					\node[point] (e3) at (0.7,1.5) {};
					\node[point] (g1) at (2.1,-1.5) {};
					\node[point] (g2) at (2.1,-0.5) {};
					\node[point] (g3) at (2.1,0.5) {};
					\node[point] (g4) at (2.1,1.5) {};
					\node[point] (x1) at (-2.1,2.5) {};
					\node[point] (x2) at (-0.7,2.5) {};
					\node[point] (x3) at (0.7,2.5) {};
					\node[point] (x4) at (2.1,2.5) {};
					\graph[use existing nodes] {
						a1 -- b1 -- b2 -- b3 -- b4 -- x1;
						a1 -- c1 -- d1 -- d2 -- d3 -- x2;
						a1 -- c1 -- e1 -- e2 -- e3 -- x3;
						a1 -- g1 -- g2 -- g3 -- g4 -- x4;
					};
					\node at (0,-3.5) {$\Tree(F)$};
				\end{scope}
				\begin{scope}[xshift=360]
					\node[point] (a1) at (0,-2.5) {};
					\node[point] (b1) at (-2.1,-1.5) {};
					\node[point] (b2) at (-2.1,-0.5) {};
					\node[point] (b3) at (-2.1,0.5) {};
					\node[point] (b4) at (-2.1,1.5) {};
					\node[point] (c1) at (0,-1.5) {};
					\node[point] (d1) at (-0.7,-0.5) {};
					\node[point] (d2) at (-0.7,0.5) {};
					\node[point] (d3) at (-0.7,1.5) {};
					\node[point] (e1) at (0.7,-0.5) {};
					\node[point] (e2) at (0.7,0.5) {};
					\node[point] (e3) at (0.7,1.5) {};
					\node[point] (g1) at (2.1,-1.5) {};
					\node[point] (g2) at (2.1,-0.5) {};
					\node[point] (g3) at (2.1,0.5) {};
					\node[point] (g4) at (2.1,1.5) {};
					\node[point] (bd1) at (-1.4,-0.5) {};
					\node[point] (bd2) at (-1.4,0.5) {};
					\node[point] (bd3) at (-1.4,1.5) {};
					\node[point] (bd4) at (-1.4,2.5) {};
					\node[point] (de1) at (0,0.5) {};
					\node[point] (de2) at (0,1.5) {};
					\node[point] (de3) at (0,2.5) {};
					\node[point] (eg1) at (1.4,-0.5) {};
					\node[point] (eg2) at (1.4,0.5) {};
					\node[point] (eg3) at (1.4,1.5) {};
					\node[point] (eg4) at (1.4,2.5) {};
					\graph[use existing nodes] {
						a1 -- b1 -- b2 -- b3 -- b4;
						a1 -- c1 -- d1 -- d2 -- d3;
						a1 -- c1 -- e1 -- e2 -- e3;
						a1 -- g1 -- g2 -- g3 -- g4;
						bd1 -- bd2 -- bd3 -- bd4;
						de1 -- de2 -- de3;
						eg1 -- eg2 -- eg3 -- eg4;
						b1 -- bd1 -- c1;
						b2 -- bd2 -- d1;
						b3 -- bd3 -- d2;
						b4 -- bd4 -- d3;
						d1 -- de1 -- e1;
						d2 -- de2 -- e2;
						d3 -- de3 -- e3;
						c1 -- eg1 -- g1;
						e1 -- eg2 -- g2;
						e2 -- eg3 -- g3;
						e3 -- eg4 -- g4;
					};
					\node at (0,-3.5) {$F'$};
				\end{scope}
			\end{tikzpicture}
		\end{equation*}
		\caption{An example of nervification, using the graded structure of $F$}
		\label{fig:nervification example}
	\end{figure}

	We now find ourselves, having suitably prepared $F$, in a position to make use of its additional graded structure. The general method of the final construction, in which we transform $F$ into a frame which nerve-validates $\SFL(\Lam)$, is the same as in \cref{thm:gradification with Scott} and \cref{thm:gradification without Scott}. We begin with the tree unravelling $\Tree(F)$, perform some alterations, then rejoin top nodes. A key difference here is that we won't rejoin every top node to every other top node whose `$\last$' value is the same. Instead, we line up all the top nodes mapping to the same element and link each top node to at most two other top nodes, which we think of as its neighbours. See \cref{fig:nervification example} for an example of the construction.

	\begin{definition}
		Let $T$ be a finite tree. Then for each $x \in T$, we have that $\ds x$ is a chain. For $k\leq \height(x)$, let $x^{(k)}$ be the element of this chain which has height $k$. Let $x^{(-k)}$ be the element which has height $\height(x)-k$.
	\end{definition}

	\begin{definition}
		For $n \in \NN$, let $\Starlikes_n \defeq \Starlikes \setminus \{1^k \mid k < n\}$.
	\end{definition}

	\begin{theorem}\label{thm:nervification}
		Take $\Lam \sse \Starlikes$ and let $F$ be a finite, graded, rooted poset of height $n$ such that $F \vD \SFL(\Lam)$. Then there is a poset $F'$ and a p-morphism $f \colon F' \to F$ such that $F' \vD \SFL(\Lam)$ and such that $F'$ is $\alpha$-diamond-connected for every $\alpha \in \Starlikes_n$.
	\end{theorem}

	\begin{proof}[Proof of \cref{thm:nervification}]
		We may assume that $\ep,1 \notin \Lam$. Further, if $2 \in \Lam$, then $\height(F)=1$, so $F$ is already $\alpha$-diamond-connected for every $\alpha \in \Starlikes_n$. Hence we may assume that $2 \notin \Lam$.

		Once more, start with $T = \Tree(F)$. Chop off the top nodes: let $T' \defeq \Trunk(T)$. For each $t \in \Top(F)$, we will add a new structure $W(t)$, which lies only above elements of $T'$. Let $\rho \colon F \to \NN$ be the rank function on $F$. Note that $\rho \circ \last \colon T \to \NN$ is the rank function on $T$. 

		Take $t \in \Top(F)$. Enumerate $\last^{-1}\{t\}= \{p_1, \ldots, p_m\}$. For each $i \leq m-1$, define:
		\begin{gather*}
			r_i \defeq p_i \wedge p_{i+1} \\
			l_i \defeq \rho(\last(r_i)) \\
			k_i \defeq \rho(t)-\rho(\last(r_i))-1
		\end{gather*}
		Note that $k_i \geq 1$ just as in the proof of \cref{thm:gradification without Scott}. Since $F$ is graded and $T$ is a tree, we have that: 
		\begin{equation*}
			\abs{\Uds{r_i,p_i}^T} = \abs{\Uds{r_i,{p_{i+1}}}^T} = k_i
		\end{equation*}
		In other words, $p_i^{(l_i)} = p_{i+1}^{(l_i)} = r_i$. We will construct a `chevron' structure which joins $p_i^{(-1)}$ to $p_{i+1}^{(-1)}$. For each $i \leq m-1$, take new elements $a(i,1), \ldots, a(i,{k_i})$, and add them to $T'$ using the following relations.
		\begin{equation*}
			a(i,1) < \cdots < a(i,{k_i}), \qquad \forall j \leq k_i \colon p_i^{(l+j)}, p_{i+1}^{(l+j)} < a(i,j)
		\end{equation*}
		Let $W(t)$ be this new structure (i.e. the chain $\{a(i,1) < \cdots < a(i,{k_i})\}$ in place). See \cref{fig:chevron structure two branches; proof:nervification} and \cref{fig:chevron structure three branches; proof:nervification} for examples of this process of adding chevrons. 

		\begin{figure}
			\begin{equation*}
				\begin{tikzpicture}[scale=0.7]
					\begin{scope}
						\node[point] (a) at (1,0) {};
						\node[point] (b1) at (0,1) {};
						\node[point] (b2) at (0,2) {};
						\node[point] (b3) at (0,3) {};
						\node[point] (c1) at (2,1) {};
						\node[point] (c2) at (2,2) {};
						\node[point] (c3) at (2,3) {};
						\node[point] (d) at (1,4) [label=right:$t$] {};
						\graph[use existing nodes] {
							a -- b1 -- b2 -- b3 -- d;
							a -- c1 -- c2 -- c3 -- d;
						};
						\node at (1,-1) {$F$};
					\end{scope}
					\begin{scope}[xshift=125]
						\node[point] (a) at (1,0) [label=right:$r_1$] {};
						\node[point] (b1) at (0,1) {};
						\node[point] (b2) at (0,2) {};
						\node[point] (b3) at (0,3) [label=left:$p_1^{(-1)}$] {};
						\node[point] (b4) at (0.5,4) [label=left:$p_1$] {};
						\node[point] (c1) at (2,1) {};
						\node[point] (c2) at (2,2) {};
						\node[point] (c3) at (2,3) [label=right:$p_2^{(-1)}$] {};
						\node[point] (c4) at (1.5,4) [label=right:$p_2$] {};
						\draw[dashed] (1,4) ellipse (1.4 and 0.7);
						\node at (2.7,4) {$t$};
						\graph[use existing nodes] {
							a -- b1 -- b2 -- b3 -- b4;
							a -- c1 -- c2 -- c3 -- c4;
						};
						\node at (1,-1) {$T$};
					\end{scope}
					\begin{scope}[xshift=275]
						\node[point] (a) at (1,0) [label=right:$r_1$] {};
						\node[point] (b1) at (0,1) [label=left:$p_1^{(1)}$] {};
						\node[point] (b2) at (0,2) [label=left:$p_1^{(2)}$] {};
						\node[point] (b3) at (0,3) [label=left:$p_1^{(3)}$] {};
						\node[point] (c1) at (2,1) [label=right:$p_{2}^{(1)}$] {};
						\node[point] (c2) at (2,2) [label=right:$p_{2}^{(2)}$] {};
						\node[point] (c3) at (2,3) [label=right:$p_{2}^{(3)}$] {};
						\node[point] (e1) at (1,2) {};
						\node[point] (e2) at (1,3) {};
						\node[point] (e3) at (1,4) {};
						\graph[use existing nodes] {
							a -- b1 -- b2 -- b3;
							a -- c1 -- c2 -- c3;
							e1 -- e2 -- e3;
							b1 -- e1 -- c1;
							b2 -- e2 -- c2;
							b3 -- e3 -- c3;
						};
						\node[above left=0.5cm of e3] (la13) {$a(1,3)$};
						\node[right=0.6cm of la13] (la12) {$a(1,2)$};
						\node[right=0.1cm of la12] (la11) {$a(1,1)$};
						\draw[-Latex] (la13) edge (e3);
						\draw[-Latex] (la12) edge (e2);
						\draw[-Latex,bend right=10] (la11) edge (e1);
						\node at (1,-1) {$T'+W(t)$};
					\end{scope}
				\end{tikzpicture}
			\end{equation*}
			\caption{The chevron structure in a case with two branches.}
			\label{fig:chevron structure two branches; proof:nervification}
		\end{figure}
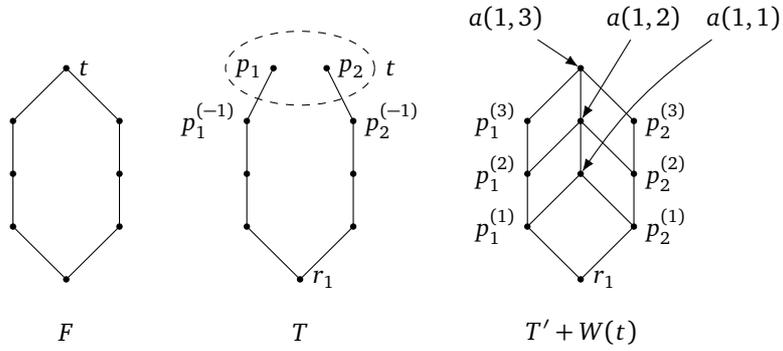

		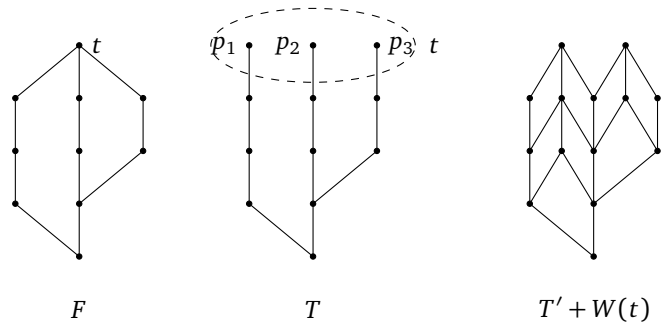
\begin{figure}
			\begin{equation*}
				\begin{tikzpicture}[scale=0.7]
					\begin{scope}[,xscale=1.2]
						\node[point] (a) at (1,0) {};
						\node[point] (b1) at (0,1) {};
						\node[point] (b2) at (0,2) {};
						\node[point] (b3) at (0,3) {};
						\node[point] (c1) at (1,1) {};
						\node[point] (c2) at (1,2) {};
						\node[point] (c3) at (1,3) {};
						\node[point] (e2) at (2,2) {};
						\node[point] (e3) at (2,3) {};
						\node[point] (d) at (1,4) [label=right:$t$] {};
						\graph[use existing nodes] {
							a -- b1 -- b2 -- b3 -- d;
							a -- c1 -- c2 -- c3 -- d;
							c1 -- e2 -- e3 -- d;
						};
						\node at (1,-1) {$F$};
					\end{scope}
					\begin{scope}[xshift=125,xscale=1.2]
						\node[point] (a) at (1,0) {};
						\node[point] (b1) at (0,1) {};
						\node[point] (b2) at (0,2) {};
						\node[point] (b3) at (0,3) {};
						\node[point] (b4) at (0,4) [label=left:$p_1$] {};
						\node[point] (c1) at (1,1) {};
						\node[point] (c2) at (1,2) {};
						\node[point] (c3) at (1,3) {};
						\node[point] (c4) at (1,4) [label=left:$p_2$] {};
						\node[point] (e2) at (2,2) {};
						\node[point] (e3) at (2,3) {};
						\node[point] (e4) at (2,4) [label=right:$p_3$] {};
						\draw[dashed] (1,4) ellipse (1.6 and 0.7);
						\node at (2.9,4) {$t$};
						\graph[use existing nodes] {
							a -- b1 -- b2 -- b3 -- b4;
							a -- c1 -- c2 -- c3 -- c4;
							c1 -- e2 -- e3 -- e4;
						};
						\node at (1,-1) {$T$};
					\end{scope}
					\begin{scope}[xshift=275,xscale=1.2]
						\node[point] (a) at (1,0) {};
						\node[point] (b1) at (0,1) {};
						\node[point] (b2) at (0,2) {};
						\node[point] (b3) at (0,3) {};
						\node[point] (c1) at (1,1) {};
						\node[point] (c2) at (1,2) {};
						\node[point] (c3) at (1,3) {};
						\node[point] (e2) at (2,2) {};
						\node[point] (e3) at (2,3) {};
						\node[point] (bc1) at (0.5,2) {};
						\node[point] (bc2) at (0.5,3) {};
						\node[point] (bc3) at (0.5,4) {};
						\node[point] (ce2) at (1.5,3) {};
						\node[point] (ce3) at (1.5,4) {};
						\graph[use existing nodes] {
							a -- b1 -- b2 -- b3;
							a -- c1 -- c2 -- c3;
							c1 -- e2 -- e3;
							bc1 -- bc2 -- bc3;
							ce2 -- ce3;
							b1 -- bc1 -- c1;
							b2 -- bc2 -- c2 -- ce2 -- e2;
							b3 -- bc3 -- c3 -- ce3 -- e3;
						};
						\node at (1,-1) {$T'+W(t)$};
					\end{scope}
				\end{tikzpicture}
			\end{equation*}
			\caption{The chevron structure in a more complex case involving three branches.}
			\label{fig:chevron structure three branches; proof:nervification}
		\end{figure}

		The process of adding $W(t)$ is independent for each $t \in \Top(F)$. Let $F'$ be the result of adding every $W(t)$ to $T'$. Define $f \colon F' \to F$ by:
		\begin{equation*}
			f(x) \defeq \left\{
			\begin{array}{ll}
				\last(x) & \text{if }x \in T' \\
				t & \text{if }x \in W(t)\text{ for some }t \in \Top(F)
			\end{array}
			\right.
		\end{equation*}
		Since we have made sure that each $W(t)$ contains, for each $p_i \in \last^{-1}\{t\}$, a node above $p_i^{(-1)}$ which maps to $t$, and that all of the new structure maps to a top node, $f$ is a p-morphism.

		Let us see that $F' \vD \SFL(\Lam)$. Take $x \in F'$. If $x \in W(t)$ for some $t$, then $\Us x$ is either empty or a chain, hence $\Us x \vD \SFL(\Lam)$. So we assume that $x \in T'$. The verification is now very similar to that in \cref{thm:gradification without Scott}, making use of \cref{lem:concomps same condition}. Take $y,z \in \Succ(x)$ such that there is a path $f(y) \rsa f(z)$ in $\Us{\last(x)}$. As in the proof of \cref{thm:gradification without Scott}, by construction of $W(t)$ we may assume that $y,z \in T'$. Just as in that proof, we can construct a path $y \rsa z$ from the path $f(y) \rsa f(z)$, using the fact that whenever $t \in \Us{\last(x)} \cap \Top(F)$, any $w,v \in \invfs t$ are connected by a path in $\Us x^{F'}$ (this is how we constructed $F'$). It is straightforward then to check that if $C \in \concomps(\Us{\last(x)})$ we have $\height(\inv f[C]) = \height(C)$, giving that:
		\begin{equation*}
			\contype(\Us x) = \contype(\Us{\last(x)})
		\end{equation*}

		To complete the proof, let us see that $F'$ is $\alpha$-diamond-connected for every $\alpha \in \Starlikes_n$. Take $x,y \in F'$ with $x < y$ and consider $\Uds{x,y}$. There are several cases.
		\begin{enumerate}[label=(\alph*)]
			\item Case $y \in T'$. We have that $\Uds{x,y}^{F'} = \Uds{x,y}^{T'}$, which is linearly-ordered since $T'$ is a tree; hence it is connected and of height at most $n-2$.
		\end{enumerate}
		Hence $y = a(i,j)$ for $a(i,j) \in W(t)$ a new element. Let $p_i,p_{i+1}, r_i, l_i$ be as above.
		\begin{enumerate}[resume*]
			\item Case $x \in W(t)$. Note that by construction $\Uds{x,y}$ is linearly-ordered.
			\item Case $x=p_i^{(l+e)}$ for some $e$. If we have $\height(\Uds{x,y})=1$, then $e=i-1$ and $\Uds{x,y}$ is the antichain on two elements, which is $\alpha$-connected. Otherwise, by construction, $a(i,{j-1}) \in \Uds{x,y}$ which is connected to everything.
			\item Case $x=p_{i+1}^{(l+e)}$ for some $e$. This is symmetric.
			\item Case $x=r_i$. Again, if $\height(\Uds{x,y}))=1$ then $j=1$ and $\Uds{x,y}$ is the antichain on two elements, otherwise $a(i,1) \in \Uds{x,y}$ which is connected to everything.
			\item Otherwise, $x < r_i$ (since $T'$ is a tree). Then $r_i \in \Uds{x,y}$ which is connected to everything. \qedhere
		\end{enumerate}
	\end{proof}

\subsection{End of Proof of \texorpdfstring{\cref{thm:starlike completeness}}{Theorem \ref{thm:starlike completeness}}}

	We can now prove our second main result:

	\begin{proof}[Proof of \cref{thm:starlike completeness}]
		By \cref{lem:poly-complete iff fmp and p-morphic of nerve-validated} and \cref{lem:starlike logics fmp}, we need to show that every finite, rooted frame of $\SFL(\Lam)$ is the up-reduction of one which nerve-validates $\SFL(\Lam)$; in fact this up-reduction is just a p-morphism. So take such a frame $F$. We may assume that $F$ is graded: when we have $2 \cdot 1 \in \Lam$, apply \cref{thm:gradification with Scott}, otherwise apply \cref{thm:gradification without Scott}. Then by \cref{thm:nervification}, there is a frame $F'$ and a p-morphism $f \colon F' \to F$ such that $F'$ is $\alpha$-nerve-connected for every $\alpha \in \Lam$ (note that by \cref{cor:BDn Jankov-Fine of starlike} we must have $\Lam \sse \Starlikes_n$ where $n = \height(F)$). Then, by \cref{thm:nerve-validates Jankov of starlike iff alpha-nerve-connected}, $F'$ nerve-validates $\SFL(\Lam)$, which completes the proof.
	\end{proof}

\section{Acknowledgements}
\label{sec:acknowledgements}

\begin{blueblock}
    We are grateful to the referee for many helpful and detailed comments which improved the presentation of the paper.

    The first author was supported by the Amsterdam Science Talent Scholarship during this research.
\end{blueblock}


	\printbibliography

\end{document}